\documentclass[12pt]{iopart}   
\usepackage[utf8]{inputenc}

\usepackage{amsthm}
\usepackage{amsfonts}
\usepackage{amssymb}
\expandafter\let\csname equation*\endcsname\relax
\expandafter\let\csname endequation*\endcsname\relax
\usepackage{amsmath}
\usepackage{mathtools}

\usepackage{graphicx}
\usepackage{macros}
\usepackage{bm}
\usepackage[usenames,dvipsnames]{xcolor}
\usepackage{tikz}
\usepackage{pgfplots} 
\usepackage{subcaption}
\usepackage{tikz}
\usetikzlibrary{arrows.meta, positioning, patterns}
\usetikzlibrary{ipe} 
\usepackage{float}
\usepackage{color}
\usepackage{hyperref}
\hypersetup{
  colorlinks,
  citecolor= Brown,
  linkcolor= link,
  urlcolor= link
}

\usepackage{appendix}
\usepackage{array}
\usepackage{footnote}
\usepackage{booktabs}
\usepackage{hhline}
\makesavenoteenv{tabular}

\usepackage{algpseudocode}
\usepackage{algorithm}

\usepackage{scrextend}
 \deffootnote{0em}{1.6em}{\enskip}

\tikzstyle{ipe stylesheet} = [
  ipe import,
  even odd rule,
  line join=round,
  line cap=butt,
  ipe pen normal/.style={line width=0.4},
  ipe pen heavier/.style={line width=0.8},
  ipe pen fat/.style={line width=1.2},
  ipe pen ultrafat/.style={line width=2},
  ipe pen normal,
  ipe mark normal/.style={ipe mark scale=3},
  ipe mark large/.style={ipe mark scale=5},
  ipe mark small/.style={ipe mark scale=2},
  ipe mark tiny/.style={ipe mark scale=1.1},
  ipe mark normal,
  /pgf/arrow keys/.cd,
  ipe arrow normal/.style={scale=7},
  ipe arrow large/.style={scale=10},
  ipe arrow small/.style={scale=5},
  ipe arrow tiny/.style={scale=3},
  ipe arrow normal,
  /tikz/.cd,
  ipe arrows, % update arrows
  <->/.tip = ipe normal,
  ipe dash normal/.style={dash pattern=},
  ipe dash dashed/.style={dash pattern=on 4bp off 4bp},
  ipe dash dotted/.style={dash pattern=on 1bp off 3bp},
  ipe dash dash dotted/.style={dash pattern=on 4bp off 2bp on 1bp off 2bp},
  ipe dash dash dot dotted/.style={dash pattern=on 4bp off 2bp on 1bp off 2bp on 1bp off 2bp},
  ipe dash normal,
  ipe node/.append style={font=\normalsize},
  ipe stretch normal/.style={ipe node stretch=1},
  ipe stretch normal,
  ipe opacity 10/.style={opacity=0.1},
  ipe opacity 30/.style={opacity=0.3},
  ipe opacity 50/.style={opacity=0.5},
  ipe opacity 75/.style={opacity=0.75},
  ipe opacity opaque/.style={opacity=1},
  ipe opacity opaque,
]

\definecolor{red}{rgb}{1,0,0}
\definecolor{green}{rgb}{0,1,0}
\definecolor{blue}{rgb}{0,0,1}
\definecolor{yellow}{rgb}{1,1,0}
\definecolor{orange}{rgb}{1,0.647,0}
\definecolor{gold}{rgb}{1,0.843,0}
\definecolor{purple}{rgb}{0.627,0.125,0.941}
\definecolor{gray}{rgb}{0.745,0.745,0.745}
\definecolor{brown}{rgb}{0.647,0.165,0.165}
\definecolor{navy}{rgb}{0,0,0.502}
\definecolor{pink}{rgb}{1,0.753,0.796}
\definecolor{seagreen}{rgb}{0.18,0.545,0.341}
\definecolor{turquoise}{rgb}{0.251,0.878,0.816}
\definecolor{violet}{rgb}{0.933,0.51,0.933}
\definecolor{darkblue}{rgb}{0,0,0.545}
\definecolor{darkcyan}{rgb}{0,0.545,0.545}
\definecolor{darkgray}{rgb}{0.663,0.663,0.663}
\definecolor{darkgreen}{rgb}{0,0.392,0}
\definecolor{darkmagenta}{rgb}{0.545,0,0.545}
\definecolor{darkorange}{rgb}{1,0.549,0}
\definecolor{darkred}{rgb}{0.545,0,0}
\definecolor{lightblue}{rgb}{0.678,0.847,0.902}
\definecolor{lightcyan}{rgb}{0.878,1,1}
\definecolor{lightgray}{rgb}{0.827,0.827,0.827}
\definecolor{lightgreen}{rgb}{0.565,0.933,0.565}
\definecolor{lightyellow}{rgb}{1,1,0.878}
\definecolor{black}{rgb}{0,0,0}
\definecolor{white}{rgb}{1,1,1}

\definecolor{link}{rgb}{0.18,0.25,0.63}
\definecolor{myred}{rgb}{0.7,0.25,0.2}
\usepackage[top=3cm,bottom=3cm,left=2.5cm,right=2.5cm]{geometry}

\numberwithin{equation}{section}

\DeclareMathOperator*{\esssup}{ess\,sup}
\DeclareMathOperator*{\essinf}{ess\,inf}

\DeclareMathOperator*{\argmin}{argmin}

\newcommand{\authorfootnotes}{\renewcommand\thefootnote{\@fnsymbol\c@footnote}}%
\makeatother
\graphicspath{{./images/}{./images/validation/}}

\begin{document}
\definecolor{link}{rgb}{0,0,0}
\definecolor{mygrey}{rgb}{0.34,0.34,0.34}
\def\blue #1{{\color{blue}#1}}

 \title[Analysis and Optimisation of TV-IC denoising model]{Analysis and optimisation of a variational model for mixed Gaussian and Salt \& Pepper noise removal}
 
 \author{Luca Calatroni$^1$, Kostas Papafitsoros$^2$}

\address{$^1$ CMAP, \'Ecole Polytechnique, Route de Saclay, 91128 Palaiseau Cedex, France}
\address{$^2$ Weierstrass Institute for Applied Analysis and Stochastics (WIAS), Mohrenstrasse 39, 10117, Berlin, Germany}
\ead{luca.calatroni@polytechnique.edu, kostas.papafitsoros@wias-berlin.de}

\begin{abstract}
We analyse a variational regularisation problem for mixed noise removal that was recently proposed in \cite{calatroni_mixed}. The data discrepancy term of the model combines $L^{1}$ and $L^{2}$ terms in an infimal convolution fashion and it is appropriate for the joint removal of Gaussian and Salt \& Pepper noise. In this work we perform a finer analysis of the model which emphasises on the balancing effect of the two parameters appearing in the discrepancy term. Namely, we study the asymptotic behaviour of the model for large and small values of these parameters and we compare it to the corresponding variational models with $L^{1}$ and $L^{2}$ data fidelity. Furthermore, we compute exact solutions for simple data functions taking the total variation as regulariser. Using these theoretical results, we then analytically study a bilevel optimisation strategy for automatically selecting the parameters of the model by means of a training set. Finally, we report some numerical results on the selection of the optimal noise model via such strategy which confirm the validity of our analysis and the use of popular data models in the case of ``blind" model selection. 
\end{abstract}

%\tableofcontents

\definecolor{link}{rgb}{0.18,0.25,0.63}

\section{Introduction}

Image denoising is a classical problem in imaging which is defined as the task of removing oscillations and interferences from a given image $f$.  Given a regular image domain $\Omega\subseteq\RR^2$ the problem can be formulated mathematically as the task of retrieving a noise-free image $u$ from its noisy measurement $f$, the latter being the result of a (possibly non-linear) degradation operator $\mathcal{T}$. In its general form, this problem can be written in the following way
$$
\text{find }u\quad \text{such that}\quad f=\mathcal{T}(u).
$$
Here the operator $\mathcal{T}$ introduces noise in the image, not only in an additive way, and it is not to be confused with the forward operator in the context of inverse problems.
In order to obtain a noise-free image, a classical technique consists of minimising an appropriate energy functional $\mathcal{J}$ over a suitable Banach space $X$ where the image functions are assumed to lie. In its general form the problem reads as follows:
\begin{equation}  \label{functional_energy}
\min_{u\in X}~\left\{ \mathcal{J}(u):= R(u) + \lambda\Phi(u,f) \right\}.
\end{equation}
Here, $R(u)$ stands for the regularisation term encoding \emph{a priori} information on the regularity of the solution, $\Phi(u,f)$ for the data-fitting measure that depends on the statistical and physical assumptions in the data and $\lambda>0$ is a scalar parameter whose magnitude balances the regularisation against trust in the data. Since the seminal work of Rudin, Osher, Fatemi \cite{rudin1992nonlinear}, a popular choice for $R$ in \eqref{functional_energy} is $R(u)=|Du|(\Omega)$, the Total Variation (TV) seminorm \cite{AmbrosioBV}, due to its ability of preserving salient structures in the image, i.e., edges, while removing noise at the same time. Here $Du$ represents the distributional derivative of the function $u\in \bv(\om)$, the space of functions of bounded variation, and $|Du|(\Omega)$ is the total variation of this measure. In the recent years, also higher-order regularisation terms that improve upon TV-induced artefacts -- notably, the creation of piecewise constant structures -- have been proposed in the literature. Among those, the most prominent is the Total Generalized Variation (TGV) \cite{TGV}, see also \cite{papafitsorosphd} for a comprehensive review.

% for instance,  the Infimal Convolution-type TV regulariser \cite{ChambolleLions,HollerKunischIC2014,journal_tvlp} and the TGV regularisation \cite{TGV} functional are nowadays popular choices combining good TV features with more natural reconstruction properties. Regarding the choice of the parameter $\lambda$, we remark that space dependent parameters, i.e. functions adapting to the underlying structures in the image have been also considered and fine studies on their reconstruction properties have been shown \cite{}. The optimal choice of (spatially dependent) weighting parameters in \eqref{functional_energy} is a challenging problems for which, depending on the availability of prior information  on the noise in the data, several strategies have been discussed. In particular, when the noise level is known \emph{a priori} discrepancy-type techniques \'a la Morozov have been proposed for general noise models in \cite{}, while in the case of blind denoising recent techniques are based on the use of training set of examples and bilevel optimisation techniques \cite{}.

In this work we will mostly focus on  the standard TV regularisation energy, as our work focuses on the choice of $\Phi$ rather than of $R$. 
% we mainly focus on a specific choice of the data fitting term $\Phi(u,f)$ in \eqref{functional_energy} modelling noise mixtures in the data. 
Classical data fidelity terms for denoising images with Gaussian or impulsive Salt \& Pepper noise are  based on the use of the $L^2$ and $L^1$ norm, respectively, i.e.,
\begin{equation*}
\Phi_{L^2}(u,f):=\frac{1}{2}\int_\Omega (f-u)^2\,dx\qquad\text{or}\qquad \Phi_{L^1}(u,f):=\int_\Omega |f-u|\,dx.
\end{equation*}
These discrepancies are statistically consistent with the assumptions on the noise since they can be derived as the MAP estimators of the underlying likelihood function \cite{benning2011error}.

There exists a considerable amount of work in the literature regarding the structures of solutions of variational problems with pure $L^{1}$ or $L^{2}$ fidelity terms. As a result, the differences between the effects that these terms have on the solutions of the corresponding denoising models are well understood.  See for instance \cite{Allard1, Allard2, Allard3, caselles2007discontinuity, poon_TV_geometric, chanL1,cristoferi, duvalL1, Jalalzai2015jmiv, meyer2001oscillating, ring2000structural, valkonen2014jump1} for TV regularisation and \cite{BrediesL1, Papafitsoros_Bredies,tgv_asymptotic, poschl2013exact, valkonen2014jump2} for TGV. For example, in the case of TV regularisation, it is known that the use of $L^{2}$ fidelity does not introduce new discontinuities in the solution,  which is not the case for the $L^{1}$ fidelity. Moreover, the $L^{1}$ model is capable of exact data recovery, in contrast to the $L^{2}$ one, where always some loss of contrast occurs.

\subsection{Image denoising for noise mixtures} \label{sec:mixture}

Due to different image acquisition and transmission faults, the given image $f$ may be corrupted by a mixture of noise statistics. This is typical, for instance, whenever the presence of noise is due to electronic faults and/or photon-counting processes (such as in astronomy and microscopy applications) combined with actual damages in the receiving sensors, resulting in a lack of information transmittance in only a few image pixels (``burned pixels''). The modelling of $\Phi$ in  \eqref{functional_energy} is therefore expected to encode such noise combination.
 %and has been studied in previous work in the literature.
  In this work we consider the special case when Gaussian and sparse Salt \& Pepper noise are present in the data.

Several authors have considered previously such noise combination. In \cite{langer,langerl1l2}, for instance, a combined model with a $L^1+L^2$ data fidelity and TV regularisation is considered for the joint removal of impulsive and Gaussian noise. Another approach is considered in \cite{impulsegauss2008} where two sequential steps having $L^1$ and $L^2$ as data fidelity are performed to remove the impulsive and the Gaussian component of the mixed noise, respectively. Framelet-based approaches combining $L^1$ and $L^2$ data fidelities in a discrete setting have also been proposed in \cite{Dong2012,Yan2013}. 
However, despite the observed good practical performance of the models described above, the use of an additive and sequential combination of $L^1$ and $L^2$ data fitting terms lacks a rigorous statistical interpretation in terms, for instance, of a MAP estimation. 

We remark that the combination of other noise distributions such as, for instance, Gaussian and Poisson is also frequent in astronomy and microscopy applications and have been studied in several works such as, for instance,  \cite{Benvenuto2008, poissongauss2013, lanza2013}. Such noise mixtures, however, are outside the scope of this work.

\subsection{Infimal convolution modelling of data discrepancies}

Recently in  \cite{calatroni_mixed}, a non-standard variational model for noise removal of mixtures of Salt \& Pepper and Gaussian, and Gaussian and Poisson noise has been studied. The model, which will be referred to as $\tv$--$\mathrm{IC}$ model, is based on the minimisation of an energy functional which is the sum of $|Du|(\om)$ and an infimal convolution of single noise data discrepancy terms. Given two positive parameters $\lambda_{1},\lambda_{2}$ it reads:
\begin{equation*}
\min_{u\in\bv(\om)} |Du|(\om)+\Phi^{\lambda_{1},\lambda_{2}} (u,f),
\end{equation*}
 %and $\mathcal{A}$ is an admissible set of functions. 
where the data fidelity $\Phi^{\lambda_{1},\lambda_{2}} (u,f)$ is defined as 
\begin{equation}\label{data_fit_def}
\Phi^{\lambda_{1},\lambda_{2}} (u,f):=\inf_{v}\; \lambda_{1}\Phi_{1}(v)+\lambda_{2} \Phi_{2}(v,f-u).
\end{equation}
Here, $\Phi_{1}, \Phi_{2}$ denote standard data fidelity terms typically used for single noise removal such as the $L^1, L^2$ norm and the Kullback-Leibler functional. In the particular case of a mixture of  Salt \& Pepper and Gaussian noise, \eqref{data_fit_def} specifies into $\Phi_{1}(v)=\|v\|_{L^{1}(\om)}$ and $\Phi_{2}(v,f-u)=\frac{1}{2}\|f-u-v\|_{L^{2}(\om)}^{2}$. In this case, the minimisation in \eqref{data_fit_def} is done over $L^{1}(\om)$ and $\Phi^{\lambda_{1},\lambda_{2}} (u,f)$ reads
\begin{equation}\label{L1L2_fid_def}
\Phi^{\lambda_{1},\lambda_{2}} (u,f)=\min_{v\in L^{1}(\om)} \lambda_{1}\|v\|_{L^{1}(\om)}+\frac{\lambda_{2}}{2} \|f-u-v\|_{L^{2}(\om)}^{2}, \quad f,u\in L^{1}(\om).
\end{equation}
It can be easily checked that the minimisation in \eqref{L1L2_fid_def} is indeed well-defined.
The reader should not be alerted by the fact that the $L^{1}$ functions $f,u$ appear in the $L^{2}$ part of \eqref{L1L2_fid_def}, as the variable $v$ takes care of any non-integrability issue, see Proposition \ref{lbl:HuberL1L2}. In fact, as we are going to remark in the same proposition,  the functional $\Phi^{\lambda_{1},\lambda_{2}}$ can be written equivalently as
\begin{equation*}
\Phi^{\lambda_{1},\lambda_{2}} (u,f)=\int_{\om} \varphi(f(x)-u(x))\,dx, 
\end{equation*}
where $\varphi$ is the well-known Huber-regularisation of the absolute value function. As a consequence, the functional $\Phi^{\lambda_{1},\lambda_{2}}$ can be simply seen as a Huberised $L^{1}$ norm. 

In \cite{calatroni_mixed}, it was shown that the data discrepancy \eqref{L1L2_fid_def} corresponds to the joint MAP estimator for a denoising problem featuring a mixture of  Laplace and Gaussian noise distributions. The effectiveness of the model for the removal of such noise mixture as well as the additional property of decomposing the noise into its sparse (Salt \& Pepper) and distributed (Gaussian) component was there confirmed with extended numerical examples.

Note that a Huber smoothing of the $L^{1}$ norm has previously been considered in order to apply fast second order minimisation algorithms such as semismooth Newton method in \cite{HintL1}. Also, in the purely discrete setting, smoothed $\tv$--$L^{1}$ models have  been studied in \cite{histo} for exact histogram specification. Similar models (among which also Huber-type) were also considered in \cite{Baus2014}, where the authors obtained  bounds on the infinity norm of the difference between data and solutions.

% In what follows, we will take a closer look at the structure of solutions of the variational model \eqref{general_model} with the choice \eqref{L1L2_fid_def}, based on an alternative formulation of it. In particular, we will first focus on the the asymptotics as $\lambda_{1}$ or $\lambda_{2}$ tend to infinity or zero and describe how the solution of the model considered vary in these cases. Next, we will focus on the one dimensional case and give a precise description of the type of solutions one can expect in terms of some quantitative bounds on the ratio of the parameters $\lambda_1$ and $\lambda_2$. \blue{Change here if we have general regulariser and if we consider Gauss-Poisson case as well}. Our analysis follows the analogous study for the total generalised variation (TGV) functional, which has been done in \cite{tgv_asymptotic}.

%In this paper, we make a finer analysis of the $\mathrm{IC}$--$\tv$ model \eqref{general_model}-\eqref{L1L2_fid_def} proposed in \cite{calatroni_mixed} for the joint removal of Gaussian and Salt \& Pepper noise. In particular, we refine the asymptotic results given in \cite[Proposition 5.1]{calatroni_mixed} and study the convergence of the minimisers of \eqref{general_model} for large values of the parameters $\lambda_1$ and $\lambda_2$ in comparison to the solution of classical models used for denoising of single denoising problems. For more insights, we study also exact solutions in the one dimensional case, and confirm our theoretical results with some numerical experiments. 

\paragraph{Our contribution.}

In this work we examine in depth the similarities and the differences between the $\tv$--$\mathrm{IC}$ model and the pure $\tv$--$L^{1}$, $\tv$--$L^{2}$ ones. 
We first provide detailed asymptotic results as $\lambda_{1}$ or $\lambda_{2}$ tend either to infinity or to zero and describe how the solution of the model  varies in these cases. Note that these results are proved for a general regularisation term. Secondly, in order to have a better insight on the type of solutions one can expect, we do a fine scale analysis of the one-dimensional TV regularised model by computing exact solutions for simple data functions $f$. Up to our knowledge, this is the first time that the effect of the Huberised $L^{1}$ fidelity term is studied in the continuous setting.

In the second part of the paper, we focus on the optimal selection of the parameters $\lambda_{1}, \lambda_{2}$ appearing in \eqref{L1L2_fid_def}. In order to do that, we consider a bilevel optimisation approach  \cite{bilevellearning, noiselearning} which in its general formulation reads
\begin{equation}\label{bilevel_unreg_intro}
\begin{aligned}
&\min_{\lambda_{1},\lambda_{2}~\ge 0}~ F(u_{\lambda_{1},\lambda_{2}})\\
\text{subject to }\quad u_{\lambda_{1},\lambda_{2}}\in~ &\underset{u\in\bv(\om)}{\operatorname{argmin}}~ |Du|(\om)+ \Phi^{\lambda_{1},\lambda_{2}}(u,f).
\end{aligned}
\end{equation}
Here $F$ denotes a cost functional which measures how far the solution $u_{\lambda_{1},\lambda_{2}}$ is from some ground truth (training) image. The parameters $\lambda_{1},\lambda_{2}$ selected within this framework are therefore those producing the closest reconstruction $u_{\lambda_{1},\lambda_{2}}$ to the training image, see Section \ref{sec:bilevel} for more details.
We perform a rigorous analysis on the existence of solutions of \eqref{bilevel_unreg_intro} as well as a regularised version of it, proving the Fr\'echet differentiability of the solution map $\mathcal{S}: (\lambda_1,\lambda_2)\mapsto u_{\lambda_1,\lambda_2}$ and the existence of adjoint states. This allows to derive a handy characterisation of the gradient of the reduced form of $F$ which can be used for efficient numerical implementations. Our analysis justifies rigorously the formal Lagrangian approach considered in \cite[Section 7]{calatroni_mixed}. 

We conclude our study with some numerical experiments connecting the analysis on the structure of solutions discussed above with the problem of learning the optimal noise model for a given noisy image with unknown noise intensity. Our numerical findings show that in case of pure Salt \& Pepper and Gaussian denoising, the bilevel optimisation approach applied to the $\tv$--$\mathrm{IC}$ model computes optimal parameters which enforces pure $L^1$ fidelity. In the case of noise mixture, a combination of $L^1$ and $L^2$ data fitting is preferred. Interestingly, in the case of pure Gaussian, it is not the pure $L^{2}$ data fitting that is selected but still a combination of $L^{1}$ and $L^{2}$, indicating the benefit of the use of $L^{1}$ discrepancy even in the case of Gaussian noise.

%This is in fact not surprising as it is in good agreement with the popular use of these models for this type of problems as discussed, e.g., in \cite{nikolova04,}. \textcolor{red}{Add reference Mila's work: TV-L1 for Gaussian}

We emphasise that the two parts of the paper are intrinsically connected, see for instance Propositions \ref{lbl:non-wellposened} and \ref{lbl:existence_bilevel}. In the former, by making use of the analytical results of the first part, we show with a help of a counterexample, that in order to show existence of solutions for the bilevel optimisation problem \eqref{bilevel_unreg_intro}, it is necessary to enforce an upper bound on the parameters $\lambda_{1}, \lambda_{2}$. The existence of solutions in this case, is shown in Proposition \ref{lbl:existence_bilevel} also by making use of the results of the first part of the paper.

Overall, this study motivates further the use of the $\tv$--$\mathrm{IC}$ model and in general the use of the infimal convolution based fidelity term, by (i) describing the structure of the expected solutions and  (ii) by proposing an automated parameter selection strategy making this model more flexible and applicable to mixed denoising problems.

\section{Analysis of the $L^{1}$--$L^{2}$ IC model: characterisation and asymptotics}   \label{sec:analysis}

%\subsection{$L^{1}$--$L^{2}$ infimal convolution as a Huberised $L^{1}$ norm.}

We start this section by observing that the $L^{1}$--$L^{2}$ infimal convolution term can be equivalently formulated as a Huberised $L^{1}$ norm. This provides an interesting motivation on its effectiveness in the removal of mixed Salt \& Pepper and Gaussian noise. As usual, $\Omega\subseteq \RR^{d}$ denotes an open, bounded, connected domain with Lipschitz boundary.

\newtheorem{HuberL1L2}{Proposition}[section]
\begin{HuberL1L2}\label{lbl:HuberL1L2}
Let $\lambda_{1},\lambda_{2}>0$, $f,u\in L^{1}(\om)$ and consider the data fitting term, $\Phi^{\lambda_{1},\lambda_{2}} (u,f)$  defined in \eqref{L1L2_fid_def}. Then 
\begin{equation}\label{HuberL1L2a}
\Phi^{\lambda_{1},\lambda_{2}} (u,f)=\int_{\om} \varphi(f(x)-u(x))\,dx, 
\end{equation}
where for $t\in\RR$
\begin{equation}\label{HuberL1L2b}
\varphi(t)=
\begin{cases}
\lambda_{1} |t|-\frac{\lambda_{1}^{2}}{2\lambda_{2}}, & \text{ if}\quad |t|\ge \frac{\lambda_{1}}{\lambda_{2}},\\
\frac{\lambda_{2}}{2} |t|^{2},						& \text{ if}\quad |t|<\frac{\lambda_{1}}{\lambda_{2}}.
\end{cases}
\end{equation}
\end{HuberL1L2}

\begin{proof}
The proof is straightforward, in view of 
\[\Phi^{\lambda_{1},\lambda_{2}} (u,f)=\min_{v\in L^{1}(\om)} \int_{\om} \Big( \lambda_{1} |v(x)|+\frac{\lambda_{2}}{2} (f(x)-u(x)-v(x))^{2}\Big) dx, \]
and noticing that the minimisation in the definition of $\Phi^{\lambda_{1},\lambda_{2}} $ can be considered pointwise. Immediate calculations show that the optimal $v$ can be computed explicitly as
% Thus we can focus on the following real-valued minimisation problem
% \begin{equation}\label{min_1d}
% \min_{s\in \RR}\; \lambda_{1} |s| +\frac{\lambda_{2}}{2} (c-d-s)^{2},
% \end{equation}
% or equivalently, by setting $t:=c-d-s$,
% \begin{equation}\label{min_1d_equi}
% \min_{t\in \RR}\; \lambda_{1} |c-d-t| +\frac{\lambda_{2}}{2} t^{2},
% \end{equation}
% After some straightforward calculations one finds
% \[\min_{t\in\RR} \lambda_{1} |c-d-t| +\frac{\lambda_{2}}{2} t^{2}=
% \begin{cases}
% \lambda_{1}|c-d| - \frac{\lambda_{1}^{2}}{2 \lambda_{2}}, & \text{ if}\quad |c-d|\ge \frac{\lambda_{1}}{\lambda_{2}},\\
% \frac{\lambda_{2}}{2}|c-d|^{2}, 						   & \text{ if}\quad |c-d|< \frac{\lambda_{1}}{\lambda_{2}}, 
% \end{cases}
% \]
% with the optimal value $t_{opt}$ being
% \[t_{opt}=
% \begin{cases}
% \frac{\lambda_{1}}{\lambda_{2}} \frac{c-d}{|c-d|}, & \text{ if}\quad |c-d|\ge \frac{\lambda_{1}}{\lambda_{2}}, \\
% c-d,										   & \text{ if}\quad 	|c-d|< \frac{\lambda_{1}}{\lambda_{2}}. 
% \end{cases}
% \]
% Then, the required form \eqref{HuberL1L2a}-\eqref{HuberL1L2b} for $\Phi^{\lambda_{1},\lambda_{2}}$ follows immediately with the following definition for the optimal $v_{opt}$ for $x\in\Omega$:
\begin{equation}\label{vopt}
v_{opt}(x)=
\begin{cases}
f(x)-u(x)-\frac{\lambda_{1}}{\lambda_{2}}\frac{f(x)-u(x)}{|f(x)-u(x)|}, & \text{ if}\quad |f(x)-u(x)|\ge \frac{\lambda_{1}}{\lambda_{2}},\\
0, & \text{ if}\quad |f(x)-u(x)|< \frac{\lambda_{1}}{\lambda_{2}}.
\end{cases}
\end{equation}
\end{proof}

From the formulation \eqref{HuberL1L2a}--\eqref{HuberL1L2b} one sees that the infimal convolution of $L^{1}$ and $L^{2}$ data fidelities coincides with a smoothed $L^{1}$ norm, see Figure \ref{fig:huber}. This is in fact well-known in the context of optimisation in Hilbert spaces \cite{bauschkecombettes}, where an explicit expression like \eqref{vopt} is often used to compute  the soft-thresholding operators. 
Furthermore, let us point out here that Proposition \ref{lbl:HuberL1L2} above is  analogous to a similar result about the Huberised total variation functional, which has also been  shown to be equivalent to a corresponding  infimal convolution functional involving $L^{1}$ and $L^{2}$ norms, see \cite{journal_tvlp} and \cite{stadler}.

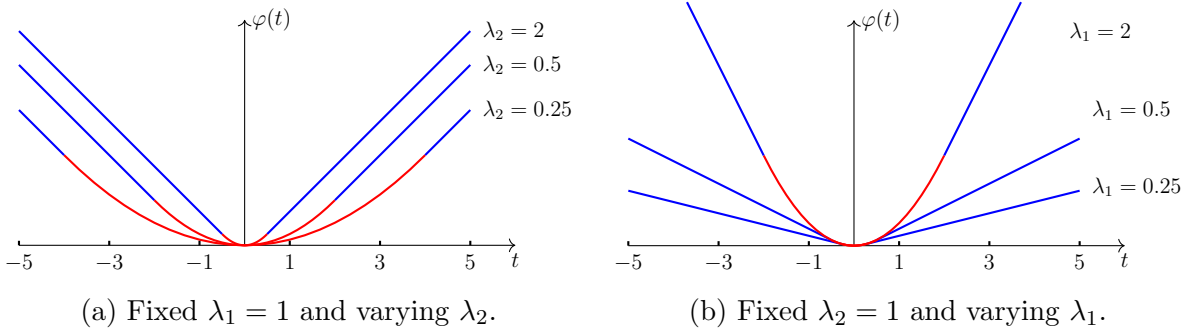
\begin{figure}[!h]
\begin{subfigure}{0.48\textwidth}
\centering
\begin{tikzpicture}[thick, scale=0.6, every node/.style={scale=0.7}]
      \draw [thin] [->] (-5,0)--(6,0) node[right, below] {$t$};
      \draw [thin] [->] (0,0)--(0,5) node[above, right] {$\varphi(t)$};
        \draw [color=red, domain=-4:4, variable=\x]  plot({\x}, {0.125*\x*\x});   %lambda2=0.25
      \draw [color=blue, domain=-4:-5, variable=\x]  plot({\x}, {abs{\x}-2});
      \draw [color=blue, domain=4:5, variable=\x]  plot({\x}, {\x-2});
       \draw (5.1,3) node[right]{$\lambda_{2}=0.25$};
      \draw [color=red, domain=-2:2, variable=\x]  plot({\x}, {0.25*\x*\x}); %lambda2=0.5
      \draw [color=blue, domain=-2:-5, variable=\x]  plot({\x}, {abs{\x}-1});
      \draw [color=blue, domain=2:5, variable=\x]  plot({\x}, {\x-1});
       \draw (5.1,4) node[right]{$\lambda_{2}=0.5$};
     \draw [color=red, domain=-0.5:0.5, variable=\x]  plot({\x}, {1*\x*\x}); %lambda2=2
      \draw [color=blue, domain=-0.5:-5, variable=\x]  plot({\x}, {abs{\x}-0.25});
      \draw [color=blue, domain=0.5:5, variable=\x]  plot({\x}, {\x-0.25});
          \draw (5.1,4.75) node[right]{$\lambda_{2}=2$};
            \foreach \i in {-5,-3,..., 3,5} {
        \draw (\i,-0.0) -- (\i,0.1) node[below=0.8mm] {$\i$};
    }
\end{tikzpicture}
\subcaption{Fixed $\lambda_{1}=1$ and varying $\lambda_2$. }
\label{fig:huber}
\end{subfigure}
\begin{subfigure}{0.48\textwidth}
\centering
\begin{tikzpicture}[thick, scale=0.6, every node/.style={scale=0.7}]
      \draw [thin] [->] (-5,0)--(6,0) node[right, below] {$t$};
      \draw [thin] [->] (0,0)--(0,5) node[above, right] {$\varphi(t)$};
        \draw [color=red, domain=-0.25:0.25, variable=\x]  plot({\x}, {0.5*\x*\x});   %lambda2=0.25
      \draw [color=blue, domain=-0.25:-5, variable=\x]  plot({\x}, {0.25*abs{\x}-0.03125});
      \draw [color=blue, domain=0.25:5, variable=\x]  plot({\x}, {0.25*\x-0.03125});
       \draw (5.1,1.3) node[right]{$\lambda_{1}=0.25$};
      \draw [color=red, domain=-0.5:0.5, variable=\x]  plot({\x}, {0.5*\x*\x}); %lambda2=0.5
      \draw [color=blue, domain=-0.5:-5, variable=\x]  plot({\x}, {0.5*abs{\x}-0.125});
      \draw [color=blue, domain=0.5:5, variable=\x]  plot({\x}, {0.5*\x-0.125});
       \draw (5.1,3) node[right]{$\lambda_{1}=0.5$};
     \draw [color=red, domain=-2:2, variable=\x]  plot({\x}, {0.5*\x*\x}); %lambda2=2
      \draw [color=blue, domain=-2:-3.7, variable=\x]  plot({\x}, {2*abs{\x}-2});
      \draw [color=blue, domain=2:3.7, variable=\x]  plot({\x}, {2*\x-2});
          \draw (4.6,4.75) node[right]{$\lambda_{1}=2$};
            \foreach \i in {-5,-3,..., 3,5} {
        \draw (\i,-0.0) -- (\i,0.1) node[below=0.8mm] {$\i$};
    }
\end{tikzpicture}
\subcaption{Fixed $\lambda_{2}=1$ and varying $\lambda_1$. }
\end{subfigure}
\caption{Example plots of the Huber function $\varphi$.}
\label{fig:huber}
\end{figure}

\begin{figure}
\begin{subfigure}{0.3\textwidth}
\centering
\begin{tikzpicture}
\draw (0,0) node{\includegraphics[height=4cm]{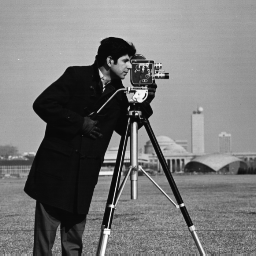}};
\end{tikzpicture}
\subcaption{Original image $u$}
\end{subfigure}
\begin{subfigure}{0.3\textwidth}
\centering
\begin{tikzpicture}
\draw (0,0) node{\includegraphics[height=4cm]{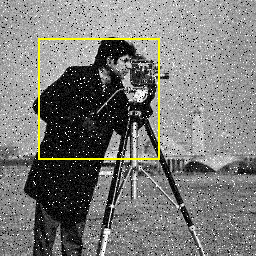}};
%\draw[red] (0.25,0.7) -- (1.25,0.7) -- (1.25,1.7) -- (0.25,1.7) -- (0.25,0.7);
\draw[yellow, line width=0.59mm] (-1.39,-0.48) -- (0.48,-0.48)-- (0.48,1.39)--(-1.39,1.39)--(-1.39,-0.5095);
\end{tikzpicture}
\subcaption{Data $f$}
\end{subfigure}
\begin{subfigure}{0.3\textwidth}
\centering
\begin{tikzpicture}
\draw (0,0) node{\includegraphics[trim=1.5cm 3.5cm 3.48cm 1.5cm, clip,height=4cm]{cameraman_noisy_selec}};
\draw[red, line width=0.3mm]  (-1.62,0.985) -- (0.33,0.985);
\end{tikzpicture}
\subcaption{Data $f$ -- detail}
\end{subfigure}

\begin{subfigure}{0.9\textwidth}
\centering
\resizebox{0.9\textwidth}{!}{
% This file was created by matlab2tikz.
%
%The latest updates can be retrieved from
%  http://www.mathworks.com/matlabcentral/fileexchange/22022-matlab2tikz-matlab2tikz
%where you can also make suggestions and rate matlab2tikz.
%
\begin{tikzpicture}

\begin{axis}[%
width=10in,
height=4.392in,
at={(1.41in,0.593in)},
scale only axis,
xmin=1,
xmax=55,
xtick={0,5,...,55},
xlabel style={font=\color{black}},
ymin=0,
ymax=1,
ylabel style={font=\color{black}},
ylabel={\LARGE{$|f(x)-u(x)|$}},
axis background/.style={fill=white}
]
\addplot [color=black, dashed, line width=0.4mm, forget plot]
  table[row sep=crcr]{%
1	0.0862745098039215\\
2	0.00392156862745097\\
3	0.0901960784313725\\
4	0.0941176470588235\\
5	0\\
6	0.129411764705882\\
7	0.180392156862745\\
8	0.0235294117647059\\
9	0.0901960784313726\\
10	0.00392156862745097\\
11	0.0549019607843138\\
12	0.00392156862745097\\
13	0.0901960784313726\\
14	0.949019607843137\\
15	0.133333333333333\\
16	0.0156862745098039\\
17	0\\
18	0.0392156862745098\\
19	0.0549019607843137\\
20	0.945098039215686\\
21	0.0666666666666667\\
22	0.0313725490196078\\
23	0.92156862745098\\
24	0.0235294117647059\\
25	0.0745098039215686\\
26	0.0549019607843137\\
27	0.0509803921568627\\
28	0.0352941176470588\\
29	0.0627450980392157\\
30	0.0156862745098039\\
31	0.0196078431372549\\
32	0.101960784313725\\
33	0.0627450980392157\\
34	0.00392156862745098\\
35	0.0470588235294118\\
36	0.0509803921568627\\
37	0.129411764705882\\
38	0.00392156862745098\\
39	0.152941176470588\\
40	0.0549019607843137\\
41	0.109803921568627\\
42	0.0352941176470588\\
43	0.00784313725490196\\
44	0.0313725490196078\\
45	0.0588235294117647\\
46	0.101960784313725\\
47	0.0627450980392157\\
48	0.0235294117647059\\
49	0.670588235294118\\
50	0.352941176470588\\
51	0.00392156862745102\\
52	0\\
53	0.0862745098039216\\
54	0.00392156862745097\\
55	0.0980392156862745\\
56	0.145098039215686\\
};
\addplot [color=blue, draw=none, mark=asterisk, mark options={solid, blue, mark size=4pt}, line width=0.6mm, forget plot]
  table[row sep=crcr]{%
14	0.949019607843137\\
};
\addplot [color=blue, draw=none, mark=asterisk, mark options={solid, blue, mark size=4pt}, line width=0.6mm, forget plot]
  table[row sep=crcr]{%
20	0.945098039215686\\
};
\addplot [color=blue, draw=none, mark=asterisk, mark options={solid, blue, mark size=4pt}, line width=0.6mm, forget plot]
  table[row sep=crcr]{%
23	0.92156862745098\\
};
\addplot [color=blue, draw=none, mark=asterisk, mark options={solid, blue, mark size=4pt}, line width=0.6mm, forget plot]
  table[row sep=crcr]{%
49	0.670588235294118\\
};
\addplot [color=black, draw=none, mark=asterisk, mark options={solid, red,mark size=4pt}, line width=0.6mm,forget plot]
  table[row sep=crcr]{%
1	0.0862745098039215\\
2	0.00392156862745097\\
3	0.0901960784313725\\
4	0.0941176470588235\\
5	0\\
6	0.129411764705882\\
7	0.180392156862745\\
8	0.0235294117647059\\
9	0.0901960784313726\\
10	0.00392156862745097\\
11	0.0549019607843138\\
12	0.00392156862745097\\
13	0.0901960784313726\\
};
\addplot [color=black, draw=none, mark=asterisk,  mark options={solid, red, mark size=4pt}, line width=0.6mm, forget plot]
  table[row sep=crcr]{%
15	0.133333333333333\\
16	0.0156862745098039\\
17	0\\
18	0.0392156862745098\\
19	0.0549019607843137\\
};
\addplot [color=black, draw=none, mark=asterisk,  mark options={solid, red,mark size=4pt}, line width=0.6mm, forget plot]
  table[row sep=crcr]{%
21	0.0666666666666667\\
22	0.0313725490196078\\
};
\addplot [color=black, draw=none, mark=asterisk,  mark options={solid, red,mark size=4pt}, line width=0.6mm, forget plot]
  table[row sep=crcr]{%
24	0.0235294117647059\\
25	0.0745098039215686\\
26	0.0549019607843137\\
27	0.0509803921568627\\
28	0.0352941176470588\\
29	0.0627450980392157\\
30	0.0156862745098039\\
31	0.0196078431372549\\
32	0.101960784313725\\
33	0.0627450980392157\\
34	0.00392156862745098\\
35	0.0470588235294118\\
36	0.0509803921568627\\
37	0.129411764705882\\
38	0.00392156862745098\\
39	0.152941176470588\\
40	0.0549019607843137\\
41	0.109803921568627\\
42	0.0352941176470588\\
43	0.00784313725490196\\
44	0.0313725490196078\\
45	0.0588235294117647\\
46	0.101960784313725\\
47	0.0627450980392157\\
48	0.0235294117647059\\
};
\addplot [color=black, draw=none, mark=asterisk,  mark options={solid, red,mark size=4pt}, line width=0.6mm, forget plot]
  table[row sep=crcr]{%
50	0.352941176470588\\
51	0.00392156862745102\\
52	0\\
53	0.0862745098039216\\
54	0.00392156862745097\\
55	0.0980392156862745\\
56	0.145098039215686\\
57	0.0745098039215686\\
};

 \addplot [color=black, draw=none, mark=asterisk,  mark options={solid, red,mark size=4pt}, line width=0.6mm, forget plot]
  table[]{%
41.5	0.9
};
\addplot [color=black, draw=none, mark=asterisk,  mark options={solid, blue,mark size=4pt}, line width=0.6mm, forget plot]
  table[]{%
41.5	0.8
}; 
 
  \node at (axis cs:45,0.9) [ black ] {\Large{Gaussian}};
 \node at (axis cs:46.25,0.8) [ black ] {\Large{Salt \& Pepper}};
 
\end{axis}
\end{tikzpicture}%
}
\subcaption{$|f(x)-u(x)|$ along the line profile}
\end{subfigure}
\caption{Interepretation of the $L^{1}$--$L^{2}$ IC term based on the formulation \eqref{HuberL1L2a}--\eqref{HuberL1L2b}. First row: original  image, noisy version and detail with line profile (red). Second row: difference $|f-u|$ along the line profile. The fidelity functional $\Phi^{\lambda_{1},\lambda_{2}}$ locally acts like $\|f-u\|_{L^{1}}$ at points with high values of $|f-u|$ (blue asterisks), thus assuming these points to be corrupted by Salt \& Pepper noise. On the other hand, $\Phi^{\lambda_{1},\lambda_{2}}$ acts locally as $\|f-u\|_{L^{2}}^{2}$ at points with low values of $|f-u|$ (red asterisks), identifying these points as corrupted by Gaussian noise.
\textbf{Parameters}: Gaussian variance $\sigma^2=0.005$, density of pixels corrupted by Salt \& Pepper noise $d=5\%$.
 %High values of the difference are likely to indicate image pixels corrupted with Salt \& Pepper noise which would then be treated consistently using an $L^1$-type data fidelity.
 }
\label{heuristic:L1L2}
\end{figure}

The formulation \eqref{HuberL1L2a}--\eqref{HuberL1L2b} of  $\Phi^{\lambda_{1},\lambda_{2}}$ provides an interesting insight on its interpretation and motivates its effectiveness for mixed noise removal. When $|f(x)-u(x)|$ is large, the noise component is interpreted as Salt \& Pepper  by the model and then $\Phi^{\lambda_{1},\lambda_{2}}$ behaves locally as  $\|f-u\|_{L^{1}(\om)}$. On the other hand, if $|f(x)-u(x)|$ is small, then the model assumes that the noise is Gaussian and enforces a data fidelity $\Phi^{\lambda_{1},\lambda_{2}}$ locally equal to $\|f-u\|_{L^{2}(\om)}^{2}$, see Figure \ref{heuristic:L1L2} for a  visualisation.

 %the direct sum of the fidelity terms, lacks this local adaptivity \com{(Langer has not done anything using this model and spatially adapted $\lambda_{}1$ and $\lambda_{2}$ right?)}.

%\blue{Need to check, he may have done something like that recently. Also, we should probably motivate the effectiveness of $\mathrm{IC}$--$\tv$ in comparison, for instance, with the two-phase approach by Nikolova et al. \cite{impulsegauss2008}. where first outliers are identified similarly as for inpainting and then removed. Then, a standard Gaussian denoising is performed on the 'remaining' pixels. The use of a Huberised $L^1$ norm removes this "pre-processing" step and the two models are just used at the same time with the local enforcement of both the $L^1$ and $L^2$ models as described above.}

\subsection{Asymptotic behaviour}
In this section, we investigate the asymptotic behaviour of the $L^{1}$--$L^{2}$ IC model. Here, we do not need to restrict to the $\tv$ regulariser, but we can consider a more general regularisation functional $J$ with the following properties:
\begin{enumerate}
\item  $J:L^{1}(\om)\to \mathbb{R}\cup\{\infty\}$ is positive, proper, convex, lower semicontinuous with respect to the strong convergence in $L^{1}(\om)$.
%\item $J$ is positively homogeneous of degree one, i.e., $J(\lambda u)=|\lambda| J(u)$. \textcolor{red}{(I am not sure if we need this...)}
\item  There exist  constants $C_{1}, C_{2}>0$ such that $C_{1} |Du|(\om)\le J(u) \le C_{2} |Du|(\om)$ for every $u\in\bv(\om)$. 
%\item \emph{(Strict convergence)} For every $u\in L^{p}(\om)\cap\bv(\om) $, $1<p<\infty$, there exists a sequence $(u_{n})_{n\in\NN}\subseteq C^{\infty}(\om)\cap \bv(\om)$ such that
%	\begin{itemize}
%	\item $u_{n}\to u$ in $L^{p}(\om)$ as $n\to\infty$.
%	\item $J(u_{n})\to J(u)$ as $n\to\infty$.5
%	\end{itemize}
%\item \com{Any other property we want for $J$?}
\end{enumerate}

Classical regularisers such as $\tv$, Huber-$\tv$ and $\tgv$ of any order, satisfy the above properties. Note that this is also true for a large class of structural TV-type functionals that are commonly used in inverse problems, see \cite{structuralTV}.
% see for instance \cite{AmbrosioBV} for $\tv$ and \cite{bredies2014regularization} for $\tgv$.

We are interested in the following general problem:
\begin{equation}\label{L1L2_J}
\min_{\substack{u\in\bv(\om) \\ v \in L^{1}(\om)}}  J(u)+ \lambda_{1} \|v\|_{L^{1}(\om)} +\frac{\lambda_{2}}{2} \|f-u-v\|_{L^{2}(\om)}^{2},
\end{equation}
which is a more general version of the $\tv$--$\mathrm{IC}$ model for Gaussian and Salt \& Pepper noise removal
\begin{equation}\label{L1L2_TV}
\min_{\substack{u\in\bv(\om) \\ v \in L^{1}(\om)}}  |Du|(\om) + \lambda_{1} \|v\|_{L^{1}(\om)} +\frac{\lambda_{2}}{2} \|f-u-v\|_{L^{2}(\om)}^{2}.
\end{equation}

The well-posedness of \eqref{L1L2_TV} has been studied in \cite{calatroni_mixed}. For the more general model \eqref{L1L2_J} existence of minimisers $(u^{\ast},v^{\ast})\in \bv(\om)\times L^{1}(\om)$ follows from a simple application of the direct method of calculus of variations. Note however, that the solution $u^{\ast}$ is not necessary unique as $\Phi^{\lambda_{1},\lambda_{2}}(\cdot,f)$ is not strictly convex ($J$ is not necessarily strictly convex either). The same holds for $v^{\ast}$ due to formula \eqref{vopt} which connects it with $u^{\ast}$.

One natural question one may ask is in what degree we can expect to recover the single noise models by sending the parameters $\lambda_{1}, \lambda_{2}$ (or their ratio) to infinity. In the following we answer this question by taking advantage of the formulation \eqref{HuberL1L2a}--\eqref{HuberL1L2b} and using some $\Gamma$-convergence arguments. Firstly, we extend a corresponding proposition that was shown in \cite[Proposition 5.1]{calatroni_mixed}, to the general regulariser case,  adjusted for our purposes.

\newtheorem{asymptotics_calatroni_mixed}[HuberL1L2]{Proposition}

\begin{asymptotics_calatroni_mixed}\label{lbl:asymptotics_calatroni_mixed}
Let $(u^{\ast},v^{\ast})\in \bv(\om)\times L^{1}(\om)$ be an optimal pair for \eqref{L1L2_J}. Then, the following assertions hold:
\begin{enumerate}
\item If $\lambda_{1}\to\infty$, $f\in L^{1}(\om)$ then $v^{\ast}\to 0$ in $L^{1}(\om)$. 
 %If $f\in\bv(\om)$ then the assumption for fixed $\lambda_{2}$ can be dropped.
\item If $\lambda_{2}\to\infty$,  $f\in L^{2}(\om)$ then $\|f-u^{\ast}-v^{\ast}\|_{L^{2}(\om)}\to 0$.
 If in addition $\lambda_{1}$ is fixed, then the same result holds with $f\in L^{1}(\om)$.
% If $f\in\bv(\om)$ then the assumption for fixed $\lambda_{1}$ can be dropped.
\item If both $\lambda_{1},\lambda_{2}\to\infty$ and $f\in L^{2}(\om)$ then (i) holds and we have that $u^{\ast}\to f$ in $L^{1}(\om)$. If $f\in \bv(\om)$ this convergence is also weakly$^{\ast}$ in $\bv(\om)$.
%\footnote{One interesting question here if we can prove something similar when merely $f\in L^{1}(\om)$. Probably not but can we make that more precise? It is not so clear what should happen?}
\end{enumerate}

%\begin{enumerate}
%%\item $v^{\ast}\to 0$ in $L^{1}(\om)$ as $\lambda_{1}\to\infty$.
%%\item $v^{\ast}\to f-u^{\ast}$ in $L^{2}(\om)$ as $\lambda_{2}\to\infty$.
%%\item $(u^{\ast},v^{\ast})\to (f,0)$ in $L^{1}(\om)\times L^{1}(\om)$ as $\lambda_{1},\lambda_{2}\to\infty$.
%\item $|Du^{\ast}|(\om)\to 0$ and $v^{\ast}\to 0$ in $L^{1}(\om)$ as $\lambda_{2}\to 0$.
%\item $|Du^{\ast}|(\om)\to 0$ and $v^{\ast}\to f-u^{\ast}$ in $L^{2}(\om)$ as $\lambda_{1}\to 0$.
%\end{enumerate} 

\end{asymptotics_calatroni_mixed}

\begin{proof}
$(i)$ %Assume that $\lambda_{1}\to \infty$, and that $f\in L^{1}(\om)$.
We notice that 
\begin{equation*}
\lambda_{1}\|v^{\ast}\|_{L^{1}(\om)}\le J(u)+\lambda_{1}\|v\|_{L^{1}(\om)}+\frac{\lambda_{2}}{2}\|f-u-v\|_{L^{2}(\om)}^{2},\quad  \forall u\in\bv(\om),\;v\in L^{1}(\om),
\end{equation*}
which by setting $v=f-u$, implies
\begin{equation}\label{v_vi_2}
\lambda_{1}\|v^{\ast}\|_{L^{1}(\om)}\le J(u)+ \lambda_{1}\|f-u\|_{L^{1}(\om)},\quad  \forall\,u\in\bv(\om),\;v\in L^{1}(\om).
\end{equation}
Given $\epsilon>0$, we can find $u_{\epsilon}\in C_{c}^{\infty}(\om)$ such that $\|f-u_{\epsilon}\|_{L^{1}(\om)}<\epsilon$. Thus \eqref{v_vi_2} becomes
\begin{equation*}
\|v^{\ast}\|_{L^{1}(\om)} \le \ \frac{1}{\lambda_{1}} J(u_{\epsilon}) +\epsilon \quad \Rightarrow 
\limsup_{\lambda_{1}\to\infty} \|v^{\ast}\|_{L^{1}(\om)}\le \epsilon.
\end{equation*}
Since $\epsilon$ was arbitrary, the result follows.
%
%-----------\\
%Setting $u=0$ and and $v$ as in  \eqref{vopt} we get
%\begin{align*}
%\|v^{\ast}\|_{L^{1}(\om)}\le \int_{|f|\ge \frac{\lambda_{1}}{\lambda_{2}}} |f|-\frac{\lambda_{1}}{2\lambda_{2}}\,dx+\frac{\lambda_{2}}{2\lambda_{1}} \int_{|f|<\frac{\lambda_{1}}{\lambda_{2}}} |f|^{2}dx.
%\le  \int_{|f|\ge \frac{\lambda_{1}}{\lambda_{2}}} |f|\,dx+\frac{\lambda_{2}}{2\lambda_{1}}\|f\|_{L^{2}(\om)}^{2}\end{align*}
%Sending $\lambda_{1}\to \infty$ and using dominated convergence we obtain that $v\to 0$ in $L^{1}(\om)$. Now by dropping the assumption $\lambda_{2}$ fixed but assuming that $f\in\bv(\om)$, we have by setting $v=0$, and $u=f$ that again $v^{\ast}\to 0$ in $L^{1}(\om)$. \\
%--------------\\

$(ii)$  %Assume now that $\lambda_{2}\to \infty$, $f\in L^{2}(\om)$. 
In this case we have that for every $u\in\bv(\om)$ and $v\in L^{1}(\om)$
\begin{equation}\label{uv_vi}
\frac{\lambda_{2}}{2}\|f-u^{\ast}-v^{\ast}\|_{L^{2}(\om)}^{2}\le J(u)+ \lambda_{1}\|v\|_{L^{1}(\om)}+\frac{\lambda_{2}}{2}\|f-u-v\|_{L^{2}(\om)}^{2},
\end{equation}
which by setting $v=0$, implies
\begin{equation*}
\frac{\lambda_{2}}{2}\|f-u^{\ast}-v^{\ast}\|_{L^{2}(\om)}^{2}\le J(u)+  \frac{\lambda_{2}}{2}\|f-u\|_{L^{2}(\om)}^{2},\quad  \forall\,u\in\bv(\om).
\end{equation*}
Then, proceeding as in step $(i)$ the result follows. Notice that if we assume that $\lambda_{1}$ is fixed (or more generally bounded from above) and by merely assuming $f\in L^{1}(\om)$, we can have the same result by setting $v=f$ and $u=0$ in \eqref{uv_vi}.

$(iii)$ Notice that if both $\lambda_{1}, \lambda_{2}\to\infty$ and $f\in L^{2}(\om)$ then from $(i)$, $(ii)$ we get that $\|v^{\ast}\|_{L^{1}(\om)}\to 0$ and $\|f-u^{\ast}-v^{\ast}\|_{L^{1}(\om)}\to 0$ and hence an application of the triangle inequality implies that $u^{\ast}\to f$ in $L^{1}(\om)$. If in addition $f\in \bv(\om)$, then we have that $|Du^{\ast}|(\om)$ is uniformly bounded by setting $v=0$, $u=f$ in
\begin{equation*}
C_{1}|Du^{\ast}|(\om)\le J(u^{\ast})\le J(u)+\lambda_{1}\|v\|_{L^{1}(\om)}+\frac{\lambda_{2}}{2}\|f-u-v\|_{L^{2}(\om)}^{2},\quad  \forall\,u\in\bv(\om),\;v\in L^{1}(\om).
\end{equation*}
From compactness in $\bv(\om)$ and the fact that $\|f-u^{\ast}\|_{L^{1}(\om)}\to 0$ we infer that $u^{\ast}\to f$ weakly$^{\ast}$ in $\bv(\om)$.
\end{proof}

In what follows, we refine the result above and prove convergence of the minimisers of \eqref{L1L2_J} to the minimisers of the single noise models. To do so, we  apply $\Gamma$-convergence arguments \cite{dalmasogamma} to the IC term $\Phi^{\lambda_1,\lambda_2}$.

\newtheorem{Gamma_convergence}[HuberL1L2]{Proposition}
\begin{Gamma_convergence}\label{lbl:Gamma_convergence}
Let $f\in L^{1}(\om)$ and let us define the functional $F^{\lambda_{1},\lambda_{2}}:L^{1}(\om)\to \RR^+$ by $F^{\lambda_{1},\lambda_{2}}(u):=\Phi^{\lambda_{1},\lambda_{2}}(u,f)$. Then
\begin{enumerate}
 \item For any fixed $\lambda_{1}$, $F^{\lambda_{1},\lambda_{2}}$\  $\Gamma$-converges to $F_{1}(\cdot):=\lambda_{1}\|f-\cdot\|_{L^{1}(\om)}$ as $\lambda_{2}\to \infty$.
 \item For any fixed $\lambda_{2}$, $F^{\lambda_{1},\lambda_{2}}$\  $\Gamma$-converges to $\overline{F}_{2}(\cdot)$ as $\lambda_{1}\to \infty$,
where for every $u\in L^1(\Omega)$, $\overline{F}_2$ is defined as
 \[\overline{F}_{2}(u):=\frac{\lambda_{2}}{2}\|f-u\|_{L^{2}(\om)}^{2}:=
 \begin{cases}
 \frac{\lambda_{2}}{2}\|f-u\|_{L^{2}(\om)}^{2}, & \text{ if}\quad f-u\in L^{2}(\om),\\
 +\infty, 								     & \text{ if}\quad f-u\in L^{1}(\om)\setminus L^{2}(\om). 
 \end{cases}
  \]
\end{enumerate}
\end{Gamma_convergence}

\begin{proof}
For $(i)$, let $\big(\lambda_{2}^{(n)}\big)_{n\in\NN}$ be a sequence with $\lambda_{2}^{(n)}\to\infty$ and set $F^{n}:= F^{\lambda_{1},\lambda_{2}^{(n)}}$.
We notice that $F^{n}$ converges uniformly to $\lambda_{1}\|f-\cdot\|_{L^{1}(\om)}$. Indeed, for $u\in L^{1}(\om)$, we have
\begin{align*}
\left | F^{n}(u)-\lambda_{1}\|f-u\|_{L^{1}(\om)} \right |
&= \left |\int_{\om} \varphi(f-u)-\lambda_{1}|f-u|\,dx \right |\\
&\le \int_{|f-u|\ge \frac{\lambda_{1}}{\lambda_{2}^{(n)}}} \frac{\lambda_{1}^{2}}{2\lambda_{2}^{(n)}}\,dx+
\int_{|f-u|<\frac{\lambda_{1}}{\lambda_{2}^{(n)}}} \left |\frac{\lambda_{2}^{(n)}}{2}|f-u|^{2} -\lambda_{1}|f-u|  \right|\,dx\\
%&\le  \int_{\om}\frac{\lambda_{1}^{2}}{2\lambda_{n}}\,dx +\int_{|f-u|<\frac{\lambda_{1}}{\lambda_{n}}}\left |\frac{\lambda_{1}}{2} |f-u|-\frac{\lambda_{1}^{2}}{\lambda_{n}} \right |\,dx\\
& \le\int_{\om}\frac{\lambda_{1}^{2}}{2\lambda_{2}^{(n)}}\,dx+ \int_{\om } \frac{\lambda_{1}^{2}}{2\lambda_{2}^{(n)}}\,dx+ \int_{\om}\frac{\lambda_{1}^{2}}{\lambda_{2}^{(n)}}\,dx\to 0 \quad \text{as} \quad n\to\infty.
\end{align*}
Since the last limit is independent of $u$, the convergence of the functionals is indeed uniform.
Moreover $F^{n}$ is continuous with respect to the $L^{1}$ topology, see for instance \cite{calatroni_mixed}. Thus from \cite[Proposition 5.2]{dalmasogamma} we immediately get that $F^{n}$ $\Gamma$-converges to  $F_{1}$ as $n\to\infty$.

For $(ii)$, we now set $F^{n}:=F^{\lambda_{1}^{(n)},\lambda_{2}}$ with $\lambda_{1}^{(n)}\to \infty$ and observe that $F^{n}$ converges pointwise to $\overline{F}_{2}$. Indeed if $f-u\in L^{1}(\om)\setminus L^{2}(\om)$ then
\[F^{n}(u)\ge \int_{|f-u|<\frac{\lambda_{1}^{(n)}}{\lambda_{2}}}\frac{\lambda_{2}}{2} |f-u|^{2}dx\to\infty\quad \text{as } n\to\infty.\]
On the other hand, if $f-u\in L^{2}(\om)$, proceeding as before we have
\begin{align*}
\left | F^{n}(u)-\frac{\lambda_{2}}{2} \|f-u\|_{L^{2}(\om)}^{2} \right | 
%&= \left | \int_{|f-u|\ge \frac{\lambda_{n}}{\lambda_{2}}} \lambda_{n}|f-u|-\frac{\lambda_{n}^{2}}{2\lambda_{2}}\,dx-\int_{|f-u|\ge \frac{\lambda_{n}}{\lambda_{2}}}\frac{\lambda_{2}}{2} |f-u|^{2}dx  \right |\\
&\le \int_{|f-u|\ge \frac{\lambda_{1}^{(n)}}{\lambda_{2}}} \lambda_{1}^{(n)}|f-u|-\frac{\big(\lambda_{1}^{(n)}\big)^{2}}{2\lambda_{2}}\,dx+ \int_{|f-u|\ge \frac{\lambda_{1}^{(n)}}{\lambda_{2}}}\frac{\lambda_{2}}{2} |f-u|^{2}dx\\
&\le \int_{|f-u|\ge \frac{\lambda_{1}^{(n)}}{\lambda_{2}}} \frac{3}{2}\lambda_{2} |f-u|^{2}dx \to 0 \quad \text{as} \quad n\to\infty.
\end{align*}
If $\big(\lambda_{1}^{(n)}\big)_{n\in\NN}$ is increasing, it can be easily verified that the sequence $(F_{n})_{n\in\NN}$ is increasing. Moreover the $L^{2}$ norm is lower semicontinuous with respect to the strong $L^{1}$ topology. Hence from \cite[Remark 5.5]{dalmasogamma}, we have that $F^{n}$ $\Gamma$-converges to $\overline{F}_{2}$. In the case where  $\big(\lambda_{1}^{(n)}\big)_{n\in\NN}$ is non-monotonically going to infinity, we can find an increasing subsequence and then the result follows from the Urysohn property of $\Gamma$-convergence, see \cite[Proposition 8.3]{dalmasogamma}.
 \end{proof}
 
As a corollary, we obtain the following result on the convergence of minimisers.

\newtheorem{u_l2_to_inf}[HuberL1L2]{Corollary}
\begin{u_l2_to_inf}[Convergence to single noise models]\label{lbl:u_l2_to_inf}
The following two results hold:
\begin{enumerate}
\item Let $f\in L^{1}(\om)$, $\lambda_{1}>0$ fixed and $\big(\lambda_{2}^{(n)}\big)_{n\in\NN}$ with $\lambda_{2}^{(n)}\to \infty$. If $(u_{n})_{n\in\NN}\subseteq \bv(\om)$ is a sequence of minimisers of \eqref{L1L2_J}, then every subsequence of $(u_{n})_{n\in\NN}$ has a weak$^{\ast}$ in $\bv(\om)$ cluster point which is a minimiser of 
\begin{equation}\label{L1J_conv}
\min_{u\in\bv(\om)} J(u)+ \lambda_{1}\|f-u\|_{L^{1}(\om)}.
\end{equation}
Moreover if the solution $u^{\ast}$ of \eqref{L1J_conv} is unique, then $u_{n}\to u^{\ast}$ weakly$^{\ast}$ in $\bv(\om)$ and 
\begin{equation}\label{convergence_of_energies_L1}
J(u_{n})+\Phi^{\lambda_{1},\lambda_{2}^{(n)}}(u_{n},f)\to J(u^{\ast})+\lambda_{1}\|f-u^{\ast}\|_{L^{1}(\om)}.
\end{equation}
\item Let $f\in L^{2}(\om)$, $\lambda_{2}>0$ fixed and $\big(\lambda_{1}^{(n)}\big)_{n\in\NN}$ with $\lambda_{1}^{(n)}\to\infty$. If $(u_{n})_{n\in\NN}\subseteq \bv(\om)$ is a sequence of minimisers of \eqref{L1L2_J}, then $u_{n}\to u^{\ast}$ weakly$^{\ast}$ in $\bv(\om)$, where $u^{\ast}$ is the unique minimiser of 
 \begin{equation}\label{L2J_conv}
\min_{u\in\bv(\om)}  J(u)+ \frac{\lambda_{2}}{2}\|f-u\|_{L^{2}(\om)}^{2}.
\end{equation}
Moreover
\begin{equation*}
J(u_{n})+\Phi^{\lambda_{1}^{(n)},\lambda_{2}}(u_{n},f)\to J(u^{\ast})+ \frac{\lambda_{2}}{2}\|f-u^{\ast}\|_{L^{2}(\om)}^{2}.
\end{equation*}
\end{enumerate}
%Let $u_{n}\in \bv(\om)$, $n\in\NN$, be a minimiser of \eqref{L1L2_TV} with data $f\in L^{1}(\om)$, $\lambda_{1}>0$ fixed and $\lambda_{2}=\lambda_{n}$ where $\lambda_{n}\to\infty$. Then every subsequence of $(u_{n})_{n\in\NN}$ has a cluster point which is a minimiser of 
%\begin{equation}\label{L1TV_conv}
%\min_{u\in\bv(\om)} \lambda_{1}\|f-u\|_{L^{1}(\om)}+|Du|(\om).
%\end{equation}
%Moreover if the solution $u^{\ast}$ of \eqref{L1TV_conv} is unique, then $u_{n}\to u$ weakly$^{\ast}$ in $\bv(\om)$.  
%%\blue{How can the solution of \eqref{L1TV_conv} be unique? You mean for large choice of parameters \'a la Meyer?} \com{sorry, typo.. I meant IF the solution is unique...}
%%\com{or also weakly$^{\ast}$ in BV?} and 
%\begin{equation}\label{convergence_of_energies}
%\Phi^{\lambda_{1},\lambda_{n}}(u_{n},f)+|Du_{n}|(\om)\to \lambda_{1}\|f-u\|_{L^{1}(\om)}+|Du|(\om).
%\end{equation}
 \end{u_l2_to_inf}
\begin{proof}
$(i)$ Since the functional $J$ is lower semicontinuous with respect to $L^{1}$, then \cite[Proposition 6.25]{dalmasogamma} we have that the minimising functionals $\Phi^{\lambda_{1},\lambda_{2}^{(n)}}(\cdot,f)+J(\cdot)$ also $\Gamma$-converge to the functional in \eqref{L1J_conv}. Moreover, note that $(u_{n})_{n\in\NN}$ is uniformly bounded in $\bv(\om)$. Since $u_{n}$ is a minimiser we have in fact that for every $u\in\bv(\om)$ and every $v\in L^1(\om)$, the following three inequalities hold
\begin{align}
C_{1}|Du_{n}|(\om)\le J(u_{n})
&\le J(u)+ \lambda_{1}\|v\|_{L^{1}(\om)}+\frac{\lambda_{2}^{(n)}}{2}\|f-u-v\|_{L^{2}(\om)}^{2}, \label{Ju_vi}\\
\frac{\lambda_{2}^{(n)}}{2}\|f-u_{n}-v_{n}\|_{L^{2}(\om)}^{2}&\le J(u)+\lambda_{1}\|v\|_{L^{1}(\om)}+\frac{\lambda_{2}^{(n)}}{2}\|f-u-v\|_{L^{2}(\om)}^{2},\quad  \label{L2u_vi}\\
\lambda_{1}\|v_{n}\|_{L^{1}(\om)}&\le J(u)+ \lambda_{1}\|v\|_{L^{1}(\om)}+\frac{\lambda_{2}^{(n)}}{2}\|f-u-v\|_{L^{2}(\om)}^{2}.\label{L1v_vi}
\end{align}
By setting $u=0$, $v=f$ in \eqref{Ju_vi} one obtains a uniform bound for the sequence $(|Du_{n}|(\om))_{n\in
\NN}$. From \eqref{L2u_vi} and from the fact that $\om$ is bounded, one obtains a uniform bound for $(\|f-u_{n}-v_{n}\|_{L^{1}(\om)})_{n\in\NN}$. Similarly, from \eqref{L1v_vi}, a uniform bound on $(\|v_{n}\|_{L^{1}(\om)})_{n\in\NN}$ is obtained and this means that $(\|u_{n}\|_{L^{1}(\om)})_{n\in\NN}$ is also bounded.
Thus every subsequence of $u_{n}$ has a cluster point in $\bv(\om)$ with respect to the weak$^{\ast}$ topology, which must be a minimiser of \eqref{L1J_conv}, \cite[Corollary 7.20]{dalmasogamma}. Furthermore if \eqref{L1J_conv} has a unique minimiser, then every subsequence of $u_{n}$ has a further subsequence that converges to $u$ weakly$^{\ast}$ in $\bv(\om)$. Thus, in this case $u_{n}$ converges to $u$ weakly$^{\ast}$ in $\bv(\om)$ and moreover \eqref{convergence_of_energies_L1} holds by \cite[Corollary 7.20]{dalmasogamma}.

$(ii)$ The $\Gamma$-convergence of the energies follows as above. By setting, $v=u=0$ in \eqref{Ju_vi} one obtains a uniform bound on $(|Du_{n}|(\om))_{n\in\NN}$ and similarly as in $(i)$ a bound on $L^{1}(\om)$ and consequently in $\bv(\om)$ is obtained for $u$. The rest of the proof follows as in $(i)$, bearing in mind that the solution of \eqref{L2J_conv} is unique.

\end{proof}

We summarise our findings so far in Table \ref{asymptotics_table}, where we have combined the results of Proposition \ref{lbl:asymptotics_calatroni_mixed} and Corollary \ref{lbl:u_l2_to_inf}.

{\footnotesize
\setlength\extrarowheight{3pt}
\begin {table}[!h]
\resizebox{\textwidth}{!}{
\begin{tabular}{| c | c | c |c | c | c|}
\cline{2-6}
    \multicolumn{1}{c|}{}
     & \begin{tabular}{@{}c@{}}$\lambda_{1}\to\infty$  \\[-3pt]  $f\in L^{1}(\om)$\\[3pt]\end{tabular}
     & \begin{tabular}{@{}c@{}}$\lambda_{1}\to\infty$  \\[-3pt] $\lambda_{2}$ fixed\\[-3pt] $f\in L^{2}(\om)$\\[3pt]\end{tabular}
     & \begin{tabular}{@{}c@{}}$\lambda_{2}\to\infty$  \\[-3pt] $f\in L^{2}(\om)$\\[-3pt] (or $f\in L^{1}(\om)$ \& $\lambda_{1}$ fixed)\\[3pt] \end{tabular}
     & \begin{tabular}{@{}c@{}}$\lambda_{2}\to\infty$  \\[-3pt] $\lambda_{1}$ fixed\\[-3pt] $f\in L^{1}(\om)$\\[-3pt] $J$--$\lambda_{1}L^{1}$ has !sol. $u^{\ast}$\\[3pt]\end{tabular}
     &\begin{tabular}{@{}c@{}}$\lambda_{1}\to\infty$  \\[-3pt] $\lambda_{2}\to\infty$\\[-3pt] $f\in L^{2}(\om)$\\[3pt] \end{tabular}
     \\
 \cline{2-6} \cline{1-1}
$v$ &$v\to 0$ in $L^{1}(\om)$    & 	$v\to 0$ in $L^{1}(\om)$ & $\|f-u-v\|_{L^{2}(\om)}\to 0$ 	& $v\to f-u^{\ast}$ in $L^{1}(\om)$
& $v\to 0$ in $L^{1}(\om)$ \\[3pt]\hline
$u$ &cannot say   &
\begin{tabular}{@{}c@{}}	$u\to$ solution $J$--$\frac{\lambda_{2}}{2}L^{2}$,\\[-2pt] w$^{\ast}$ in $\bv(\om)$  \end{tabular}& cannot say & 
	$u\to u^{\ast}$,  w$^{\ast}$ in $\bv(\om)$ & 
	\begin{tabular}{@{}c@{}}	$u\to f$ in $L^{1}(\om)$,\\[-2pt] (w$^{\ast}$ in $\bv(\om)$\\[-2pt] if $f\in \bv(\om)$)\end{tabular}\\[3pt]\hline  	
\end{tabular}
}
%& $u\to f$ in $L^{1}(\om)$ 
\vspace{10pt}
\caption{Summary of all the asymptotic results  concerning the solution pair $u, v$ of \eqref{L1L2_J} when one or both parameters $\lambda_{1}$ and $\lambda_{2}$ are let to infinity.}
\label{asymptotics_table}
\end{table}
}

\medskip

In the case of bounded data and TV regularisation, the results obtained above can be refined. We first recall the following well-known result, see \cite[Lemma 3.5]{Zwicknagl}. 

\newtheorem{TVmaximum}[HuberL1L2]{Proposition}
\begin{TVmaximum}\label{lbl:TVmaximum}
Let $u$ be a solution of \eqref{L1L2_TV}, with $f\in L^{\infty}(\om)$. Then the following maximum principle holds:
\begin{equation*}
\essinf_{x\in\om} f(x)\le \essinf_{x\in\om} u(x)\le \esssup_{x\in\om} u(x)\le \esssup_{x\in\om} f(x).
\end{equation*}
\end{TVmaximum}
% \begin{proof}
%  Suppose that \eqref{maximum_princ} does not hold, and let $u$ be the solution of  \eqref{L1L2_TV}. Define the truncation operator $\;\hat{} :\bv(\om)\to L^{\infty}(\om)\cap\bv(\om)$ as
% \[
% \hat{u}(x):= \max\big (\min (\esssup_{t\in\om} f(t),u(x)), \essinf_{t\in\om} f(t)\big ).
% \]
% Then, it follows immediately that $\Phi^{\lambda_{1},\lambda_{2}}(\hat{u},f)< \Phi^{\lambda_{1},\lambda_{2}}(u,f)$. Moreover, it can be checked that $\hat{}$ is 1-Lipschitz continuous. Using now the coarea formula in $\bv(\om)$ \cite{AmbrosioBV}, we get that $|D\hat{u}|(\om)\le |Du|(\om)$, and thus $u$ is suboptimal, which is a contradiction.
% \end{proof}

We can now prove the following result for the $\tv$--$\mathrm{IC}$ minimisation problem \eqref{L1L2_TV}.

\newtheorem{large_ratio}[HuberL1L2]{Proposition}
\begin{large_ratio}\label{lbl:large_ratio}
Suppose that $f\in L^{\infty}(\om)$ and  the parameters $\lambda_{1},\lambda_{2}>0$ satisfy the following condition
\begin{equation}\label{large_ratio_param}
\frac{\lambda_{1}}{\lambda_{2}}\ge 2\|f\|_{\infty}.
\end{equation}
Then, if $u$ is a solution of \eqref{L1L2_TV}, there holds
\begin{equation}\label{Phi_L2}
\Phi^{\lambda_{1},\lambda_{2}}(u,f)=\frac{\lambda_{2}}{2}\|f-u\|_{L^{2}(\om)}^{2}.
\end{equation}
As a result, the problem \eqref{L1L2_TV} is equivalent to a standard $\tv$--$L^{2}$ minimisation problem.
 \end{large_ratio}

\begin{proof}
This is a direct consequence of Proposition \ref{lbl:TVmaximum} and the formulation  \eqref{HuberL1L2a}--\eqref{HuberL1L2b} of $\Phi^{\lambda_{1},\lambda_{2}}(u,f)$. Indeed, using a translation argument and Proposition \ref{lbl:TVmaximum} one shows directly that for the solution $u$ of \eqref{L1L2_TV}, it holds $\|f-u\|_{\infty}\le 2\|f\|_{\infty}$. Thus if \eqref{large_ratio_param} holds, \eqref{HuberL1L2a}  implies that $\varphi(f-u)=\frac{\lambda_{2}}{2}|f-u|^{2}$ so that \eqref{Phi_L2} holds as well.
\end{proof}

We note here that the adaptation of Proposition \ref{lbl:TVmaximum} and, consequently, of Proposition \ref{lbl:large_ratio} to other widely used regularisers is not immediate. For instance, it remains an open problem to show that the solution $u$ of the $\tgv$--$L^{2}$ problem with data $f\in L^{\infty}(\om)$ is also an $L^{\infty}$ function, see for instance the corresponding discussion in \cite{valkonen2014jump2}. However, in dimension one this fact is true  when $f\in \bv(\om)$, by taking advantage of the estimate $\|u\|_{L^{\infty}(\om)}\le C \|u\|_{\bv(\om)}$, see for instance \cite[Proposition 2]{tgv_asymptotic}.

\begin{figure}[!h]
\centering
\begin{tikzpicture}
  \fill [gray!50!white, domain=0:8, variable=\x]
      (0, 0)
      -- plot ({\x}, {0.5*\x})
      -- (0, 4)
      -- cycle;
      \draw [thin] [->] (0,0)--(10,0) node[right, below] {$\lambda_{2}$};
      \draw [thin] [->] (0,0)--(0,4) node[above, left] {$\lambda_{1}$};
      \draw [domain=0:8, variable=\x]
      plot({\x}, {0.5*\x})node[right] at (4.6,2) {$\lambda_{1}=2\|f\|_{\infty}\lambda_{2}$};
      \draw (2,2.5) node[right]{$\tv$--$L^{2}$};
      \draw (9,2) node[right]{$\tv$--$L^{1}$};
      \draw [->] (8,0.5) -- (8.5,0.5);
      \draw [->] (8,1) -- (8.5,1);
      \draw [->] (8,1.5) -- (8.5,1.5);
      \draw [->] (8,2) -- (8.5,2);
      \draw [->] (8,2.5) -- (8.5,2.5);
      \draw [->] (8,3) -- (8.5,3);
      \draw [->] (8,3.5) -- (8.5,3.5);
\end{tikzpicture}
\caption{If $\frac{\lambda_{1}}{\lambda_{2}}\ge 2\|f\|_{\infty}$  then  \eqref{L1L2_TV} is equivalent to the $\tv$--$L^{2}$ problem, see Proposition \ref{lbl:large_ratio}. By fixing $\lambda_{1}$ and sending $\lambda_{2}$ to infinity the solution $u$ converges to a solution of an $\tv$--$L^{1}$ in the sense of Corollary \ref{lbl:u_l2_to_inf}.}
\label{fig:largeratio}
\end{figure}
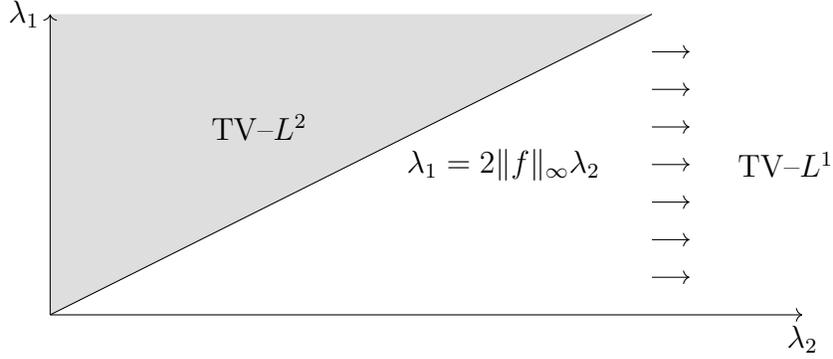

In view of the Proposition above, one sees that in the case of TV regularisation, the Gaussian noise model can be recovered simply  by  fixing either $\lambda_{2}$ and setting $\lambda_{1}$ large enough or by fixing $\lambda_{1}$ and setting $\lambda_{2}$ small enough. In Figure \ref{fig:largeratio} we graphically depict this behaviour.

%\com{About limits to zero: It is clear that if we fix $\lambda_{1}$ and we send $\lambda_{2}$ to zero at some point we will have the mean value of $f$. Could it be the case that if we fix $\lambda_{2}$ and send $\lambda_{1}$ to zero after some point we get the median of $f$? This is not so clear, but some numerics could shed light\ldots}
%\blue{Just to understand: the first claim comes from the expression of $v_{opt}$ right?} \com{that comes from Proposition \ref{lbl:large_ratio} and the fact that in the L2-ROF if we sent the parameter to 0 we get the mean value } \blue{ Basically when we have $TV-L^2$  model we know that we would get the mean value for small values of the parameter $\lambda^2$.  Basically the same for the other way around a part from the ratio $-\frac{\lambda_1^2}{\lambda_2}$, right? Yeah let's check though.}

\subsection{Convergence of the  parameters to zero and non-exact recovery of the data}

We now study the asymptotic behaviour of the model when the parameters are sent to zero. For this analysis and for the sake of simplicity, we focus on the TV minimisation model \eqref{L1L2_TV} but the results can be easily extended to the general regulariser case. First, we recall the definition of the \emph{mean} $u_{\om}$ and \emph{median} values $\left\{u^{\om}\right\}$ of an $L^{1}$ function $u$ defined by:
\begin{align*}
u_{\om}&:=\int_{\om}u\,dx,\\
u^{\om}&\in\underset{c\in\RR}{\operatorname{argmin}}\; \int_{\om}|u-c|\,dx.
\end{align*}

\newtheorem{uniquenessmean}[HuberL1L2]{Remark}

\begin{uniquenessmean}
Note that the median value is not necessarily unique. Moreover, if $u\in L^{2}(\om)$ then $u_{\om}=\underset{c\in\RR}{\operatorname{argmin}}\; \int_{\om}|u-c|^{2}\,dx$.
\end{uniquenessmean}

 We have the following result:

\newtheorem{l1tozero}[HuberL1L2]{Proposition}
\begin{l1tozero}\label{lbl:l1tozero}
Let $f\in L^{1}(\om)$ and $\big(\lambda_{1}^{(n)}\big)_{n\in\NN}$, $\big(\lambda_{2}^{(n)}\big)_{n\in\NN}$ two sequences such that
\begin{equation}\label{l1l2_to_0}
\lambda_{1}^{(n)}\to 0 \quad\text{and}\quad\frac{\lambda_{1}^{(n)}}{\lambda_{2}^{(n)}}\to 0,\quad\text{as }n\to\infty.
\end{equation}
 Then, denoting by $(u_{n})_{n\in\NN}$ the sequence of the corresponding solutions of \eqref{L1L2_TV} we have that 
\[u_{n}\to f^{\om}\quad \text{weakly}^{\ast} \text{ in }\bv(\om).\]
By this we mean that every subsequence of $(u_{n})_{n\in\NN}$ has a further subsequence converging to a median of $f$.
\end{l1tozero}

\begin{proof}
Let $(u_{n})_{n\in\NN}\subseteq \bv(\om)$ be a sequence of the corresponding solutions for the parameters $\big(\lambda_{1}^{(n)},\lambda_{2}^{(n)}\big)$. Notice that the sequence $(u_{n})_{n\in\NN}$ is uniformly bounded in $\bv(\om)$. Indeed the bound on $\tv$ is obtained again from \eqref{Ju_vi}, while for  the $L^{1}$ bound, we first observe that:
\begin{equation}\label{frombelowL1}
-\frac{\lambda_{1}^{2}}{2\lambda_{2}}|\om| +\lambda_{1}\|f-u\|_{L^{1}(\om)}\le\int_{\om}\varphi(f-u)\,dx,\quad \forall u\in L^{1}(\om).
\end{equation}
Using $\eqref{frombelowL1}$, we get
\begin{align*}
-\frac{\big(\lambda_{1}^{(n)}\big)^{2}}{2\lambda_{2}^{(n)}}|\om| +\lambda_{1}^{(n)}\|f-u_{n}\|_{L^{1}(\om)}&\le\int_{\om}\varphi(f-u_{n})\,dx\le \lambda_{1}^{(n)} \|f\|_{L^{1}(\om)} \quad \Rightarrow\\
-\frac{\lambda_{1}^{(n)}}{2\lambda_{2}^{(n)}}|\om| +\|f-u_{n}\|_{L^{1}(\om)}&\le\|f\|_{L^{1}(\om)},
\end{align*}
where the bound is obtained using \eqref{l1l2_to_0}.
Thus every subsequence of $(u_{n})_{n\in\NN}$ has a further (not relabelled) subsequence converging to an element $u^{\ast}\in\bv(\om)$. We will show that $u^{\ast}$ is a median of $f$. Notice first that since $\lambda_{1}^{(n)}\to 0$, then from \eqref{Ju_vi} we get that $|Du_{n}|(\om)\to 0$. Thus, from the lower semicontinuity of total variation and the fact that $\om$ is connected we get that $u^{\ast}$ is a constant. 
%Recall the Poincar\'e inequality for $\bv$ functions,
%\[\|u-u_{\om}\|_{L^{1}(\om)}  \le C\|Du\|(\om), \quad \text{for all }u\in\bv(\om),\]
%where the constant $C>0$ depends only on the domain $\om$.
 By thus setting $u=c\in\RR$ and $v=f-c$ along with \eqref{frombelowL1} in 
 \[|Du_{n}|(\om)+\int_{\om}\varphi(f-u_{n})\,dx \le  |Du|(\om)+ \lambda_{1}^{(n)} \|v\|_{L^{1}(\om)}+\frac{\lambda_{2}^{(n)}}{2} \|f-u-v\|_{L^{2}(\om)}^{2},\]   
we get
\begin{align*}
-\frac{\big (\lambda_{1}^{(n)}\big)^{2}}{2\lambda_{2}^{(n)}} |\om| +\lambda_{1}^{(n)}\|f-u_{n}\|_{L^{1}(\om)} &\le  \lambda_{1}^{(n)} \|f-c\|_{L^{1}(\om)} & \Rightarrow &\\
%-\frac{\lambda_{1}^{(n)}}{\lambda_{2}^{(n)}}|\om| +\|f-u_{n}\|_{L^{1}(\om)} +\frac{1}{\lambda_{1}^{(n)} C}\|u_{n}-(u_{n})_{\om}\|_{L^{1}(\om)}&\le  \|f-c\|_{L^{1}(\om)} & \Rightarrow &\text{ (for small enough $\lambda_{1}^{(n)}$)}\\
%-\frac{\lambda_{1}^{(n)}}{\lambda_{2}^{(n)}}|\om| +\|f-u_{n}\|_{L^{1}(\om)} +\|u_{n}-(u_{n})_{\om}\|_{L^{1}(\om)}&\le  \|f-c\|_{L^{1}(\om)} & \Rightarrow  &\\
-\frac{\lambda_{1}^{(n)}}{\lambda_{2}^{(n)}}|\om|  +\|f-u_{n}\|_{L^{1}(\om)}&\le  \|f-c\|_{L^{1}(\om)} & \Rightarrow & \text{ (by taking limits)}\\
\|f-u^{\ast}\|_{L^{1}(\om)}&\le  \|f-c\|_{L^{1}(\om)}. &  & 
\end{align*}
Hence since $u^{\ast}$ is constant and $c\in\RR$ was arbitrary, we have that $u$ is a median of $f$.
\end{proof}

%\emph{Remark.} Analogous results to Proposition \ref{lbl:l1tozero} can be obtained for other regularisers. For instance concerning (the second order) $\tgv$, note that the following Poincar\'e type of inequality holds, see \cite{BredValk}
%\begin{equation}\label{TGV_poincare}
%\|u-u^{\mathrm{aff}}\|_{L^{1}(\om)}\le C \tgv_{\beta/\alpha,1}^{2}(u),\quad \text{for all }u\in\bv(\om),
%\end{equation}
%where here the constant $C$ depends only on the domain $\om$ and the ratio $\beta/\alpha$ of the $\tgv$ parameters. Moreover, here $u^{\mathrm{aff}}$ is the $L^{2}$-linear regression of $u$, i.e., the affine function which is closer to $u$ in the $L^{2}$ distance. Then by following exactly the same steps as in the proof of Proposition \ref{lbl:l1tozero}, we deduce that when \eqref{l1l2_to_0} holds then the solution of $u$ to the problem \eqref{L1L2_J} with $J(u)=\tgv(u)$ is converging weakly$^{\ast}$ in $\bv(\om)$ to an $L^{1}$-linear regression of the data $f$.

%A closer look to the proof of Proposition \ref{lbl:l1tozero}, actually shows that $\lambda_{2}$ does not have to remain fixed as along as the ratio $\lambda_{1}/\lambda_{2}$ is still going to zero. 
 Similarly, we have the following result: 

\newtheorem{l2tozero}[HuberL1L2]{Proposition}
\begin{l2tozero}\label{lbl:l2tozero}
Let $f\in L^{2}(\om)$ and $\big(\lambda_{1}^{(n)}\big)_{n\in\NN}$, $\big(\lambda_{2}^{(n)}\big)_{n\in\NN}$ two sequences such that
\begin{equation*}
\lambda_{2}^{(n)}\to 0 \quad\text{and}\quad\frac{\lambda_{2}^{(n)}}{\lambda_{1}^{(n)}}\to 0,\quad\text{as }n\to\infty.
\end{equation*}
Then for the corresponding solutions $(u_{n})_{n\in\NN}$ of \eqref{L1L2_TV} we have that 
\[u_{n}\to f_{\om}\quad \text{weakly}^{\ast} \text{ in }\bv(\om).\]
\end{l2tozero}
\begin{proof}
%\com{TODO! check if this can be done if $\lambda_{1}$ is not fixed}
The proof follows the same steps as in Proposition \ref{lbl:l1tozero}. First, observe that  the sequence of solutions $(u_{n})_{n\in\NN}$ is bounded in $\bv(\om)$. Indeed, from \eqref{Ju_vi} we have that $|Du_{n}|(\om)\to 0$. From \eqref{L1v_vi} we further get that $(v_{n})_{n\in\NN}$ is bounded in $L^{1}$ and from \eqref{L2u_vi} we get that $(f-u_{n}-v_{n})_{n\in\NN}$ is bounded in $L^{1}$. Thus, from the triangle inequality we have that $(u_{n})_{n\in\NN}$ is bounded in $L^{1}$. Hence, there exists  a subsequence of $(u_{n})_{n\in\NN}$ that converges to a function $u^{\ast}$ weakly$^{\ast}$ in $\bv(\om)$ with $u^{\ast}$ being a constant.  It remains to show that $u^{\ast}$ is the mean value of $f$. As before, we have for an arbitrary $c\in\RR$
\begin{align}
|Du_{n}|(\om)+\int_{\om}\varphi(f-u_{n})\,dx 
&\le  \frac{\lambda_{2}^{(n)}}{2} \|f-c\|_{L^{2}(\om)}^{2}\quad &\Rightarrow& \nonumber\\
\frac{1}{2}\int_{|f-u_{n}|<\frac{\lambda_{1}^{(n)}}{\lambda_{2}^{(n)}}}  |f-u_{n}|^{2}dx 
& \le \frac{1}{2} \|f-c\|_{L^{2}(\om)}^{2} & \Rightarrow &\;\;(\text{using Fatou's Lemma}) \nonumber\\
\frac{1}{2}\|f-u^{\ast}\|_{L^{2}(\om)}^{2}&\le \frac{1}{2}\|f-c\|_{L^{2}(\om)}^{2}.\nonumber
\end{align}
%Denoting $\alpha_{n}:=\lambda_{1}^{(n)}/\lambda_{2}^{(n)}$ we have that $\alpha_{n}\to\infty$ and 
%\begin{align*}
%\int_{|f-u_{n}|\ge\alpha_{n}} \alpha_{n}|f-u_{n}|-\frac{\alpha_{n}^{2}}{2}\,dx
%&\le \int_{|f-u_{n}|\ge\alpha_{n}} \alpha_{n}|f-u_{n}|\,dx\\
%&\le \int_{|f-u_{n}|\ge\alpha_{n}} |f-u_{n}|^{2}dx \to 0,
%\end{align*}
%where for the last 
Since $u$ is a constant and $c\in\RR$ was arbitrary the proof is complete.
\end{proof}

%\subsection{The one-homogeneous analogue}

\medskip

The following proposition states that with the $L^{1}$--$L^{2}$ infimal convolution fidelity model we can never expect exact recovery of the data. This is similar to the pure $L^{2}$ model, see also \cite[Proposition 4.1]{Zwicknagl}.

\newtheorem{not_exact}[HuberL1L2]{Proposition}
\begin{not_exact}\label{lbl:not_exact}
Let $f\in L^{1}(\om)$ and  $u^{\ast}$ to be a solution of the minimisation problem \eqref{L1L2_TV}. Then $u^{\ast}=f$ if and only if $f$ is a constant.
\end{not_exact}

\begin{proof}
%The proof is rather straightforward. Suppose that $u^{\ast}=f$, then we have
%\begin{align*}
%|Df|(\om)&\le |Du|(\om)+\lambda_{1}\|v\|_{L^{1}(\om)}+\frac{\lambda_{2}}{2}\|f-u-v\|_{L^{2}(\om)}^{2}\quad \text{for all }v\in L^{1}(\om), u\in \bv(\om)\\
%|Df|(\om)&\le |Du|(\om)+ \frac{\lambda_{2}}{2}\|f-u\|_{L^{2}(\om)}^{2}\quad \text{for all }u\in\bv(\om).
%\end{align*}
%This means that $f$ is a minimiser of a $L^{2}$--$\tv$ problem with data $f$, which can happen only if $f$ is a constant, see again  \cite[Proposition 4.1]{Zwicknagl}.

One direction is straightforward. Suppose now that $f$ is a solution of \eqref{L1L2_TV}. Note that in this case necessarily we must have $f\in \bv(\om)\subseteq L^{d^{\ast}}(\om)$, where $d^{\ast}=d/(d-1)$, see \cite{AmbrosioBV}. It follows that for every $0<\epsilon<1$, the function $f_{\epsilon}:=\epsilon f$ is  suboptimal. Thus we have
\begin{equation*}
|Df|(\om)\le \epsilon |Df|(\om)+\Phi^{\lambda_{1},\lambda_{2}}(u,\epsilon f) \quad \Longrightarrow \quad 0\le (\epsilon-1) |Df|(\om)+ \Phi^{\lambda_{1},\lambda_{2}}(u,\epsilon f) .
\end{equation*}
We continue
\begin{align}
0&\le (\epsilon-1) |Df|(\om)+ \Phi^{\lambda_{1},\lambda_{2}}(u,\epsilon f) \nonumber\\
0&\le (\epsilon-1) |Df|(\om)+ \int_{|f-\epsilon f|\ge \frac{\lambda_{1}}{\lambda_{2}}} \lambda_{1} |f-\epsilon f|-\frac{\lambda_{1}^{2}}{2\lambda_{2}}\, dx+\int_{|f-\epsilon f| < \frac{\lambda_{1}}{\lambda_{2}}} \frac{\lambda_{2}}{2} |f-\epsilon f|^{2}\, dx+ \;\Longrightarrow\nonumber\\
0&\ge|Df|(\om) - \int_{|f|\ge \frac{\lambda_{1}}{\lambda_{2}|1-\epsilon|}} \lambda_{1} |f|\, dx  + \int_{|f|\ge \frac{\lambda_{1}}{\lambda_{2}|1-\epsilon|}} \frac{\lambda_{1}^{2}}{2\lambda_{2}|1-\epsilon|} \, dx
- |1-\epsilon| \int_{|f|< \frac{\lambda_{1}}{\lambda_{2}|1-\epsilon|}} \frac{\lambda_{2}}{2} |f|^{2}\, dx.\label{3terms}
\end{align}
Now working with each one of the first three terms in \eqref{3terms},  using dominated convergence, we have
\begin{align}
\lim_{\epsilon\to 1} \int_{|f|\ge \frac{\lambda_{1}}{\lambda_{2}|1-\epsilon|}} |f|\, dx=0,\label{term1}
\end{align}
\begin{align}
\lim_{\epsilon\to 1} \int_{|f|\ge \frac{\lambda_{1}}{\lambda_{2}|1-\epsilon|}} \frac{\lambda_{1}^{2}}{2\lambda_{2}|1-\epsilon|} \, dx= \frac{\lambda_{1}}{2} \lim_{t\to \infty}  \int_{|f|\ge t} t\, dx
\le  \frac{\lambda_{1}}{2} \lim_{t\to \infty}   \int_{|f|\ge t} |f|\, dx
=0, \label{term2}
\end{align}
and finally
\begin{align}
& \lim_{\epsilon\to 1} |1-\epsilon| \int_{|f|< \frac{\lambda_{1}}{\lambda_{2}|1-\epsilon|}} \frac{\lambda_{2}}{2} |f|^{2}\, dx
= \frac{\lambda_{2}}{2}\lim_{\epsilon\to 1} |1-\epsilon| \int_{|f|< \frac{\lambda_{1}}{\lambda_{2}|1-\epsilon|}}  |f|^{\frac{d}{d-1}} |f|^{\frac{d-2}{d-1}}\nonumber\\
& \le  \frac{\lambda_{2}}{2}  \left (\frac{\lambda_{1}}{\lambda_{2}} \right)^{\frac{d-2}{d-1}} \lim_{\epsilon\to 1} 
\frac{|1-\epsilon|}{|1-\epsilon|^{\frac{d-2}{d-1}}} \int_{|f|< \frac{\lambda_{1}}{\lambda_{2}|1-\epsilon|}}   |f|^{\frac{d}{d-1}} \, dx
\le  \frac{\lambda_{2}}{2}  \left (\frac{\lambda_{1}}{\lambda_{2}} \right)^{\frac{d-2}{d-1}} \|f\|_{L^{d^{\ast}}(\om)}^{d^{\ast}} \lim_{\epsilon\to 1} |1-\epsilon|^{\frac{1}{d-1}}=0 \label{term3}.
\end{align}
By combining \eqref{term1}, \eqref{term2}, \eqref{term3} with \eqref{3terms} we get that $|Df|(\om)=0$, and since $\om$ is connected, $f$ is a constant function. Note that the calculations above assumed that $d>1$ but if $d=1$ then $f\in L^{\infty}(\om)\subseteq L^{2}(\om)$ and the analogous calculation to \eqref{term3} follows more easily.
\end{proof}

It is clear from all the analysis above that the $\tv$--$\mathrm{IC}$ model, at least when $f\in L^{\infty}(\om)$, can reproduce the $\tv$--$L^{2}$ solutions, but as far as the $\tv$--$L^{1}$ solutions are concerned, these are only (guaranteed to be) recovered  in the limit $\lambda_{2}\to\infty$, see again Figure \ref{fig:largeratio}. 

\subsection{The one-homogeneous analogue}

We conclude this section by briefly presenting an alternative form of \eqref{L1L2_fid_def}, i.e., its one-homogeneous analogue, by which the $\tv$--$L^{1}$ solutions can also be recovered for finite parameters. This discussion is motivated by some analogous results in \cite{journal_tvlp}.
We  define:
\begin{equation*}
\Phi_{1-hom}^{\lambda_{1},\lambda_{2}}(u,f):=\min_{v\in L^{1}(\om)} \lambda_{1} \|v\|_{L^{1}(\om)} +\lambda_{2} \|f-u-v\|_{L^{2}(\om)}, 
\end{equation*}
which, via a straightforward computation gives 
\begin{align*}
& \min\left(\lambda_{1},\frac{\lambda_{2}}{2|\om|^{1/2}}\right)\|f-u\|_{L^{1}(\om)}
\le \min_{v\in L^{1}(\om)} \lambda_{1} \|v\|_{L^{1}(\om)}+\frac{\lambda_{2}}{2|\om|^{1/2}}\|f-u-v\|_{L^{1}(\om)}\\
&\le\min_{v\in L^{1}(\om)}  \lambda_{1} \|v\|_{L^{1}(\om)}+\frac{\lambda_{2}}{2}\|f-u-v\|_{L^{2}(\om)} \le \lambda_{1}\|f-u\|_{L^{1}(\om)}.
\end{align*} 
Hence, it is clear that 
\[\Phi_{1-hom}^{\lambda_{1},\lambda_{2}}(u,f)=\lambda_{1}\|f-u\|_{L^{1}(\om)}, \quad \text{if} \quad \frac{\lambda_{1}}{\lambda_{2}}\le \frac{1}{2|\om|^{1/2}}, \]
Thus, under such choice the $\tv$--$L^{1}$ model can be recovered. We state here without any proof the relationship between these two versions. Let us define the following sets for $f\in L^{1}(\om)$:
\begin{align*}
S_{IC}&= \Big \{u^{\ast}\in \bv(\om):\; u^{\ast}= \underset{u\in \bv(\om)}{\operatorname{argmin}}\; |Du|(\om)+\Phi^{\lambda_{1},\lambda_{2}}(u,f),\;\text{ for some }\lambda_{1},\lambda_{2}>0\Big \},\\
S_{IC}^{1-hom}&= \Big \{u^{\ast}\in \bv(\om):\; u^{\ast}= \underset{u\in \bv(\om)}{\operatorname{argmin}}\; |Du|(\om)+\Phi_{1-hom}^{\lambda_{1},\lambda_{2}}(u,f),\;\text{ for some }\lambda_{1},\lambda_{2}>0\Big \},\\
S_{L^{1}}&= \Big \{u^{\ast}\in \bv(\om):\; u^{\ast}= \underset{u\in \bv(\om)}{\operatorname{argmin}}\; |Du|(\om)+\lambda_{1}\|f-u\|_{L^{1}(\om)},\;\text{ for some }\lambda_{1}>0\Big \}.
\end{align*} 
Then one can show by using similar techniques as in \cite{journal_tvlp} that
\begin{equation*}
S_{IC}^{1-hom}=S_{IC}\cup S_{L^{1}}\quad\text{and that, in general,}\quad S_{L^{1}}\setminus S_{IC}\ne \emptyset. 
\end{equation*}

\section{Exact solutions}   \label{sec:exact}
%\blue{OK so this answers my question above. We should structure the paper in two parts where we have analysis and numerics for each section} \com{yeah as I said let's see what we have in end and then we decide. I just wrote these section headers more as a list of what we could do :)}

In order to get more insights about the relationship of the $\tv$--$\mathrm{IC}$ model with the pure $\tv$--$L^{1}$ and $\tv$--$L^{2}$ models, we compute in this section some exact solutions for simple one dimensional data functions $f$. In particular, we set here $\om=(-2L,2L)$ for some $L>0$, and we consider as data $f$ the following step function
\begin{equation}\label{data_f}
f(x)=
\begin{cases}
0, & \text{ if}\quad x\in (-2L,-L)\cup (L,2L),\\
h, & \text{ if}\quad x\in [-L,L],
\end{cases}
\end{equation}
where $h>0$. Using similar techniques as in \cite{BrediesL1, journal_tvlp, Papafitsoros_Bredies, ring2000structural}, we can easily show using primal-dual optimality conditions, that a function $u\in\bv(\om)$ is a solution of \eqref{L1L2_TV} if and only if there exists a function $v\in H_{0}^{1}(\om)$ such that
\begin{align}
v'&=
\begin{cases}
\lambda_{1} \frac{f-u}{|f-u|}, & \text{ if}\quad |f-u|\ge \frac{\lambda_{1}}{\lambda_{2}}, \\
\lambda_{2} (f-u), & \text{ if}\quad |f-u|<\frac{\lambda_{1}}{\lambda_{2}},
\end{cases} \label{opt1}\\
v&\in \mathrm{Sgn}(Du),\label{opt2}
\end{align}
where 
\[\mathrm{Sgn}(Du)= \left \{ v\in L^{\infty}(\om)\cap L^{\infty}(\om,Du):\; \|v\|_{\infty}\le 1, \; v=\frac{dDu}{d|Du|},\; |Du|\text{--a.e.} \right \}.\]
Here $\frac{dDu}{d|Du|}$ denotes the Radon--Nikod\'ym density of $Du$ with respect to $|Du|$. 
Compared to the aforementioned references, the only difference here is the right-hand side of \eqref{opt1} which is the subdifferential of $\Phi^{\lambda_{1},\lambda_{2}}$ evaluated at $f-u$. With the help of the optimality conditions above, we are able to compute analytically all the solutions to the problem \eqref{L1L2_TV} for the data \eqref{data_f}, and for all combinations of the parameters $\lambda_{1}, \lambda_{2}$. Note that similarly to pure $L^{1}$ and $L^{2}$ cases one can show that no new jump discontinuities are created for the solution $u$, which will be constant in the areas where $\{f\ne u\}$. Thus all the solutions will be either constants or piecewise constants with jumps at $x=-L$ and $x=L$ which must also have the same orientations with the jumps of $f$.

We first examine the case $u=\frac{h}{2}$, i.e., the mean value of $f$. In such case we have $|f-u|=\frac{h}{2}$ everywhere and thus if $\frac{h}{2}\le \frac{\lambda_{1}}{\lambda_{2}}$ then $v'=\lambda_{2}(f-u)$ everywhere. In order for the condition \eqref{opt1} to hold we must also have $\lambda_{2}\le \frac{2}{hL}$. Note that if $\frac{h}{2}\le\frac{\lambda_{1}}{\lambda_{2}}$ and $\lambda_{2}< \frac{2}{hL}$, then the solution must be  constant otherwise one can check that in every case the condition \eqref{opt2} would be violated. One can further show in this case that if $u$ is a constant with $\frac{h}{2}<\frac{\lambda_{1}}{\lambda_{2}}$, the only possibility is $u=\frac{h}{2}$. 

Suppose now that $\frac{h}{2}\ge \frac{\lambda_{1}}{\lambda_{2}}$ and also $\lambda_{1}<\frac{1}{L}$. Observe that in this case, every constant function $u=c$ with $\frac{\lambda_{1}}{\lambda_{2}}\le c \le h-\frac{\lambda_{1}}{\lambda_{2}}$ satisfies \eqref{opt1}--\eqref{opt2}. In that case we have $|f-u|\ge \frac{\lambda_{1}}{\lambda_{2}}$ everywhere. Notice again that no other constant function is a solution. By contradiction, that would  mean that $|f-u|>\frac{\lambda_{1}}{\lambda_{2}}$ and $|f-u|<\frac{\lambda_{1}}{\lambda_{2}}$ on $(-2L,L)\cup (L,2L)$ and $(-L,L)$ respectively (or vice versa). With the help of \eqref{opt1} and the fact that  $v\in H_{0}^{1}(\om)$ one would then arrive to  a contradiction.  Furthermore, one can  check that discontinuous solutions cannot occur in this case either.

We concentrate now on the case $\frac{h}{2}\ge \frac{\lambda_{1}}{\lambda_{2}}$ and $\lambda_{1}=\frac{1}{L}$. Note that the constant functions $u=c$ with $\frac{\lambda_{1}}{\lambda_{2}}\le c \le h-\frac{\lambda_{1}}{\lambda_{2}}$ are  solutions in this case as well. However, one can also verify that the following family of discontinuous functions are also solutions:
\[u(x)=
\begin{cases}
c_{1}, & \text{ if}\quad x\in (-2L,-L),\\
h-d,    & \text{ if}\quad x\in [-L,L],\\
c_{2}, & \text{ if}\quad x\in (L,2L),
\end{cases}
\qquad \frac{\lambda_{1}}{\lambda_{2}} \le c_{i}<h-d\le h-\frac{\lambda_{1}}{\lambda_{2}}, \;\; i=1,2.
\]
It can be checked similarly as before that no other solutions can occur.

Finally, we consider the case $\lambda_{1}>\frac{1}{L}$ and $\lambda_{2}\ge \frac{2}{hL}$. We claim that in that case the unique solution is given by
\[u(x)=
\begin{cases}
\frac{1}{L\lambda_{2}}, & \text{ if}\quad x\in (-2L,-L)\cup (L,2L),\\
h-\frac{1}{L\lambda_{2}},    & \text{ if}\quad x\in [-L,L].\\
\end{cases}
\]
One can similarly check that no other solution is possible. We summarise our findings in the following proposition.

\newtheorem{prop_exact}[HuberL1L2]{Proposition}
\begin{prop_exact}\label{lbl:prop_exact}
Let $\om=(-2L,2L)$ and $f\in\bv(\om)$ being the jump function given by \eqref{data_f}. Then the solutions $u$ to the $\tv$--$\mathrm{IC}$ minimisation problem
\[\min_{u\in\bv(\om)} |Du|(\om)+\Phi^{\lambda_{1},\lambda_{2}}(u,f),\]
are given by the following formulae:
\begin{enumerate}
\item If $\frac{h}{2}<\frac{\lambda_{1}}{\lambda_{2}}$ and $\lambda_{2}\le \frac{2}{hL}$, then the solution is unique and given by
\[u=\frac{h}{2}.\]
\item If $\frac{h}{2}\ge \frac{\lambda_{1}}{\lambda_{2}}$ and $\lambda_{1}<\frac{1}{L}$, then there exist infinitely many constant solutions given by
\[u=c,\qquad \frac{\lambda_{1}}{\lambda_{2}}\le c\le h -\frac{\lambda_{1}}{\lambda_{2}}.\]
\item If $\frac{h}{2}\ge \frac{\lambda_{1}}{\lambda_{2}}$ and $\lambda_{1}=\frac{1}{L}$, then there there exist infinitely many solutions given by
\[u(x)=
\begin{cases}
c_{1}, & \text{ if}\quad x\in (-2L,-L),\\
h-d,    & \text{ if}\quad x\in [-L,L],\\
c_{2}, & \text{ if}\quad x\in (L,2L),
\end{cases}
\qquad \frac{\lambda_{1}}{\lambda_{2}} \le c_{i}\le h-d\le h-\frac{\lambda_{1}}{\lambda_{2}}, \;\; i=1,2.
\]
\item If $\lambda_{2}>\frac{2}{hL}$ and $\lambda_{1}>\frac{1}{L}$, then the solution is unique and given by
\[u(x)=
\begin{cases}
\frac{1}{L\lambda_{2}}, & \text{ if}\quad x\in (-2L,-L)\cup (L,2L),\\
h-\frac{1}{L\lambda_{2}},    & \text{ if}\quad x\in [-L,L].\\
\end{cases}
\]
\end{enumerate}
\end{prop_exact}

%\begin{figure}[t]
%\includegraphics[width=0.8\textwidth]{images/exact.eps}
%\caption{Visualisation of all the possible solutions to the $L^{1}$--$L^{2}$ infimal convolution fidelity total variation problem \eqref{L1L2_TV} for data \eqref{data_f} and for all the possible combinations of the parameters $\lambda_{1}$ and $\lambda_{2}$, see Proposition \ref{lbl:prop_exact}.}
%\label{fig:exact}
%\end{figure}

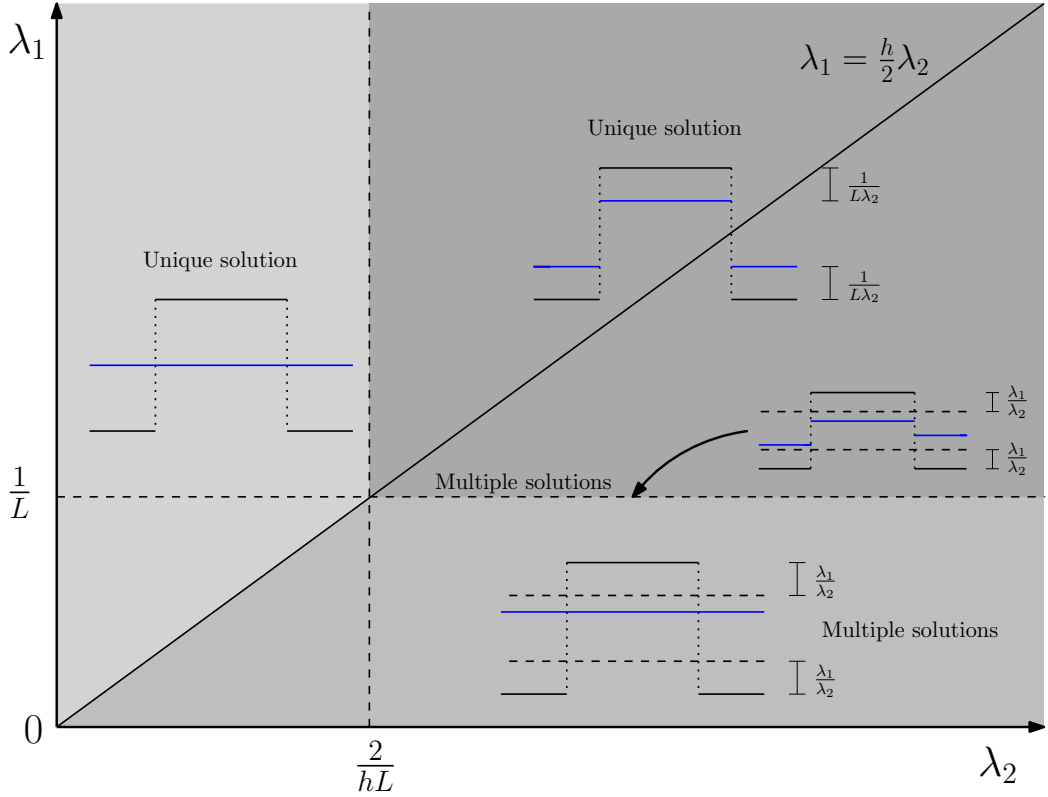
\begin{figure}[t!]
\centering
\resizebox{1\textwidth}{!}{
\begin{tikzpicture}[ipe stylesheet]
  \useasboundingbox (0, 390) rectangle (595, 842);
  \draw[line cap=round]
    (96, 544)
     -- (96, 544);
  \node[ipe node, font=\LARGE]
     at (40, 744) {$\lambda_{1}$};
  \node[ipe node, font=\LARGE]
     at (512, 392) {$\lambda_{2}$};
  \draw[line cap=round]
    (216, 768)
     -- (216, 768);
  \node[ipe node, font=\LARGE]
     at (48, 408) {$0$};
  \node[ipe node, font=\LARGE]
     at (209, 392) {$\frac{2}{hL}$};
  \node[ipe node, font=\LARGE]
     at (40, 525) {$\frac{1}{L}$};
  \fill[lightgray]
    (64, 416)
     -- (216, 528)
     -- (216, 768)
     -- (64, 768)
     -- cycle;
  \fill[gray]
    (64, 416)
     -- (216, 528)
     -- (544, 528)
     -- (544, 416)
     -- cycle;
  \fill[darkgray]
    (216, 528) rectangle (544, 768);
  \draw[line cap=round]
    (544, 528)
     -- (544, 528)
     -- cycle;
  \draw[ipe pen fat, ->]
    (64, 416)
     -- (64, 416)
     -- (544, 416);
  \draw[ipe pen fat, <-]
    (64, 768)
     -- (64, 416);
  \draw[ipe pen heavier, ipe dash dashed]
    (216, 416)
     -- (216, 768);
  \draw[ipe pen heavier, ipe dash dashed]
    (64, 528)
     -- (544, 528);
  \draw[ipe pen heavier]
    (64, 416)
     -- (544, 768);
  \node[ipe node, font=\Large]
     at (425, 736) {$\lambda_{1}=\frac{h}{2}\lambda_{2}$};
  \node[ipe node, font=\large]
     at (106, 640) {\small{Unique solution}};
  \draw[ipe pen heavier]
    (80, 560)
     -- (112, 560)
     -- (112, 560);
  \draw[ipe pen heavier]
    (112, 624)
     -- (176, 624)
     -- (176, 624);
  \draw[ipe pen heavier]
    (176, 560)
     -- (208, 560)
     -- (208, 560);
  \draw[ipe pen heavier, ipe dash dotted]
    (112, 624)
     -- (112, 560);
  \draw[ipe pen heavier, ipe dash dotted]
    (176, 624)
     -- (176, 560);
  \draw[blue, ipe pen heavier]
    (208, 592)
     -- (80, 592);
  \draw[ipe pen heavier]
    (280, 432)
     -- (312, 432)
     -- (312, 432);
  \draw[ipe pen heavier]
    (312, 496)
     -- (376, 496)
     -- (376, 496);
  \draw[ipe pen heavier]
    (376, 432)
     -- (408, 432)
     -- (408, 432);
  \draw[ipe pen heavier, ipe dash dotted]
    (312, 496)
     -- (312, 432);
  \draw[ipe pen heavier, ipe dash dotted]
    (376, 496)
     -- (376, 432);
  \draw[ipe pen heavier, ipe dash dashed]
    (408, 480)
     -- (280, 480)
     -- (280, 480);
  \draw[ipe pen heavier, ipe dash dashed]
    (408, 448)
     -- (280, 448)
     -- (280, 448);
  \draw[blue, ipe pen heavier]
    (408, 472)
     -- (280, 472);
  \draw
    (424, 480)
     -- (424, 496)
     -- (424, 496);
  \draw
    (420, 496)
     -- (428, 496)
     -- (428, 496);
  \draw
    (420, 480)
     -- (428, 480)
     -- (428, 480);
  \node[ipe node]
     at (432, 484) {$\frac{\lambda_{1}}{\lambda_{2}}$};
  \draw
    (424, 432)
     -- (424, 448)
     -- (424, 448);
  \draw
    (420, 448)
     -- (428, 448)
     -- (428, 448);
  \draw
    (420, 432)
     -- (428, 432)
     -- (428, 432);
  \node[ipe node]
     at (432, 436) {$\frac{\lambda_{1}}{\lambda_{2}}$};
  \draw[shift={(405.345, 541.658)}, xscale=0.7889, yscale=0.5792, ipe pen heavier]
    (0, 0)
     -- (32, 0)
     -- (32, 0);
  \draw[shift={(430.59, 578.73)}, xscale=0.7889, yscale=0.5792, ipe pen heavier]
    (0, 0)
     -- (64, 0)
     -- (64, 0);
  \draw[shift={(481.08, 541.658)}, xscale=0.7889, yscale=0.5792, ipe pen heavier]
    (0, 0)
     -- (32, 0)
     -- (32, 0);
  \draw[shift={(430.59, 578.73)}, xscale=0.7889, yscale=0.5792, ipe pen heavier, ipe dash dotted]
    (0, 0)
     -- (0, -64);
  \draw[shift={(481.08, 578.73)}, xscale=0.7889, yscale=0.5792, ipe pen heavier, ipe dash dotted]
    (0, 0)
     -- (0, -64);
  \draw[shift={(506.325, 569.462)}, xscale=0.7889, yscale=0.5792, ipe pen heavier, ipe dash dashed]
    (0, 0)
     -- (-128, 0)
     -- (-128, 0);
  \draw[shift={(506.325, 550.927)}, xscale=0.7889, yscale=0.5792, ipe pen heavier, ipe dash dashed]
    (0, 0)
     -- (-128, 0)
     -- (-128, 0);
  \draw[shift={(518.947, 569.462)}, xscale=0.7889, yscale=0.5792]
    (0, 0)
     -- (0, 16)
     -- (0, 16);
  \draw[shift={(515.792, 578.73)}, xscale=0.7889, yscale=0.5792]
    (0, 0)
     -- (8, 0)
     -- (8, 0);
  \draw[shift={(515.792, 569.462)}, xscale=0.7889, yscale=0.5792]
    (0, 0)
     -- (8, 0)
     -- (8, 0);
  \node[ipe node]
     at (525.258, 571.779) {$\frac{\lambda_{1}}{\lambda_{2}}$};
  \draw[shift={(518.947, 541.659)}, xscale=0.7889, yscale=0.5792]
    (0, 0)
     -- (0, 16)
     -- (0, 16);
  \draw[shift={(515.792, 550.927)}, xscale=0.7889, yscale=0.5792]
    (0, 0)
     -- (8, 0)
     -- (8, 0);
  \draw[shift={(515.792, 541.659)}, xscale=0.7889, yscale=0.5792]
    (0, 0)
     -- (8, 0)
     -- (8, 0);
  \node[ipe node]
     at (525.258, 543.976) {$\frac{\lambda_{1}}{\lambda_{2}}$};
  \draw[shift={(430.59, 564.828)}, xscale=0.7889, yscale=0.5792, blue, ipe pen heavier]
    (0, 0)
     -- (64, 0)
     -- (64, 0);
  \draw[shift={(481.08, 557.877)}, xscale=0.7889, yscale=0.5792, blue, ipe pen heavier]
    (0, 0)
     -- (32, 0)
     -- (28, 0);
  \draw[shift={(405.346, 553.243)}, xscale=0.7889, yscale=0.5792, blue, ipe pen heavier]
    (0, 0)
     -- (32, 0)
     -- (28, 0);
  \draw[ipe pen fat, <-]
    (344, 528)
     .. controls (360, 546.6667) and (378.6667, 557.3333) .. (400, 560);
  \node[ipe node, font=\large]
     at (436, 460) {\small{Multiple solutions}};
  \node[ipe node]
     at (364, 568) {};
  \node[ipe node, font=\large]
     at (248, 532) {\small{Multiple solutions}};
  \draw[ipe pen heavier]
    (296, 624)
     -- (328, 624)
     -- (328, 624);
  \draw[ipe pen heavier]
    (328, 688)
     -- (392, 688)
     -- (392, 688);
  \draw[ipe pen heavier]
    (392, 624)
     -- (424, 624)
     -- (424, 624);
  \draw[ipe pen heavier, ipe dash dotted]
    (328, 688)
     -- (328, 624);
  \draw[ipe pen heavier, ipe dash dotted]
    (392, 688)
     -- (392, 624);
  \draw[blue, ipe pen heavier]
    (328, 672)
     -- (392, 672)
     -- (392, 672);
  \draw[blue, ipe pen heavier]
    (424, 640)
     -- (392, 640)
     -- (392, 640);
  \draw[blue, ipe pen heavier]
    (328, 640)
     -- (296, 640)
     -- (304, 640);
  \draw
    (440, 672)
     -- (440, 688)
     -- (440, 688);
  \draw
    (436, 688)
     -- (444, 688)
     -- (444, 688);
  \draw
    (436, 672)
     -- (444, 672)
     -- (444, 672);
  \node[ipe node]
     at (448, 676) {$\frac{1}{L\lambda_{2}}$};
  \draw
    (440, 624)
     -- (440, 640)
     -- (440, 640);
  \draw
    (436, 640)
     -- (444, 640)
     -- (444, 640);
  \draw
    (436, 624)
     -- (444, 624)
     -- (444, 624);
  \node[ipe node]
     at (448, 628) {$\frac{1}{L\lambda_{2}}$};
  \node[ipe node, font=\large]
     at (322, 704) {\small{Unique solution}};
\end{tikzpicture}
}
\caption{Visualisation of all the possible solutions to the $\tv$--$\mathrm{IC}$ problem \eqref{L1L2_TV} for data \eqref{data_f} and for all the possible combinations of the parameters $\lambda_{1}$ and $\lambda_{2}$, see Proposition \ref{lbl:prop_exact}.}
\label{fig:exact}
\end{figure}

A visualisation of these solutions is depicted in Figure \ref{fig:exact}. Observe that for large enough ratio $\frac{\lambda_{1}}{\lambda_{2}}$ all the $\mathrm{TV}$--$L^{2}$ solutions are recovered as Proposition \ref{lbl:large_ratio}  predicts -- compare also  Figures \ref{fig:largeratio} and \ref{fig:exact}. Moreover, observe that as $\lambda_{1}$ and $\frac{\lambda_{1}}{\lambda_{2}}$ goes to zero, the solutions indeed converge to a median of $f$, as shown in Proposition \ref{lbl:l1tozero}. Note however, that apart from some medians,  the solutions of the  $\mathrm{TV}$--$L^{1}$ model are not recovered. More precisely, the ones that perfectly fit the data in the whole domain or part of it, i.e., $u=h$, $u=0$ and $u=f$ cannot be obtained here. This is also in accordance to Proposition \ref{lbl:not_exact}.

%\subsection{Characteristic function of a disk in a plane}

%subsection{One dimensional examples}

%\section{Rule of thumb for choosing $\lambda_{1}$, $\lambda_{2}$}

%\section{$\Phi^{\lambda_{1},\lambda_{2}}$ and higher order regularisers}

\section{Automatic selection of parameters} \label{sec:optimisation}

We describe now a \emph{bilevel optimisation} strategy for the estimation of optimal parameters $\lambda_1$ and $\lambda_2$ in the $\tv$--$\mathrm{IC}$ model \eqref{L1L2_TV} based on the use of training sets, see \cite{bilevellearning,noiselearning,DeLosReyes2017}. This approach has been heuristically considered for the $\tv$--$\mathrm{IC}$  model in \cite[Section 7]{calatroni_mixed} with little theoretical justification. To fill this gap, we  prove in this section existence results for the solution of the bilevel minimisation problem and for the corresponding adjoint problem, thus making the derivation the optimality systemin \cite{calatroni_mixed} rigorous.
%Finally this section is furnished with a series of numerical examples. 

We point out that in \cite{Langer2017, langerl1l2} an adaptive optimisation approach has been proposed for the automatic selection of parameters when a linear combination of $L^1$ and $L^2$ data fidelities is considered. However, differently to our setting, in that approach the noise level is assumed to be known. 
%On the other hand, our approach does not require such prior knowledge, and it is rather suitable for set-ups that are described next.

\subsection{Bilevel optimisation}   \label{sec:bilevel}

Learning approaches have become very popular over the recent years due to their ability of combining data- and model-driven algorithms for the optimal design of imaging models. In particular, bilevel optimisation techniques have been proposed by several different authors in discrete \cite{pockSVM, pockbilevel,pockBilevelNonsmooth,pockReactionDiffusion} and functional \cite{bilevellearning, noiselearning, interiorpaper2015, DeLosReyes2017, hintermuellerPartI, hintermuellerPartII} settings as a tool to estimate the ``best" variational image restoration model within a certain class by means of training examples. Typically, such examples consist of  images obtained in standard acquisition settings, and thus corrupted by noise with equal (unknown) intensity, paired with their corresponding versions ideally acquired in a very low-noise setting. In medical imaging, for instance, such training set can be provided by means of real and/or simulated phantoms. Note that in order to make the estimation robust, a large training set is often desirable; for that, stochastic optimisation techniques and sampling approaches can be used to reduce the computational costs, see, e.g., \cite{calatroni2014dynamic}. 

%In the following we stick for simplicity with the TV regularisation and consider the problem of estimating the optimal parameters $\lambda_1$ and $\lamdba_2$ in \eqref{L1L2_TV}. 

The general bilevel optimisation problem can be formulated as:
\begin{equation}   \label{bilevel_general1}
\min_{\bm{\lambda}\in [0,\infty)^m}~F(u_{\bm{\lambda}})
\end{equation}
subject to:
\begin{equation}  \label{bilevel_general2}
u_{\bm{\lambda}}\in\argmin_{u\in X} \left\{ J(u,\bm{\lambda}):= |Du|(\Omega)+ \Phi^{\bm{\lambda}}(u,f)   \right\},
\end{equation}
where $\bm{\lambda}=(\lambda_1,\ldots,\lambda_m)\in [0,\infty)^m$ are the parameters to optimise and $F\ge 0$ is an appropriate quality measure which is minimised under the constraint that the function $u_{\bm{\lambda}}$ is a solution of the denoising model \eqref{bilevel_general2} in a suitable function space $X$.

In \cite{bilevellearning, noiselearning} this approach is used to estimate the optimal parameters $\bm{\lambda}$ in the case when single data models $\Phi_i$, $i=1,\ldots,m$ are linearly combined, i.e., when
$$
\Phi^{\bm{\lambda}}(u,f) = \sum_{i=1}^m ~\lambda_i\Phi_i(u,f).
$$
\noindent There, theoretical results showing existence of minima and the adjoint state for the problem \eqref{bilevel_general1}--\eqref{bilevel_general2} are  shown and Newton-type methods are proposed for its efficient numerical solution. Similar results and algorithms are further studied in \cite{interiorpaper2015, DeLosReyes2017} for the estimation of optimal parameters of higher-order regularisers (e.g., TGV) combined with Gaussian fidelity.

\medskip 
In the following we set $m=2$ and
%we work for simplicity with TV regularisation and 
 consider the problem of estimating the optimal parameters $\lambda_1$ and $\lambda_2$ in \eqref{L1L2_TV}.

\paragraph{General framework:}
%For a proper, convex and BV weak$^*$ lower semicontinuous 
The non-smooth $\tv$--$\mathrm{IC}$ bilevel problem reads:
\begin{equation}\label{bilevel_unreg}
\begin{aligned}
&\min_{\lambda_{1},\lambda_{2}~\ge 0}~ F(u_{\lambda_{1},\lambda_{2}})\\
\text{subject to }\quad u_{\lambda_{1},\lambda_{2}}\in~ &\underset{u\in\bv(\om)}{\operatorname{argmin}}~ |Du|(\om)+ \Phi^{\lambda_{1},\lambda_{2}}(u,f),
\end{aligned}
\end{equation}
where $\Phi^{\lambda_1,\lambda_2}$ is the IC fidelity \eqref{L1L2_fid_def}.
%We will also require the following two assumptions for the cost functional $F$:
%\begin{align}
%&F(u_{\lambda_{1},\lambda_{2}})<\infty \text{ for some } \lambda_{1}, \lambda_{2}\ge %0. \label{assumption_F_1}\\
%&\underset{c\in \RR}{\operatorname{argmin}}~ F(c)\ne \emptyset.\label{assumption_F_2}
%\end{align}
%Note that if $\lambda_{1}=0$ and/or $\lambda_{2}=0$, then $\Phi^{\lambda_{1},\lambda_{2}}(u,f)=0$ for every $u\in L^{1}(\om)$. In that case every constant function is a minimiser of the lower level problem of \eqref{bilevel_unreg} and the condition \eqref{assumption_F_1} is satisfied if merely $F(c)<\infty$ for some constant $c\in\RR$.
We follow \cite{bilevellearning,noiselearning, interiorpaper2015, DeLosReyes2017} and introduce an appropriate smoothing of the TV semi-norm combined with a further quadratic smoothing. This is crucial for the following proofs and for the design of the gradient-based optimisation algorithms we intend to use.

For $\epsilon\ll 1$, we then consider the following regularised version of \eqref{bilevel_unreg}:
\begin{equation}\label{bilevel_reg}
\begin{aligned}
&\min_{\lambda_{1},\lambda_{2}\ge 0} F(u_{\lambda_{1},\lambda_{2}})\\
\text{subject to }\quad u_{\lambda_{1},\lambda_{2}}\in &~\underset{u\in H^{1}(\om)}{\operatorname{argmin}}~\frac{\epsilon}{2}\|u\|_{H^{1}(\om)}^{2}+ \|\nabla u\|_{\gamma, L^{1}(\om)} + \Phi_{\epsilon,\gamma}^{\lambda_{1},\lambda_{2}}(u,f)  .
\end{aligned}
\end{equation}
Here, we denote by $\|\nabla u\|_{\gamma, L^{1}(\om)}=\int_{\om} |\nabla u|_\gamma \,dx$, with $|\cdot|_\gamma$ being a smooth Huber-type regularisation depending on a parameter $\gamma >0$ whose $C^1$-derivative reads:
\begin{equation}\label{eq:local reg. of q}
h_{\gamma}(z):=
\begin{cases}
\frac{z}{|z|} &\text{ if }~\gamma |z|-1 \geq \frac{1}{2\gamma},\\
 \frac{z}{|z|} (1- \frac{\gamma}{2} (1- \gamma |z|+\frac{1}{2\gamma})^2) &\text{ if }~\gamma |z|-1 \in (-\frac{1}{2\gamma}, \frac{1}{2\gamma}),\\
\gamma z &\text{ if }~\gamma |z|-1 \leq -\frac{1}{2\gamma}.
\end{cases}
\end{equation}
Note that $| \cdot |_{\gamma}$ has one degree higher regularity than the classical Huber function $\varphi$ in \eqref{HuberL1L2b}.
This higher-order Huber-type smoothing has been previously used in  \cite{bilevellearning,DeLosReyes2017} for similar bilevel problems since it endows the problem \eqref{bilevel_reg} with further regularity, compare Theorem \ref{theo:lin_state}.
% \begin{equation}\label{g_gamma}
% g_{\gamma}(z)=
% \begin{cases}
% |z| - \frac{1}{2\gamma} -\frac{1}{24\gamma^{3}} & \text{ if}\quad |z|\ge \frac{1}{\gamma}+\frac{1}{2\gamma^{2}},\\
% |z|+\frac{1}{6}(1-\gamma |z|+\frac{1}{2\gamma})^{3}- \frac{1}{2\gamma} -\frac{1}{24\gamma^{2}}
% & \text{ if}\quad \frac{1}{\gamma}-\frac{1}{2\gamma^{2}}<|z|< \frac{1}{\gamma}+\frac{1}{2\gamma^{2}},\\
% \frac{\gamma}{2}|z|^{2} & \text{ if}\quad |z|\le \frac{1}{\gamma}-\frac{1}{2\gamma^{2}}.
% \end{cases}
% \end{equation}
% Note that its $C^1$ derivative $h_{\gamma}: \RR^{d} \to \RR^{d}$ reads:
% \begin{equation}\label{h_gamma}
% h_{\gamma}(z)=
% \begin{cases}
% \frac{z}{|z|} & \text{ if}\quad |z|\ge \frac{1}{\gamma}+\frac{1}{2\gamma^{2}},\\
% \frac{z}{|z|} & \text{ if}\quad \frac{1}{\gamma}-\frac{1}{2\gamma^{2}}<|z|< \frac{1}{\gamma}+\frac{1}{2\gamma^{2}},\\
% \gamma z  & \text{ if}\quad |z|\le \frac{1}{\gamma}-\frac{1}{2\gamma^{2}}.
% \end{cases}
% \end{equation}

We then similarly regularise the IC fidelity term as:
\begin{equation}\label{Phi_gamma}
\Phi_{\epsilon,\gamma}^{\lambda_{1},\lambda_{2}}(u,f):= \min_{v\in L^{2}(\om)}~\left\{ \mathcal{G}_{\epsilon,\gamma}^{\lambda_1,\lambda_2}(v,u,f):=\frac{\epsilon}{2}\| v \|_{L^2(\om)}^2+ \lambda_1\|v\|_{\gamma, L^{1}(\om)}+\frac{\lambda_{2}}{2} \|f-u-v\|_{L^{2}(\om)}^{2}\right\},
\end{equation}
where $\|v\|_{\gamma, L^{1}(\om)}$ is defined analogously as above.
For simplicity, from now on, we will assume  that $f\in L^{2}(\om)$. 
%which, combined with the  $H^{1}$ regularisation in \eqref{bilevel_reg}, ensures that the optimum in \eqref{Phi_gamma} lies in $ L^{2}(\om)$.

%We recall that the $H^1$ regularisation on $u$ is needed to prove the differentiability of the solution map $\mathcal{S}:\bm{\lambda}\mapsto u_{\bm{\lambda}}$, see \cite{bilevellearning, noiselearning}.

%
%Let $Du=\nabla u~dx+ D_s u$ be the Lebesgue decomposition of the two-dimensional distributional gradient $Du$ into its absolutely continuous $\nabla u~dx$ and  singular $D_s u$ parts, where $dx$ is the Lebesgue measure in $\RR^2$. We consider for $\gamma\gg 1$ the Huber regularisation of  the TV regularisation \textcolor{red}{metti referenza a def Huber } $|Du|$ :
%\begin{equation}   \label{huberregularTV}
%| Du |_\gamma (\om) := \int_\om |\nabla u|_\gamma~dx + \int_\om |D_s u|.
%\end{equation}

%Summarising, given a noisy image $f\in L^2(\Omega)$, we then seek parameters $\lambda_1$ and $\lambda_2$ in \eqref{L1L2_fid_def} solving the following problem defined for a suitable cost functional $F:H^1(\Omega)\to\RR^+$:
%\begin{equation}   \label{bilevelTVIC1}
%\min_{\lambda_1, ~\lambda_2\geq 0}~ F(u_{\lambda_1,\lambda_2})
%\end{equation}
%subject to:
%\begin{equation} \label{bilevelTVIC2}
%u_{\lambda_1,\lambda_2} \in {\operatorname{argmin}}_{u\in H^1(\om)} \left\{ J^{\gamma,\varepsilon}(u,\lambda_1,\lambda_2) := \frac{\varepsilon}{2} \| u \|_{H^1(\Omega)}^2+|Du|_\gamma(\om) + \Phi^{\lambda_1,\lambda_2}(u,f) \right\}, 
%\end{equation}
%where $\Phi^{\lambda_1,\lambda_2}$ is defined as in \eqref{L1L2_fid_def}.
%

Inspired by \cite{bilevellearning, DeLosReyes2017}, we focus on two main  choices of $F$.
% related to the standard PSNR and SSIM quality measures. 
Namely, we consider
 the $L^2$ cost corresponding to Peak Signal to Noise Ratio (PSNR) optimisation
\begin{equation}   \label{costfunctL2}
F_{L^{2}}(u_{\lambda_1,\lambda_2}):=\| u_{\lambda_1,\lambda_2} - \tilde{u} \|_{L^2(\om)}^2,\quad \text{for training data $\tilde{u}\in L^2(\Omega)$,}
\end{equation}
and the Huberised TV cost, which is related to quality measures that are more adjusted  to actual human perception, such as the Structural Similarity Index (SSIM):
\begin{equation}   \label{costfunctHuber}
 \qquad F_{L_{\gamma}^{1}D}(u_{\lambda_1,\lambda_2}):=\| D(u_{\lambda_1,\lambda_2} - \tilde{u}) \|_{\cM,\gamma},\quad \text{for training data $\tilde{u}\in \bv(\Omega)$.}
\end{equation}
%\end{itemize}
Here, $\|Du\|_{\cM,\gamma}:= \int_{\om} |\nabla u|_\gamma~dx +|D^{s}u|(\om)$ so that the smooth Huber-type regularisation is applied on the absolutely continuous part of $|Du|$. 
For the abstract formulation of \eqref{bilevel_reg} in terms of a general $F$, we refer the reader to \cite{DeLosReyes2017}.

\newtheorem{proper_cost}[HuberL1L2]{Remark}

\begin{proper_cost}
Note that if $\lambda_{1}=0$ and/or $\lambda_{2}=0$, then $\Phi^{\lambda_{1},\lambda_{2}}(u,f)=0$ for every $u\in L^{1}(\om)$. In that case every constant function $c$ is a minimiser of the lower level problem of \eqref{bilevel_unreg}. Thus we trivially have $\inf_{\lambda_{1},\lambda_{2}\ge 0}~ F(u_{\lambda_{1},\lambda_{2}})< \infty$ where $F$ is taken to be either $F_{L^{2}}$ or $F_{L_{\gamma}^{1}D}$.
\end{proper_cost}

\subsection{Well-posedness of the $\tv$--$\mathrm{IC}$ bilevel problem}

We now discuss the well-posedness of the bilevel problems \eqref{bilevel_unreg} and \eqref{bilevel_reg}.
For this type of problems, it is a common practice to impose an extra box constraint on the parameters $\lambda_{1}, \lambda_{2}$ in order to ensure existence of solutions (see \cite{bilevellearning}), although generalisations to unbounded intervals are also possible \cite{DeLosReyes2017}. The following proposition says that for problem \eqref{bilevel_unreg} (i.e., with no upper bound constraints on the parameter domain) existence of solutions may fail.
%Note that the proofs of the next two propositions make use of the results of Section \ref{sec:analysis}.

\newtheorem{non-wellposened}[HuberL1L2]{Proposition}
\begin{non-wellposened}\label{lbl:non-wellposened}
%Suppose that the cost function $F$ has the property $F(u)=0$ if and only if $u=\tilde{u}$.
There exist data $f\in L^{2}(\om)$ and training data $\tilde{u}\in L^{2}(\om)$ such that the non-smooth bilevel problem \eqref{bilevel_unreg} does not have a solution for  the cost function $F_{L^2}$.
\end{non-wellposened}

%\textcolor{red}{The $F_{L^2}$ is defined for training data in $L^2$, here $\tilde{u}$ is $L^1$?}

\begin{proof}
Take $f$ to be any non-constant function in $\bv(\om)\cap L^{2}(\om)$ with the property that there exists $\lambda^{\ast}$ such that 
\begin{equation*}
f=\underset{u\in\bv(\om)}{\operatorname{argmin}}~|Du|(\om)+\lambda_{1}\|f-u\|_{L^{1}(\om)} 
\end{equation*}
for all $\lambda_{1}\ge \lambda^{\ast}$. We note that there are a plethora of such functions, in particular this holds for any one-dimensional function in $\bv(\om)\cap L^{\infty}(\om)$, see also \cite{chanL1}. Now set $\tilde{u}=f$. Then from Proposition \ref{lbl:not_exact}, we have that
 $F_{L^2}(u_{\lambda_{1},\lambda_{2}})>0$ for every $\lambda_{1},\lambda_{2}\ge 0$. Now fix $\lambda_{1}\ge \lambda^{\ast}$ and let $\lambda_{2}^{(n)}\to\infty$ as $n\to\infty$. Then according to Corollary \ref{lbl:u_l2_to_inf} we have that $u_{\lambda_{1},\lambda_{2}^{(n)}}\to f$ weakly$^{\ast}$ in $\bv(\om)$ as $n\to\infty$. The function $f$ can in fact be chosen such that this convergence is even uniform, see for instance the example in Proposition \ref{lbl:prop_exact}. Then it is clear that $F_{L^2}(u_{\lambda_{1},\lambda_{2}^{(n)}})\to 0$, which means that
\[\inf_{\lambda_{1},\lambda_{2}\ge 0} F_{L^2}(u_{\lambda_{1},\lambda_{2}})=0.\]
Since  $F(u_{\lambda_{1},\lambda_{2}})>0$ for every $\lambda_{1},\lambda_{2}\ge 0$, we have that for such choices of $f$ and $\tilde{u}$, the bilevel problem \eqref{bilevel_unreg} does not have solutions.
\end{proof}
 
Similarly, we can use the example of Proposition \ref{lbl:prop_exact} to  show that the bilevel problem \eqref{bilevel_unreg} with  $F_{L^1_\gamma D}$ cost may not have a solution either, in the case when $\lambda_{1}$ and $\lambda_{2}$ are unbounded.

On the other hand, in the following Proposition we show that the existence of solutions of \eqref{bilevel_unreg} is always guaranteed whenever box constraints on $\lambda_{1},\lambda_{2}$ are considered.

%Note that its proof also uses results from the first part of the paper.
%The proof follows standard methods.% see \cite{}.

% that the problem \eqref{bilevelTVIC1}-\eqref{bilevelTVIC2} has a solution. 
%Similarly as in \cite{DeLosReyes2017}, we restrict to the simplified case when box constraints on the parameters $\lambda_1$ and $\lambda_2$ are assumed. The extension to the unbounded case can be done using the general existence result proved in \cite{interiorpaper2015}. For the proof of this first result we make use of the implicit formulation of the IC fidelity term $\Phi$ where the dependence on the auxiliary minimisation variable $v$ is left hidden. This simplifies the analysis due to the standard properties of $\Phi$ which can be found, e.g., in \cite{bauschkecombettes}.

\newtheorem{existence_bilevel}[HuberL1L2]{Proposition}
\begin{existence_bilevel}[Well-posedness of  \eqref{bilevel_unreg}  with box constraints]\label{lbl:existence_bilevel}
The bilevel problem \eqref{bilevel_unreg} with the extra box constraints $0\leq \lambda_i\leq L_i$, $L_i>0$ for $i=1,2$ admits an optimal solution $(\hat{\lambda}_1,\hat{\lambda}_2)$ for both choices of cost functionals \eqref{costfunctL2}--\eqref{costfunctHuber}.
\end{existence_bilevel}

\begin{proof}
Let $(\lambda_{1}^{(n)},\lambda_{2}^{(n)})\in C:=\left\{ (\lambda_1,\lambda_2): 0\leq \lambda_1\leq L_1, 0\leq \lambda_2\leq L_2\right\}$ be a  minimising sequence for \eqref{bilevel_unreg}.  Let us denote by $u_n:=u_{\lambda_{1}^{(n)},\lambda_{2}^{(n)}}$ the corresponding solution to the lower level problem corresponding to the parameter pair $(\lambda_{1}^{(n)}, \lambda_{2}^{(n)})$.
%, as well as by $J(u,\lambda_{1}^{(n)},\lambda_{2}^{(n)})$ the corresponding lower level energy. Then we have
%$$
%J (u_n,\lambda_{1}^{(n)},\lambda_{2}^{(n)}) \leq J(u,\lambda_{1}^{(n)},\lambda_{2}^{(n)}),\quad \text{for all } u\in \bv(\Omega),
%$$
%and by choosing $u= 0$ we get that for every $n\in\NN$
%\begin{equation}   \label{unif:bound:bilevel}
%J(u_n,\lambda_{1}^{(n)},\lambda_{2}^{(n)}) \leq \Phi^{\lambda_{1}^{(n)},\lambda_{2}^{(n)}}(0,f) \leq \frac{L_2}{2}\|f\|_{L^2(\om)}^2 = : K,
%\end{equation}

First suppose that after some index $n_{0}$ at least one of the terms $(\lambda_{1}^{(n)})_{n\in\NN}$ and $(\lambda_{2}^{(n)})_{n\in\NN}$ is zero. This means that  $u_{n}=c_{n}$ are constants for $n\ge n_{0}$. Now if $F=F_{L^{2}}$, due to the coercivity of this functional we have that $(c_{n})_{n\in\NN}$ is bounded, so $(u_{n})_{n\in\NN}$ is bounded in $\bv(\om)$. If $F=F_{L_{\gamma}^{1}}$ then it is obvious that $F(u_{n})=\|D\tilde{u}\|_{\mathcal{M},\gamma}$ for $n\ge 0$. Thus in this case  every constant function trivially solves the bilevel problem \eqref{bilevel_unreg}. 

We can then assume that $\lambda_{1}^{(n)},\lambda_{2}^{(n)}>0$ for every $n\in\NN$.
We claim again that the sequence $u_{n}$ is bounded in $\bv(\om)$. Indeed, we have  that for every $n
\in\NN$ we can bound the $\mathrm{TV}$ term as
\begin{equation*}
|Du_{n}|(\om)\le |Du_{n}|(\om)+ \Phi^{\lambda_{1}^{(n)},\lambda_{2}^{(n)}}(u_{n},f)\le  \Phi^{\lambda_{1}^{(n)},\lambda_{2}^{(n)}}(0,f) \leq \frac{L_2}{2}\|f\|_{L^2(\om)}^2. 
\end{equation*}
To bound $u_{n}$ in $L^{1}(\om)$, we separate the two cases depending on whether the sequence $\frac{\lambda_{1}^{(n)}}{\lambda_{2}^{(n)}}$ is bounded or not.

\noindent If this sequence is bounded by some $K>0$, we observe that $u_{n}$ is also a miminiser of 
\[\min_{u\in\bv(\om)}\Phi^{1,\tfrac{\lambda_{2}^{(n)}}{\lambda_{1}^{(n)}}}(u,f)+\frac{1}{\lambda_{1}^{(n)}}|Du|(\om),\]
which implies that
\[\Phi^{1,\tfrac{\lambda_{2}^{(n)}}{\lambda_{1}^{(n)}}}(u_{n},f)\le \Phi^{1,\tfrac{\lambda_{2}^{(n)}}{\lambda_{1}^{(n)}}}(0,f).\]
Then we have the following successive bounds
\begin{align*}
\int_{\om}|f-u_{n}|\,dx
&= \int_{|f-u_{n}|<\frac{\lambda_{1}^{(n)}}{\lambda_{2}^{(n)}}} |f-u_{n}|\, dx 
+   \int_{|f-u_{n}|\ge\frac{\lambda_{1}^{(n)}}{\lambda_{2}^{(n)}}} \frac{\lambda_{1}^{(n)}}{2\lambda_{2}^{(n)}}\,dx \\
& + \int_{|f-u_{n}|\ge\frac{\lambda_{1}^{(n)}}{\lambda_{2}^{(n)}}} |f-u_{n}| - \frac{\lambda_{1}^{(n)}}{2\lambda_{2}^{(n)}}\, dx 
\le K |\om| +\frac{1}{2} \int_{\om} |f-u_{n}|\,dx + \Phi^{1,\tfrac{\lambda_{2}^{(n)}}{\lambda_{1}^{(n)}}}(u_{n},f),
\end{align*}
whence we get:
\begin{align*}
& \frac{1}{2} \int_{\om}|f-u_{n}|\,dx \le K |\om| +   \Phi^{1,\tfrac{\lambda_{2}^{(n)}}{\lambda_{1}^{(n)}}}(0,f)\\
&\le K |\om| + \int_{|f|\ge \frac{\lambda_{1}^{(n)}}{\lambda_{2}^{(n)}}} |f|\,dx
+ \frac{\lambda_{2}^{(n)}}{2 \lambda_{1}^{(n)}}\int_{|f|< \frac{\lambda_{1}^{(n)}}{\lambda_{2}^{(n)}}} |f|^{2}\,dx 
 \le K |\om| + \int_{\om} |f|\,dx +\frac{K}{2}  |\om|. 
 \end{align*}

\noindent Suppose now that the sequence $\frac{\lambda_{1}^{(n)}}{\lambda_{2}^{(n)}} $ is unbounded. This means that there exists a (non-relabelled) subsequence $\frac{\lambda_{2}^{(n)}}{\lambda_{1}^{(n)}} \to 0$. Since $\lambda_{1}^{(n)}$ is bounded, this further implies that $\lambda_{2}^{(n)}\to 0$. Then from Proposition \ref{lbl:l2tozero} we get that $u_{n}\to f_{\om}$ weakly$^{\ast}$ in $\bv(\om)$ and that, in particular, $u_{n}$ is bounded in $L^{1}(\om)$. So in both cases we have that $u_n$ is bounded in $\bv(\om)$.

Having this combined with the boundedness of the sequence $(\lambda_{1}^{(n)},\lambda_{2}^{(n)})_{n\in\NN}$ implies that we can then extract a further non-relabelled subsequence $ (\lambda_1^{(n)}, \lambda_{2}^{(n)}, u_n)_{n\in\NN}$ converging weakly$^\ast$ in $\RR \times \RR \times \bv(\Omega)$ to a limit point $(\hat{\lambda}_1, \hat{\lambda}_2, \hat{u} )$. In particular, this entails that $u_n\to \hat{u}$ strongly in $L^1(\Omega)$. Similarly as in the proof of Proposition \ref{lbl:Gamma_convergence} one can now show that the sequence $J(\cdot,\lambda_{1}^{(n)},\lambda_{2}^{(n)})$ $\Gamma$-converges to  $J(\cdot,\hat{\lambda}_1, \hat{\lambda}_2)$, with respect to the strong topology in $L^{1}(\om)$, which means that $\hat{u}$ is a minimiser of the lower level problem with parameters $(\hat{\lambda}_{1},\hat{\lambda}_{2})$.
%and therefore, again up to a subsequence, $u_n$ converges to $\hat{u}$ pointwise almost everywhere.
% By standard continuity property of the data term $\Phi^{\lambda_1,\lambda_2}$ (see, e.g., \cite{bauschkecombettes}) and weakly lower semicontinuity of the Huberised-TV term  (see \cite{Dem}) we then get:

Now, using \cite[Corollary 7.20]{dalmasogamma} we finally have:
\begin{align*}
  \Phi^{\hat{\lambda}_1,\hat{\lambda}_2}(\hat{u},f) + |D\hat{u}|(\Omega) 
& = \lim_{n\to\infty}   \Phi^{\lambda_{1}^{(n)},\lambda_{2}^{(n)}}(u_{n},f) +|Du_{n}|(\Omega), 
\end{align*}
which completes the proof.
%where we have used the notation $v_n : = \argmin_{v\in L^2(\Omega)}~\lambda_{1}^{(n)}\|v\|_{L^1(\Omega)}+\frac{\lambda_{2}^{(n)}}{2}\|f-u_n-v\|_{L^2}^2$ and where the second inequality holds by the convergence of the sequence $(\lambda_{1}^{(n)},\lambda_{2}^{(n)})$ and by the uniform bounds \eqref{unif:bound:bilevel} on the $L^1$ and $L^2$ terms. 
%
%To conclude, observe that in both cases \eqref{costfunctL2} and \eqref{costfunctHuber}, the cost functional $F$ is weakly lower semicontinuous and the set $C$ is closed, hence weakly closed by Mazur's Lemma.

\end{proof}

\newtheorem{rmk:regTOnonreg}[HuberL1L2]{Remark}

\begin{rmk:regTOnonreg}
Similarly, one can prove the existence of solutions to the regularised bilevel problem \eqref{bilevel_reg}. Note that in \cite{interiorpaper2015} similar techniques are used to prove analogous results for general regularisers and data fidelities. There, the authors proved also the outer-semicontinuity property of the solution map $\mathcal{S}$ which guarantees that the minimisers of the regularised problem converge towards the minimisers of the one where $\epsilon=0$.
%Furthermore, under suitable assumptions on the initial data, strict positivity of the optimal parameters is therein proved even in the case of box constraints.
\end{rmk:regTOnonreg}

\newtheorem{rmk:compactness}[HuberL1L2]{Remark}

The quadratic $H^1$ and $L^2$ regularisations in \eqref{bilevel_reg} and \eqref{Phi_gamma} are required to ensure the differentiability of the solution map, as we are going to highlight in the following. Note that in this section we make use of the formulation \eqref{bilevel_reg} where the two variables $u_{\lambda_{1},\lambda_{2}}, v_{\lambda_1,\lambda_2}$ are treated jointly so that the lower-level problem actually reads:
\begin{align}  \label{eq:lower_level_joint}
( u_{\lambda_{1},\lambda_{2}}, v_{\lambda_1,\lambda_2} ) & \in \argmin_{\substack{u\in H^{1}(\om)\\v\in L^2(\om)}}~ \frac{\epsilon}{2}\Big( \|u\|_{H^{1}(\om)}^{2}+\|v\|_{L^2(\om)}^2\Big)+  \|\nabla u\|_{\gamma, L^{1}(\om)} \\ & +  \lambda_1\|v\|_{\gamma, L^{1}(\om)}+\frac{\lambda_{2}}{2} \|f-u-v\|_{L^{2}(\om)}^{2}. \notag
\end{align}
An alternative analysis could exploit the characterisation of the IC data fidelity term given by Proposition \ref{lbl:HuberL1L2} and consider the minimisation over $u$ only. In such case, we believe that only an $H^1$-regularisation on $u$ would be enough for the following proofs. 
Here, however, we stick with the joint approach to compare our results with the ones derived formally in \cite{calatroni_mixed}. 

Note that whenever $u_{\lambda_1,\lambda_2}$ is given, one can compute the corresponding $v_{\lambda_1,\lambda_2}$ by simply solving the optimisation problem \eqref{Phi_gamma}. In particular, the following proposition makes explicit a property of $v_{\lambda_1,\lambda_2}$ which will be needed in the following. We refer the reader to \cite[Remark 2.1]{calatroni_mixed} for a similar characterisation in the case non-regularised case.
%We omit its proof since it follows directly by the definition of proximal operators as observed in \cite[Remark 2.1]{calatroni_mixed}.

\newtheorem{prop:proximal_map}[HuberL1L2]{Proposition}

\begin{prop:proximal_map}  \label{prop:prox_map}
Let $u, f$ in $L^2(\om)$, $0\leq \lambda_i\leq L_i$  with $L_i>0$, $i=1,2$ and $\epsilon>0$. Let $v_{\lambda_1,\lambda_2}$ the minimiser of the functional $\mathcal{G}_{\epsilon,\gamma}^{\lambda_1,\lambda_2}(\cdot,u,f)$ defined in \eqref{Phi_gamma}.
%let $v_{\lambda_1,\lambda_2}$ be defined as:
%$$
%v_{\lambda_1,\lambda_2} = \argmin_{v\in L^2(\om)}~ \lambda_1\|v\|_{\gamma, L^{1}(\om)}+\lambda_{2} \|f-u-v\|_{L^{2}(\om)}^{2}.
%$$
%Then, 
There holds:
\begin{equation}  \label{eq:firm_nonexpansiveness}
v_{\lambda_1,\lambda_2}= v_{\lambda_1,\lambda_2}(u)= \argmin_{v\in L^2(\om)}~\mathcal{G}_{\epsilon,\gamma}^{\lambda_1,\lambda_2}(v,u,f)= \mathrm{prox}_{\frac{\lambda_1}{\epsilon+\lambda_2}\|\cdot\|_{\gamma,L^1(\Omega)}}\left(\frac{\lambda_2(f-u)}{\epsilon+\lambda_2}\right),
\end{equation}
where for any $z\in L^2(\Omega)$, $\mathrm{prox}_{\tau g}(z)$ denotes the proximal-mapping operator in $L^2(\Omega)$ of the function $g$ with parameter $\tau$. In particular, $v_{\lambda_1,\lambda_2}: z\mapsto \mathrm{prox}_{\frac{\lambda_1}{\epsilon+\lambda_2}\|\cdot\|_{\gamma,L^1(\Omega)}} (z)$ is a firmly non-expansive operator and it is therefore $1$-Lipschitz continuous, i.e.:
$$
\|v_{\lambda_1,\lambda_2}(u_1)-v_{\lambda_1,\lambda_2}(u_2)\|_{L^2(\om)}\leq \|u_1 - u_2\|_{L^2(\om)},\quad\text{for all }u_1,u_2\in L^2(\om).
$$
\end{prop:proximal_map}

\begin{proof}
Straightforward calculations give:
\begin{align} 
v_{\lambda_1,\lambda_2}& =\argmin_{v\in L^2(\Omega)}~ \frac{\epsilon}{2}\| v \|_{L^2(\om)}^2+ \lambda_1\|v\|_{\gamma, L^{1}(\om)}+\frac{\lambda_{2}}{2} \|f-u-v\|_{L^{2}(\om)}^{2} \notag \\
& = \argmin_{v\in L^2(\Omega)}~ \|v\|_{\gamma, L^{1}(\om)} + \frac{\epsilon+\lambda_2}{2\lambda_1}\| v \|_{L^2(\om)}^2 + \frac{1}{2\frac{\lambda_1}{\lambda_2}\left(1+\frac{\epsilon}{\lambda_2}\right)}\|f-u\|^2_{L^2(\Omega)} - \frac{\lambda_2}{\lambda_1} \int_{\om}(f-u)v\,dx \notag \\
& = \argmin_{v\in L^2(\Omega)}~\|v\|_{\gamma, L^{1}(\om)} + \frac{1}{2\frac{\lambda_1}{\epsilon+\lambda_2}}\left\Vert\frac{\lambda_2(f-u)}{\epsilon+\lambda_2} - v \right\Vert_{L^2(\Omega)}^2 \notag\\
& = \mathrm{prox}_{\frac{\lambda_1}{\epsilon+\lambda_2}\|\cdot\|_{\gamma,L^1(\Omega)}}\left(\frac{\lambda_2(f-u)}{\epsilon+\lambda_2}\right).  \notag
\end{align}
The firm non-expansiveness property \eqref{eq:firm_nonexpansiveness} follows then directly from \cite[Proposition 12.27]{bauschkecombettes}. 
\end{proof}

\subsection{Optimality system}

We now study in more detail the bilevel problem \eqref{bilevel_reg} with \eqref{eq:lower_level_joint} as lower level problem and prove the existence of Lagrange multipliers by deriving the optimality system characterising its stationary points. As a by product, we find an easy formula to compute the gradient of the cost functional in terms of its adjoint state, which simplifies the design of the gradient-based algorithm which we are going to use to solve \eqref{bilevel_reg} in an efficient way.

We start defining the Hilbert space $\mathcal{H}:=H^1(\Omega)\times L^2(\Omega)$ endowed with the scalar product $(z,w)_{\mathcal{H}}:=(z_1,w_1)_{H^1(\om)} + (z_2,w_2)_{L^2(\om)}$ for all $z,w\in \mathcal{H}$.
We further denote by $\mathcal{S}:\RR^2\to \mathcal{H}$ the solution map $\mathcal{S}:(\lambda_1,\lambda_2) \mapsto (u_{\lambda_1,\lambda_2}, v_{\lambda_1,\lambda_2})$ which assigns to the optimal parameters $(\lambda_1,\lambda_2)$ the corresponding solution pair of \eqref{eq:lower_level_joint}. 
To avoid heavy notations, we will omit the explicit dependence of the pair $(u,v)$ on the parameters $(\lambda_1, \lambda_2)$ unless explicitly needed. Note that we take the whole space $\RR^2$ as differentiability set although  from an imaging point of view the use of negative parameters clearly does not make sense. However, the following differentiability result holds in such case as well. 

\medskip

Recalling the definition of $h_\gamma$ in \eqref{eq:local reg. of q}, we have that in correspondence with an optimal pair $y:=(u,v)\in \mathcal{H}$, we have that the following variational equality holds true for all test functions $\Psi:=(\psi_1,\psi_2) \in \mathcal{H}$:
\begin{multline}   \label{opt:var:eq}
\epsilon\left(y,\Psi \right)_{\mathcal{H}} + \int_\Omega h_\gamma(\nabla u)\nabla\psi_1\,dx + \lambda_1\int_\Omega h_\gamma(v)\psi_2\,dx  + \lambda_2\int_\Omega (u+v-f)(\psi_1 + \psi_2)\,dx = 0.
\end{multline}
We now prove the main differentiability result.

%where thanks to Lemma \ref{lemma:diff:IC} the derivative of $\Phi^{\lambda_1,\lambda_2}$ is well-defined.

\newtheorem{differentiability_solution_map}[HuberL1L2]{Theorem}

\begin{differentiability_solution_map}[Fr\'echet differentiability of the solution map]   \label{theo:lin_state}
The solution operator $\mathcal{S}: \RR^2\to \mathcal{H}$ which assigns to each parameter pair $(\lambda_1,\lambda_2)$ the element $y:=(u,v)=\mathcal{S}(\lambda_1,\lambda_2)$, solution of the $\tv$--$\mathrm{IC}$ denoising problem  \eqref{bilevel_reg}--\eqref{Phi_gamma} is Fr\'echet differentiable. In particular, for any $\theta=(\theta_1,\theta_2) \in \RR^2$ its Fr\'echet derivative is the unique solution $z:= \mathcal{S}'(\lambda_1,\lambda_2)[\theta_1,\theta_2] \in \mathcal{H}$  of the following linearised equation:
\begin{multline}  \label{linearised:state}
\epsilon\left(z,\Psi \right)_{\mathcal{H}} +  \int_\Omega h'_\gamma(\nabla u)  \nabla z_1 \nabla\psi_1\,dx + \lambda_1 \int_\Omega h'_\gamma(v) z_2 \psi_2\,dx + \theta_1\int_\Omega h_\gamma(v)\psi_2\, dx \\+ \lambda_2 \int_\Omega  (z_1+z_2) (\psi_1+\psi_2)\,dx  + \theta_2 \int_\Omega (u+v-f) (\psi_1 + \psi_2)\,dx = 0,
\end{multline}
for all $\Psi=(\psi_1,\psi_2)\in \mathcal{H}$.
\end{differentiability_solution_map}

\begin{proof}
Thanks to the ellipticity of the scalar product in $\mathcal{H}$ and the monotonicity of $h_\gamma$, existence and uniqueness of $z$ of \eqref{linearised:state} are guaranteed by Lax-Milgram theorem.

Now, we want to show that $z$ is the Fr\'echet derivative of $\mathcal{S}$. To do that, given $\theta=(\theta_1,\theta_2)\in\RR^2$ let us define $y^{+}=(u^{+},v^{+}):=\mathcal{S}(\lambda_1+\theta_1,\lambda_2+\theta_2)$. We aim to show that $\xi=(\xi_1,\xi_2):=y^+-y-z\in\mathcal{H}$ satisfies $\|\xi\|_{\mathcal{H}}=o(|\theta|)$. 

Writing \eqref{opt:var:eq} for $y$ and $y^+$, and \eqref{linearised:state} for $z$ and combining them together, we get that for every $\Psi\in \mathcal{H}$ there holds:
\begin{multline}
\epsilon\left(\xi,\Psi\right)_{\mathcal{H}} + \int_\Omega \left( h_\gamma(\nabla u^+) - h_\gamma( \nabla u) \right) \nabla\psi_1\,dx - \int_\Omega h'_\gamma(\nabla u) \nabla z_1 \nabla\psi_1\,dx \notag\\
%+ \left( \lambda_1 + \theta_1 \right) \int_\Omega h_\gamma(v^+) \psi_2~dx - \lambda_1\int_\Omega h_\gamma(v) \psi_2~dx - \lambda_1 \int_\Omega h'_\gamma(v) z_2 \psi_2~dx - \theta_1\int_\Omega h_\gamma(v)\psi_2~ dx \\ 
+ \lambda_1 \int_\Omega \left( h_\gamma(v^+)-h_\gamma(v) \right) \psi_2\,dx + \theta_1 \int_\Omega \left( h_\gamma(v^+)-h_\gamma(v) \right) \psi_2\,dx  - \lambda_1 \int_\Omega h'_\gamma(v) z_2 \psi_2\,dx \\
+ \lambda_2 \int_\Omega (\xi_1 +\xi_2)( \psi_1 +\psi_2) \,dx +  \theta_2 \int_\Omega \left( (u^+ - u) + (v^+-v) \right)( \psi_1 + \psi_2 )\,dx=0. \notag
\end{multline}
We now add and subtract the terms:
\begin{align*}
& \int_\Omega h'_\gamma(\nabla u)\left(\nabla(u^+-u)\right)\nabla \psi_1\,dx,\quad\text{and}\quad
\lambda_1\int_\Omega h'_\gamma(v)\left(v^+-v\right)\psi_2\,dx, 
\end{align*}
thus getting:
\begin{multline*}
\epsilon\left(\xi,\Psi\right)_{\mathcal{H}} + \int_\Omega h'_\gamma(\nabla u)\nabla\xi_1 \nabla\psi_1\,dx + \lambda_1\int_\Omega h'_\gamma(v)\xi_2 \psi_2\,dx 
+ \lambda_2 \int_\Omega (\xi_1+\xi_2)( \psi_1 +\psi_2)~dx  
 \\
= -\int_\Omega \left( h_\gamma(\nabla u^+) - h_\gamma(\nabla u) - h'_\gamma(\nabla u)(\nabla(u^+ -u)  )\right)\nabla\psi_1\,dx  
\\ - \lambda_1\int_\Omega \left( h_\gamma(v^+)-h_\gamma(v)-h'_\gamma(v)(v^+-v)\right)\nabla\psi_2\,dx -\theta_1\int_\Omega \left( h_\gamma(v^+)-h_\gamma(v)\right)\psi_2\,dx \\
- \theta_2\int_\Omega \left( (u^+-u)+(v^+-v) \right) (\psi_1+\psi_2)\,dx,\quad\text{for all }\Psi\in\mathcal{H}.
\end{multline*}
We now choose $\Psi=\xi$. Under this choice, by monotonicity of $h'_\gamma$ we have that the last three terms on the left-hand side become non-negative. We further deduce:
\begin{align*}
\|\xi\|_{\mathcal{H}} & \leq C\Big( \|h_\gamma(\nabla u^+) - h_\gamma(\nabla u) - h'_\gamma(\nabla u)(\nabla (u^+ -u))\|_{L^2(\Omega)} \\ 
& + |\lambda_1|~\| h_\gamma(v^+)-h_\gamma(v)-h'_\gamma(v)(v^+-v)\|_{L^2(\Omega)} + |\theta_1|~\| h_\gamma(v^+)-h_\gamma(v) \|_{L^2(\Omega)} \\
& + |\theta_2|~\|u^+-u\|_{L^2(\Omega)} + |\theta_2|~\|v^+-v\|_{L^2(\Omega)}\Big).
\end{align*}
where $C$ is a generic positive and finite constant which may change from line to line.
By the differentiability and Lipschitz continuity of $h_\gamma$ and $h'_{\gamma}$, we have:
\begin{align} \label{ineq:proof}
\|\xi\|_{\mathcal{H}} & \leq  C\Big( o\left( \|\nabla u^+ - \nabla u\|_{L^2(\Omega)}\right) + |\lambda_1|~o\left( \|v^+-v\|_{L^2(\Omega)}\right) \\ & + |\theta_2|~\|u^+-u\|_{L^2(\Omega)} 
 + \left(|\theta_1|+|\theta_2|\right)~\|v^+-v\|_{L^2(\Omega)} \Big).  \notag
%\leq C\Big( o\left( \|u^+ - u\|_{W^{1,p}(\Omega)}\right) + |\lambda_1|~o\left( \|v^+-v\|_{L^p(\Omega)}\right) \\
%&+ |\theta|\|u^+-u\|_{L^2(\Omega)} \Big),
\end{align}
%for $p>2$ by standard Sobolev embeddings. \textcolor{red}{How can I conclude?}
We now focus on the terms depending on the difference between $v^+=v_{\lambda_1+\theta_1,\lambda_2+\theta_2}(u^+)$ and $v=v_{\lambda_1,\lambda_2}(u)$.  By triangle inequality we have:
\begin{align*}
\|v^+-v\|_{L^2(\Omega)}& \leq \|v_{\lambda_1+\theta_1,\lambda_2+\theta_2}(u^+) -v_{\lambda_1+\theta_1,\lambda_2+\theta_2}(u) \|_{L^2(\Omega)} + \| v_{\lambda_1+\theta_1,\lambda_2+\theta_2}(u) - v_{\lambda_1,\lambda_2}(u) \|_{L^2(\Omega)}\\
& \leq \| u^+ - u\|_{L^2(\om)} + |\theta|,
\end{align*}
where the last inequality follows from Proposition \ref{prop:prox_map} and from the continuity property of the proximal mapping with respect to its parameter which can be easily checked in our case recalling that $|h_\gamma(\cdot)|\leq 1$.
We then deduce in \eqref{ineq:proof} that:
\begin{equation}  \label{eq:small_o}
\|\xi\|_{\mathcal{H}}\leq C\Big(o(|\theta|) + o(\|u^+ - u\|_{H^1(\Omega)}) \Big)
\end{equation}
Thus, to conclude it only remains to show that $o(\| u^+ - u\|_{H^1(\Omega)}) = o(|\theta|)$. To do that, we use standard Sobolev embeddings and the regularity result of Gr\"{o}ger for second-order systems \cite[Theorem 1]{Groger1989} and get that for $p>2$:
\begin{align*}
\|u^+ - u\|_{H^1(\Omega)}& \leq \|u^+ - u\|_{W^{1,p}(\om)} \leq C|\theta| \Big( \| \mathrm{div}(h_\gamma(\nabla u))\|_{W^{-1,p}(\om)} + \| h_\gamma(\nabla u)\|_{W^{-1,p}(\om)}\Big)\notag\\
& \leq C|\theta| \| h_\gamma(\nabla u)\|_{L^\infty(\om)}\leq C|\theta|,
\end{align*}
since $|h_\gamma(\cdot)|\leq 1$. Combining with \eqref{eq:small_o} this finishes the proof.
\end{proof}

%\textcolor{red}{I have to check last estimate: dependence on $v$? Check Groeger.}

\newtheorem{rmk:groger}[HuberL1L2]{Remark}

\begin{rmk:groger}
The regularity result by Gr\"{o}ger is a classical argument for the proof of Fr\'echet differentiability in similar bilevel problems (see, e.g., \cite{bilevellearning}). Note, however, that the original result in \cite{Groger1989} was proved for $C^2$-regular domains, while its extension to convex Lipschitz domains (such as image domains) has been proved by Dauge in \cite{Dauge1992}.
\end{rmk:groger}

%Inspired by the Lagragian derivation in \cite[Section 7]{calatroni_mixed}, 
We now prove the existence and uniqueness of the adjoint state of the problem \eqref{bilevel_reg}.

\newtheorem{adjoint_states}[HuberL1L2]{Theorem}

\begin{adjoint_states}[Adjoint equation]  \label{theo:adjoint}
Let $(u,v)\in \mathcal{H}$. There exists a unique solution $\Pi:=(p_1,p_2)\in \mathcal{H}$ to the adjoint PDE:
\begin{multline} \label{adj:syst}
\epsilon(\Pi,W)_{\mathcal{H}}+\int_\Omega h'_\gamma(\nabla u)\nabla w_1\nabla p_1\,dx + \lambda_1\int_\Omega h'_\gamma(v)w_2 p_2\,dx  \\
+ \lambda_2\int_\Omega (p_1+p_2)w_1\,dx + \lambda_2\int_\Omega (p_1+p_2)w_2\,dx
= - \int_\Omega F'(u)w_1\,dx,
\end{multline}
for any $W:=(w_1,w_2)\in\mathcal{H}$.
\end{adjoint_states}

\begin{proof}
For $w\in\mathcal{H}$, let us consider the following bilinear form on $\mathcal{H}\times\mathcal{H}$:
\begin{multline}  
a(\Pi,W):= \epsilon(\Pi,W)_{\mathcal{H}}+\int_\Omega h'_\gamma(\nabla u)\nabla w_1\nabla p_1~dx \notag \\ + \lambda_1\int_\Omega h'_\gamma(v)w_2 p_2~dx 
+ \lambda_2\int_\Omega (p_1+p_2)w_1~dx + \lambda_2\int_\Omega (p_1+p_2)w_2~dx. \notag
\end{multline}
The form $a(\cdot,\cdot)$ is trivially symmetric and coercive, since by taking $W=\Pi$, we get:
\begin{multline}
a(\Pi,\Pi) = \epsilon\|\Pi\|^2 _{\mathcal{H}} + \int_\Omega h'_\gamma(\nabla u)\nabla p_1\nabla p_1 
+ \lambda_1\int_\Omega h'_\gamma(v)p_2 p_2~dx + \lambda\int_\Omega (p_1+p_2)^2~dx \geq \epsilon\|\Pi\|^2 _{\mathcal{H}}. \notag
\end{multline}
by monotonicity of $h'_\gamma$. By Lax-Milgram theorem, we infer that there exists a unique solution of  \eqref{adj:syst}.
\end{proof}

Note that by taking $w_2=0$ in \eqref{adj:syst}, we get the optimality condition for $p_1$, i.e.:
\begin{equation*}  
\epsilon(p_1,w_1)_{H^1(\om)}+\int_\Omega h'_\gamma(\nabla u)\nabla w_1\nabla p_1\,dx  
+ \lambda_2\int_\Omega (p_1+p_2)w_1\,dx = - \int_\Omega F'(u)w_1\,dx,
\end{equation*}
for any $w_1\in H^1(\Omega)$. Similarly, for $w_1=0$, we get the optimality condition for $p_2$:
\begin{equation*}
\epsilon\int_\Omega w_2 p_2\,dx + \lambda_1\int_\Omega h'_\gamma(v)w_2 p_2\,dx  + \lambda_2\int_\Omega (p_1+p_2)w_2\,dx
=0,
\end{equation*}
for any $w_2\in L^2(\Omega)$. 

\medskip

Finally, we now combine the results above to derive the optimality system of the bilevel problem \eqref{bilevel_reg}. We recall that by $y=(u,v)\in\mathcal{H}$ we denote the solution pair.

\newtheorem{optimality_system}[HuberL1L2]{Theorem}

\begin{optimality_system}[Optimality system] \label{theo:optimality}
Let $(\bar{\lambda}_1,\bar{\lambda}_2)\in \RR^+\times\RR^+$ an optimal solution of the problem \eqref{bilevel_reg}. Then, there exists a Lagrange multiplier  $\Pi:=(p_1,p_2)\in\mathcal{H}$ and $\mu_1,\mu_2\in\RR^+$ such that the following system holds:
\begin{align}  \label{eq:optimality_system}
&\epsilon\left(y,\Psi\right)_{\mathcal{H}} + \int_\Omega h_\gamma(\nabla u)\nabla \psi_1\,dx + \lambda_1\int_\Omega h_\gamma(v)\psi_2\,dx \notag \\ 
&+ \lambda_2\int_\Omega (u+v-f)(\psi_1 + \psi_2)\, dx = 0, \qquad \text{for all }\Psi=(\psi_1,\psi_2) \in \mathcal{H},  \notag\\
&\epsilon(\Pi,W)_{\mathcal{H}}+\int_\Omega h'_\gamma(\nabla u)\nabla w_1\nabla p_1\,dx + \lambda_1\int_\Omega h'_\gamma(v)w_2 p_2\,dx  \\ 
 &+ \lambda_2\int_\Omega (p_1+p_2)(w_1+w_2)\,dx
= - \int_\Omega F'(u)w_1\,dx,\qquad\text{for all } W=(w_1,w_2) \in \mathcal{H}, \notag  \\ 
&\mu_1:= \int_\Omega h'_\gamma(v)p_2\,dx,\quad \mu_2:=\int_\Omega (f-v-u)(p_1+p_2)\,dx,  \notag \\
&\mu_1\geq 0, \quad \mu_2\geq 0,\quad \mu_1\cdot \bar{\lambda}_1 = \mu_2\cdot \bar{\lambda}_2 = 0. \notag
\end{align}
\end{optimality_system}

\begin{proof}
We can write the bilevel problem \eqref{bilevel_reg} in a reduced form as:
$$
\min_{\lambda_1,\lambda_2\geq 0}~\mathcal{F}(\lambda_1,\lambda_2):= F(u_{\lambda_1,\lambda_2}).
$$
Using \cite[Theorem  3.1]{Zowe1979}, we deduce the existence of multipliers $\mu_1, \mu_2\in\RR^+$ such that:
\begin{align}
&\mu_1 = \nabla_{\lambda_1} \mathcal{F}(\bar{\lambda}_1,\bar{\lambda}_2), \notag\\
&\mu_2 = \nabla_{\lambda_2} \mathcal{F}(\bar{\lambda}_1,\bar{\lambda}_2), \notag\\
&\mu_1\geq 0,\quad \mu_2\geq 0,\quad \mu_1\cdot\bar{\lambda}_1=  \mu_2\cdot\bar{\lambda}_2 =0 . \notag
\end{align}

Computing the gradient of $\mathcal{F}$ by using the chain rule we get:
$$
\nabla \mathcal{F}(\lambda_1,\lambda_2)[\theta_1,\theta_2]= \big( F'(u_{\lambda_1,\lambda_2}), \mathcal{S}'(\lambda_1,\lambda_2)[\theta_1,\theta_2]\big)_{L^2(\om)}=\int_\Omega F'(u_{\lambda_1,\lambda_2})z\,dx
$$
where $z\in\mathcal{H}$ is the linearised state provided by Theorem \ref{theo:lin_state}. Theorem \ref{theo:adjoint} ensures that there exist a Lagrange multiplier $\Pi:=(p_1,p_2)$ satisfying the adjoint equation \eqref{adj:syst}, which entails:
\begin{align}
\int_\Omega F'(u_{\lambda_1,\lambda_2})z\,dx  = & -\epsilon(\Pi,z)_{\mathcal{H}} -\int_\Omega h'_\gamma(\nabla u)\nabla z_1\nabla p_1\,dx  \notag \\
& - \lambda_1\int_\Omega h'_\gamma(v)z_2 p_2\,dx- \lambda_2\int_\Omega (p_1+p_2)(z_1+z_2)\,dx  \notag \\
& = \theta_1\int_\Omega h_\gamma(v)p_2\,dx + \theta_2\int_\Omega (u+v-f)(p_1+p_2)\,dx, \notag
\end{align}
which completes the proof.
\end{proof}

\newtheorem{gradient_formula}[HuberL1L2]{Remark}

\begin{gradient_formula}  \label{rmk:gradient_formula}
Theorem \ref{theo:optimality} provides also the following handy formula for the computation of the gradient of the bilevel problem \eqref{bilevel_reg} with data fidelity \eqref{Phi_gamma} in a reduced form:
\begin{equation}  \label{formula:gradient}
\nabla \mathcal{F}(\lambda_1,\lambda_2)[\theta_1,\theta_2] = \theta_1\int_\Omega h_\gamma(v)p_2~dx + \theta_2\int_\Omega (u+v-f)(p_1+p_2)~dx.
\end{equation}
\end{gradient_formula}

\newtheorem{rmk:comparison}[HuberL1L2]{Remark}

\begin{rmk:comparison}
Compared to the optimality system derived in \cite[Section 7]{calatroni_mixed} via Lagrangian formalism, we note that the optimality system \eqref{eq:optimality_system} presents an additional quadratic $\epsilon$-regularisation on $v$. This is needed in the proof of Theorem \ref{theo:lin_state} to get a uniform, finite  estimate for $\xi_2$ not depending on $\lambda_2$ which may be zero.
\end{rmk:comparison}

\section{Numerical experiments}  \label{sec:numres}

In this section we report some numerical results on the computation of the optimal parameters $(\bar{\lambda}_1,\bar{\lambda}_2)$ of the $\tv$--$\mathrm{IC}$ model by means of the bilevel optimisation strategy described in the previous section. For given training images with a mixture of Gaussian and Salt \& Pepper noise with various intensities, we describe in Algorithm \ref{alg:bilevel} the main steps to compute $(\bar{\lambda}_1,\bar{\lambda}_2)$ by solving the optimality system \eqref{eq:optimality_system} via a second-order BFGS optimisation approach. A general review of second-order numerical methods for PDE-constrained optimisation models can be found in the book \cite{JCbook}, while more details on the numerical realisation of similar bilevel models can be found in \cite[Section 8.3]{bilevellearning} and in \cite[Section 4.1]{DeLosReyes2017}.

In the following numerical computations:
\begin{itemize}
\item[-] We consider test images of size $N\times N\equiv256\times 256$ pixels. The differential operators are discretised using finite difference schemes with mesh step size $h=1/N$. Standard forward/backward differences are used for the discretisation of the divergence/gradient operator, respectively.
\item[-] For illustrative purposes, we report the results obtained using a training set consisting of one training pair $(\tilde{u},\tilde{f})$ only. We recall that $\tilde{f}$ represents the noisy version of $\tilde{u}$ corrupted by a mixture of Gaussian and Salt \& Pepper noise of various intensities which we will specify in each case. The efficient extension to multiple constraints can be done similarly as in \cite{bilevellearning,calatroni2014dynamic}.
\item[-] The lower-level regularised  $\tv$--$\mathrm{IC}$ problem \eqref{bilevel_reg} is solved by means of the SemiSmooth Newton (SSN) algorithm described in \cite{calatroni_mixed} with a warm start. The regularisation parameters are chosen as $\epsilon=10^{-10}$ and $\gamma=10^3$. 
%A warm initialisation of the algorithm is used, i.e. for each iteration the denoised image computed in the previous iteration is used as initial guess for the SSN solver.
The algorithm is stopped if either the difference between two consecutive iterates is below $\texttt{tol}=10^{-6}$ or if the maximum number of iterations $\textsf{maxiter}=35$ is reached.
\item[-] For the outer BFGS iterations, an Armijo line-search with parameter $\eta=10^{-4}$ is employed together with a curvature verification. The Armijo rule:
$$
\mathcal{F}(\bm{\lambda_k} + \alpha_k \bm{d_k}) - \mathcal{F}(\bm{\lambda_k}) \leq \eta \alpha_k \nabla \mathcal{F}(\bm{\lambda_k})^\top \bm{d_k}
$$
is checked at any iteration $k\geq 2$ Here $\bm{\lambda_k}=(\lambda_k^1,\lambda_2^k)$ stands for the parameter pair updated along the iterations, $\bm{d_k}$ for the quasi-Newton descent direction and $\alpha_k$ for the line step. The expression of $\nabla\mathcal{F}$ is given in Remark \ref{rmk:gradient_formula}.
The outer algorithm is stopped when the maximum between the norm of the gradient of the cost functional and the difference of two subsequent iterates is smaller than $\texttt{tol}_1=10^{-6}$.
\item[-] The adjoint equation is solved by means of standard sparse linear solvers.
\item[-] Our validation images 
% have been selected from the Berkeley dataset, freely available online\footnote{\url{https://www2.eecs.berkeley.edu/Research/Projects/CS/vision/bsds/BSDS300/html/dataset/images.html}} and 
are taken from the public domain.
% see Figure \ref{fig:berkeley}.
% \begin{figure}[!h] 
% \begin{center}
% \includegraphics[height=1.5cm]{im20}
% \includegraphics[height=1.5cm]{im33}
% \includegraphics[height=1.5cm]{test47}
% \includegraphics[height=1.5cm]{test52}
% \includegraphics[height=1.5cm]{test61}
% \includegraphics[height=1.5cm]{test71}
% \end{center}
% \caption{Some images from the Berkeley dataset used for validation.}
%  \label{fig:berkeley}
% \end{figure}
\end{itemize}

\begin{algorithm}[h!]
  \caption{Bilevel optimisation algorithm for computing optimal $\lambda_1$ and $\lambda_2$ in \eqref{bilevel_reg}}
 \label{alg:bilevel}
  \begin{algorithmic}
    \Statex {\textbf{Input}: Training pair $(\tilde{u},\tilde{f})$. Regularisation parameters: $\gamma\gg 1$, $\epsilon\ll 1$.}
    \Statex {\textbf{Output}: Optimal parameters $\bar{\lambda}_1$ and $\bar{\lambda}_2$.}
    \Statex {\textbf{Initialise}: $\lambda_1^0$, $\lambda_2^0$, $n=1$.}
    \While{not converging}
    \Statex {SSN algorithm to compute $(u^n, v^n)$ by solving \eqref{opt:var:eq} with parameters $(\lambda_{1}^{(n)},\lambda_{2}^{(n)})$;}
     \Statex {compute adjoint states $(p_1^n, p_2^n)$;}
     \Statex {compute $F'(u^n)$ using \eqref{formula:gradient};}
     \Statex {BFGS update to compute new $(\lambda_1^{n+1},\lambda_2^{n+1})$;}
     \Statex {Armijo line-search with parameter $\eta\ll 1$;}
      \Statex {$n=n+1$;} 
  \EndWhile
  \end{algorithmic}
\end{algorithm}

\begin{figure}[!h]
\centering
\begin{subfigure}[t]{0.24\textwidth}\centering
\includegraphics[height=4cm]{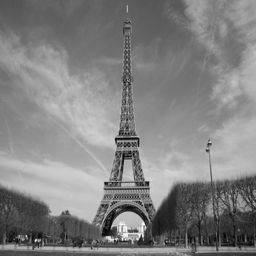}
\caption{$\tilde{u}$}
\label{fig:opt_den_gt}
\end{subfigure}
\begin{subfigure}[t]{0.24\textwidth}\centering
\includegraphics[height=4cm]{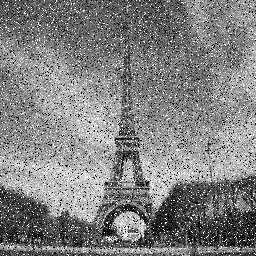}
\subcaption{$\tilde{f}$}
\label{fig:opt_den_f}
\end{subfigure}
\begin{subfigure}[t]{0.24\textwidth}\centering
\includegraphics[height=4cm]{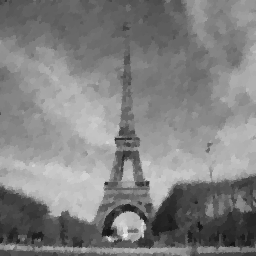}
\caption{$u_{\bar{\lambda}_1,\bar{\lambda}_2}$, $F_{L^2}$ cost.}
\label{fig:opt_den_L2}
\end{subfigure}
\begin{subfigure}[t]{0.24\textwidth}\centering
\includegraphics[height=4cm]{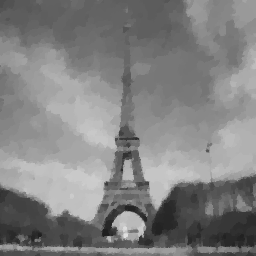}
\caption{$u_{\bar{\lambda}_1,\bar{\lambda}_2}$, $F_{L^1_{\gamma}D}$ cost.}
\label{fig:opt_den_Huber}
\end{subfigure}
\caption{Optimal $\tv$--IC denoising results for initial guess $(\lambda_1^0,\lambda_2^0)=(1,10)$ w.r.t. $F_{L^2}$ and $F_{L^1_{\gamma}D}$ costs \eqref{costfunctL2}--\eqref{costfunctHuber}. Noisy image corrupted with Gaussian noise with variance $\sigma^2=0.01$ and percentage of missing pixels $d=10\%$. \ref{fig:opt_den_f} Noisy image $\tilde{f}$: PSNR=14.40 dB, SSIM=0.17. \ref{fig:opt_den_L2} Optimal $u_{\bar{\lambda}_1,\bar{\lambda}_2}$ w.r.t. $F_{L^2}$ cost: PSNR=28.35 dB, SSIM=0.81. \ref{fig:opt_den_Huber} Optimal $u_{\bar{\lambda}_1,\bar{\lambda}_2}$ w.r.t. $F_{L^1_{\gamma}D}$ cost: PSNR=27.91 dB, SSIM=0.83.}
\label{fig:optimal_denoising1}
\end{figure}

In Figure \ref{fig:optimal_denoising1} we report a numerical experiment confirming the effectiveness of the bilevel optimisation approach on images corrupted by a mixture of Gaussian and Salt \& Pepper noise with Gaussian variance $\sigma^2=0.01$ and percentage of missing pixels $d=10\%$. We report the result obtained with respect to both the $L^2$ cost functional \eqref{costfunctL2} and the Huberised $L^1$ gradient cost \eqref{costfunctHuber}. As observed in \cite{DeLosReyes2017}, we remark that minimising with respect to the $L^2$ cost is indeed equivalent to PSNR optimisation, while the optimisation with respect to the Huberised $L^1$ gradient cost produces better visual results which is similar to optimising the SSIM.

In order to validate numerically the theoretical insights given by the analysis performed in Section \ref{sec:analysis}, we report in the following Figure \ref{fig:cvx_combination1} a plot of the optimal parameters $(\bar{\lambda}_1,\bar{\lambda}_2)$ computed by solving the bilevel system \eqref{eq:optimality_system} with $L^2$ cost \eqref{costfunctL2} and in correspondence of a training pair $(\tilde{u},\tilde{f}_\theta)$ where the noisy image is corrupted by a mixture of Gaussian and Salt \& Pepper noise with varying intensity. Namely, for $\theta\in[0,1]$, we corrupt the training  image $\tilde{u}$ in \ref{fig:cvx_comb_gt} with Gaussian noise with distribution $\mathcal{N}(0,\theta\sigma^2), \sigma^2=0.005$ and Salt \& Pepper noise with a percentage of corrupted pixels equal to $d=(1-\theta)10\%$. Consequently, the noisy image $f_\theta$ is corrupted by pure impulsive noise for $\theta=0$  (Figure \ref{fig:cvx_comb_theta0}), by pure Gaussian noise for $\theta=1$ (Figure \ref{fig:cvx_comb_theta1}) and by a mixture of the two for $\theta\in(0,1)$ (Figures \ref{fig:cvx_comb_theta025}-\ref{fig:cvx_comb_theta075}). As suggested by the theory (see, in particular, Proposition \ref{lbl:Gamma_convergence} and Corollary \ref{lbl:u_l2_to_inf}), we observe that when $\theta=0$ the bilevel strategy selects a large optimal parameter $\bar{\lambda}_2$, enforcing a $\tv$--$L^1$ denoising model which is well known to be optimal for this type of noise (see, e.g., \cite{duvalL1, nikolova04}). Furthermore, by Corollary \ref{lbl:u_l2_to_inf} we also have that the lower-level solution of the mixed noise problem approximates (in the sense of $\Gamma$-convergence) a solution of the corresponding $\tv$--$L^1$ model. It is then interesting to notice that  in the case $\theta=1$ the bilevel optimisation strategy does not enforce a $\tv$--$L^2$ model (i.e., a large $\bar{\lambda}_1$), but rather a combination of $\tv$--$L^1$ and $\tv$--$L^2$.
This might be an indication that in practice the $\tv$--$L^{1}$ model works well also in the case of pure Gaussian noise removal. 
%This is in fact not surprising, since in previous works (see, e.g., \textcolor{red}{REF}) people have noticed how the $\tv$--$L^2$ may be not indeed optimal for Gaussian noise removal and how the use of a $L^1$ discrepancy may produce better results.
This is reflected in our experiment where the estimated optimal data model turns out to be indeed a combination of the two  discrepancies. Note that the two parameters $\lambda_1$ and $\lambda_2$ scale differently, with $\bar{\lambda}_1$ only slightly varying across the different simulations.

We report in Table \ref{table:PSNR} the numerical values of $(\bar{\lambda}_1, \bar{\lambda}_2)$ and the corresponding PSNR values of the noisy images $\tilde{f}_\theta$ and of the optimal reconstructions $u_{\bar{\lambda}_1, \bar{\lambda}_2}$.

%\begin{figure}[!h]
%\centering
%\includegraphics[height=7cm]{berlin256x256_cvx_combination}
%\caption{Optimal parameters $(\bar{\lambda}_1, \bar{\lambda}_2)$ computed by solving the bilevel system \eqref{eq:optimality_system} with $L^2$ cost and training pair $(\tilde{u},\tilde{f}_\theta)$ depending on a parameter $\theta\in[0,1]$ which controls the amount of Gaussian and Salt \& Pepper noise in the data, see Figure \ref{fig:training_pair_cvx_cmb}. For $\theta=0$ only Salt \& Pepper noise is present and the model correctly computes optimal parameters enforcing $\tv$-$L^1$ denoising model. For $0<\theta<1$ the data discrepancies are weighted by smaller parameters $\bar{\lambda}_1$ and $\bar{\lambda}_2$. For each $\theta$ the bilevel Algorithm \ref{alg:bilevel} is initialised with $(\lambda_1^0,\lambda_2^0)=(1,1)$.}
%\label{fig:cvx_combination1}
%\end{figure}

\begin{figure}%{0.9\textwidth}
\centering
\resizebox{0.8\textwidth}{!}{
% This file was created by matlab2tikz.
%
%The latest updates can be retrieved from
%  http://www.mathworks.com/matlabcentral/fileexchange/22022-matlab2tikz-matlab2tikz
%where you can also make suggestions and rate matlab2tikz.
%
\begin{tikzpicture}

\begin{axis}[%
width=10.15in,
height=4.879in,
at={(1.703in,0.658in)},
scale only axis,
xmin=0,
xmax=1,
xtick={0,0.25,0.5,0.75,1},
xlabel style={font=\color{white!15!black}},
xlabel={\LARGE{$\theta$}},
ymin=0,
ymax=140,
ytick={0,20,40,60,80,100,120,140},
axis background/.style={fill=white},
title style={font=\bfseries},
title={\Large{Optimal parameters for weighted Gaussian and Salt \& Pepper mixture}},
legend style={at={(0.796,0.661)}, anchor=south west, legend cell align=left, align=left, draw=white!15!black, row sep=0.2cm, column sep=0.1cm, inner sep=6pt, only marks}
]

\addplot [color=blue,  dash pattern={on 8pt off 6pt } ,  line width=2.0pt, mark size=5.0pt, mark=asterisk, mark options={solid, blue}]
  table[row sep=crcr]{%
0	1.95\\
0.25	2.14548394769629\\
0.5	2.31134002231474\\
0.75	2.50862636579191\\
1	2.47323109541657\\
};

\addlegendentry{\Large{$\bar{\lambda}_1$}}

\addplot [color=red, dash pattern={on 8pt off 6pt } , line width=2.0pt, mark size=5.0pt, mark=asterisk, mark options={solid, red}]
  table[row sep=crcr]{%
0	123.39\\
0.25	61.0386095488387\\
0.5	39.9384458566552\\
0.75	45.8965093493385\\
1	57.3484914154739\\
};
\addlegendentry{\Large{$\bar{\lambda}_2$}}

% \addplot [color=black, draw=none, mark=asterisk,  mark options={solid, blue,mark size=4pt}, line width=0.6mm, forget plot]
%  table[]{%
%0.77	99.5
%};
%\addplot [color=black, draw=none, mark=asterisk,  mark options={solid, red,mark size=4pt}, line width=0.6mm, forget plot]
%  table[]{%
%0.77	89.5
%}; 
%
%  \node at (axis cs:0.8,100) [ black ] {\Large{$\bar{\lambda}_1$}};
% \node at (axis cs:0.8,90) [ black ] {\Large{$\bar{\lambda}_2$}};

\end{axis}
\end{tikzpicture}%
}
\caption{Optimal parameters $(\bar{\lambda}_1, \bar{\lambda}_2)$ computed by solving the bilevel system \eqref{eq:optimality_system} with $L^2$ cost and training pair $(\tilde{u},\tilde{f}_\theta)$ depending on a parameter $\theta\in[0,1]$ which controls the amount of Gaussian and Salt \& Pepper noise in the data, see Figure \ref{fig:training_pair_cvx_cmb}. For $\theta=0$ only Salt \& Pepper noise is present and the model  computes optimal parameters enforcing $\tv$-$L^1$ denoising model. For $0<\theta<1$ the data discrepancies are weighted by smaller parameters $\bar{\lambda}_1$ and $\bar{\lambda}_2$. For each $\theta$ the bilevel Algorithm \ref{alg:bilevel} is initialised with $(\lambda_1^0,\lambda_2^0)=(1,1)$.}
\label{fig:cvx_combination1}
\end{figure}
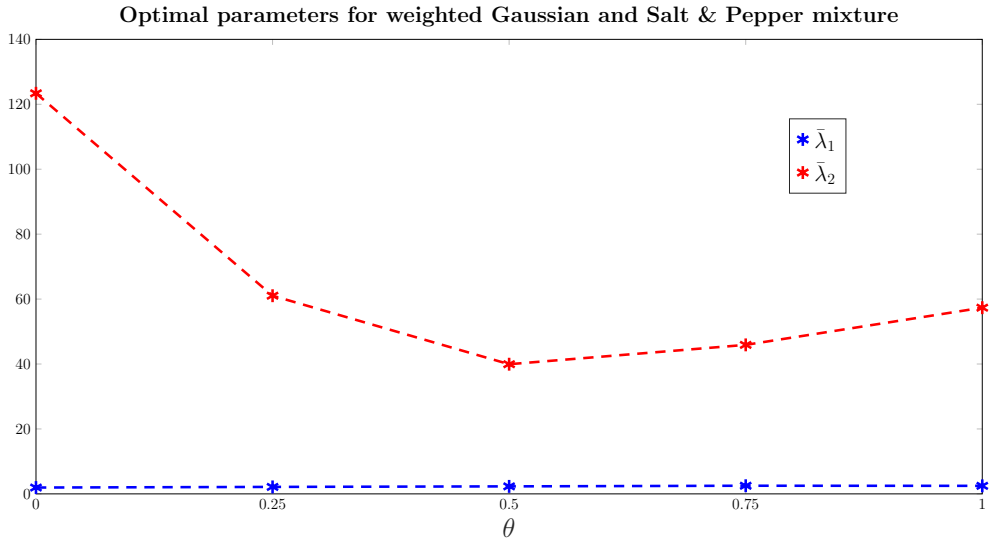

\begin{figure}[!h]
\centering
\begin{subfigure}[t]{0.31\textwidth}\centering
\includegraphics[height=5cm]{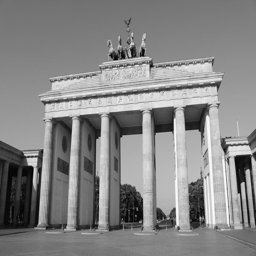}
\caption{$\tilde{u}$}
\label{fig:cvx_comb_gt}
\end{subfigure}
\begin{subfigure}[t]{0.31\textwidth}\centering
\includegraphics[height=5cm]{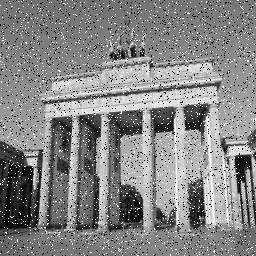}
\subcaption{$\tilde{f}_\theta$, $\theta=0$ (S \& P noise).}
\label{fig:cvx_comb_theta0}
\end{subfigure}
\begin{subfigure}[t]{0.31\textwidth}\centering
\includegraphics[height=5cm]{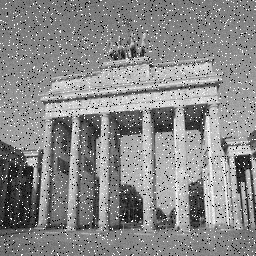}
\caption{$\tilde{f}_\theta$, $\theta=0.25$.}
\label{fig:cvx_comb_theta025}
\end{subfigure}\\
\begin{subfigure}[t]{0.31\textwidth}\centering
\includegraphics[height=5cm]{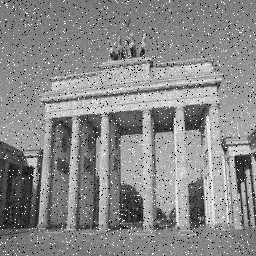}
\caption{$\tilde{f}_\theta$, $\theta=0.5$.}
\label{fig:cvx_comb_theta05}
\end{subfigure}
\begin{subfigure}[t]{0.31\textwidth}\centering
\includegraphics[height=5cm]{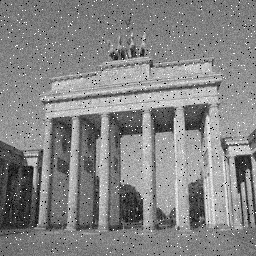}
\caption{$\tilde{f}_\theta$, $\theta=0.75$.}
\label{fig:cvx_comb_theta075}
\end{subfigure}
\begin{subfigure}[t]{0.31\textwidth}\centering
\includegraphics[height=5cm]{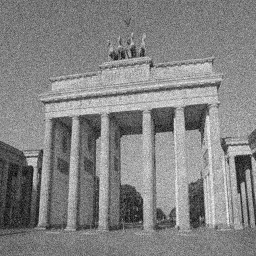}
\caption{$\tilde{f}_\theta$, $\theta=1$ (Gaussian noise).}
\label{fig:cvx_comb_theta1}
\end{subfigure}
\caption{Training images used for the computation of the optimal parameters $(\bar{\lambda}_1, \bar{\lambda}_2)$ whose values are reported in Figure \ref{fig:cvx_combination1}; The noise-free image $\tilde{u}$ is corrupted with a mixture of Gaussian and Salt \& Pepper noise of different intensities, ranging from pure impulsive noise when $\theta=0$ to pure Gaussian when $\theta=1$.}
\label{fig:training_pair_cvx_cmb}
\end{figure}

%\begin{center}
\begin {table}[!h]
\centering
\begin{tabular}{|c||c|c|c|c|}
\hline 
\rule{0pt}{13pt}
 $\theta$ & PSNR $\tilde{f}_\theta$ & $\bar{\lambda}_1$ &$\bar{\lambda}_2$ & PSNR $u_{\bar{\lambda}_1,\bar{\lambda}_2}$ \\
 \hline
 0  &  15.41 dB & 1.95 & 123.39 & 29.71 dB\\
0.25&   16.19 dB & 2.15  & 61.04 & 24.66 dB\\
0.5 & 17.49 dB & 2.31 &  39.94 & 25.35 dB\\
0.75& 19.28 dB & 2.51 &  45.90 & 24.01 dB\\
1 &   22.74 dB & 2.47 & 57.35 & 24.25 dB\\
 \hline
\end{tabular}
\caption{}
 \label{table:PSNR}
 \end{table}
%\end{center}

\section{Conclusions and outlook}  \label{sec:conclusion}

In this paper we have presented a fine analysis of the $\tv$--IC denoising model originally proposed in \cite{calatroni_mixed} for mixed Gaussian and Salt \& Pepper noise removal. Our study started with a characterisation of the IC data-discrepancy as a Huber-regularised $L^1$ discrepancy. We then studied in detail the asymptotic behaviour of the solutions of the  $\tv$--IC model using $\Gamma$-convergence arguments, showing that the solutions of the single $\tv$--$L^1$ and $\tv$--$L^2$ denoising models can be retrieved in the limit as the parameters tend to infinity. We gained more insights on the model by calculating some exact solutions for simple one-dimensional data functions.  Using these theoretical results we then formulated and rigorously analysed a bilevel optimisation approach in function spaces for the estimation of the optimal parameters of the model. 
%We proved existence of a the non-smooth bilevel problem in the case where box constraints are imposed on the parameters.
With the use of a counterexample motivated by our theoretical work, we showed that box constraints on the parameters are indeed necessary to obtain existence.

In the spirit of \cite{bilevellearning,DeLosReyes2017}, after a suitable regularisation of the non-smooth problem, we then proved the differentiability of the solution map via standard arguments as well as the existence of the adjoint states by which we can derive the corresponding optimality system in a compact form.  Thanks to the handy characterisation of the gradient of the bilevel cost functional in terms of its adjoint states, efficient numerical schemes can be easily implemented. In particular, in Section \ref{sec:numres} we considered the second-order BFGS Algorithm \ref{alg:bilevel} to numerically compute the optimal parameters of the $\tv$--IC model for noise mixtures with different noise levels. The numerical results show good agreement with the analytical study of the first sections, making this strategy appealing for blind image denoising applications, where the intensity of each noise component is unknown. 
%We used these numerics to verify in practice the fact that the use $L^{1}$ fidelity is also beneficial in the case of Gaussian noise as well.

Further research could address the validation of the $\tv$--IC model over a set of images in order to estimate the preferred noise model with respect to the denoising application at hand (i.e., the noise intensity, the structure of the image etc.), similarly as done in \cite{DeLosReyes2017}. The IC discrepancy could further be combined with higher-order regularisers (such as $\tgv$ \cite{TGV}) for more visually pleasing reconstructions. Finally, a similar study could be done for more general noise mixtures such as Gaussian \& Poisson noise, for which the IC discrepancy has been shown in \cite{calatroni_mixed} to be a statistically consistent model. 
%Finally, spatially varying parameters can also be considered.

\ack{The authors thank Prof. Juan Carlos De Los Reyes for his useful comments on Section \ref{sec:optimisation}. LC acknowledges the joint ANR/FWF Project ``Efficient Algorithms for Nonsmooth Optimization in Imaging" (EANOI) FWF n. I1148 / ANR-12-IS01-0003,  the Fondation Mathématique Jacques Hadamard (FMJH) and the RISE EU project NoMADS. KP acknowledges the support of the Einstein Foundation Berlin within the ECMath project CH12.  Both authors would like to thank the Isaac Newton Institute for Mathematical Sciences for support and hospitality during the program “Variational Methods and Effective Algorithms for Imaging and Vision” when work on this paper was undertaken. This work was supported by EPSRC grant number EP/K032208/1.}

\section*{References}

\small{
\bibliographystyle{spmpsci}
\bibliography{kostasbib}

\begin{thebibliography}{10}
\providecommand{\url}[1]{{#1}}
\providecommand{\urlprefix}{URL }
\expandafter\ifx\csname urlstyle\endcsname\relax
  \providecommand{\doi}[1]{DOI~\discretionary{}{}{}#1}\else
  \providecommand{\doi}{DOI~\discretionary{}{}{}\begingroup
  \urlstyle{rm}\Url}\fi

\bibitem{Allard1}
Allard, W.: {Total Variation Regularization for Image Denoising, {I}.
  {G}eometric Theory}.
\newblock SIAM Journal on Mathematical Analysis \textbf{39}(4), 1150--1190
  (2008).
\newblock \url{http://dx.doi.org/10.1137/060662617}

\bibitem{Allard2}
Allard, W.: {Total Variation Regularization for Image Denoising, {II}.
  {E}xamples}.
\newblock SIAM Journal on Imaging Sciences \textbf{1}(4), 400--417 (2008).
\newblock \url{http://dx.doi.org/10.1137/070698749}

\bibitem{Allard3}
Allard, W.: {Total Variation Regularization for Image Denoising, {III}.
  {E}xamples.}
\newblock SIAM Journal on Imaging Sciences \textbf{2}(2), 532--568 (2009).
\newblock \url{http://dx.doi.org/10.1137/070711128}

\bibitem{AmbrosioBV}
Ambrosio, L., Fusco, N., Pallara, D.: {Functions of Bounded Variation and Free
  Discontinuity Problems}.
\newblock Oxford University Press, USA (2000)

\bibitem{Baus2014}
Baus, F., Nikolova, M., Steidl, G.: Fully smoothed $\ell_{1}$-${TV}$ models:
  Bounds for the minimizers and parameter choice.
\newblock Journal of Mathematical Imaging and Vision \textbf{48}(2), 295--307
  (2014).
\newblock \url{https://doi.org/10.1007/s10851-013-0420-0}

\bibitem{bauschkecombettes}
Bauschke, H., Combettes, P.: {Convex analysis and monotone operator theory in
  {H}ilbert spaces}.
\newblock CMS Books in Mathematics/Ouvrages de Math{\'e}matiques de la SMC.
  Springer, New York (2011).
\newblock \url{https://doi.org/10.1007/978-3-319-48311-5}

\bibitem{benning2011error}
Benning, M., Burger, M.: Error estimates for general fidelities.
\newblock Electronic Transactions on Numerical Analysis \textbf{38}, 44--68
  (2011)

\bibitem{Benvenuto2008}
Benvenuto, F., La~Camera, A., Theys, C., Ferrari, A., Lant{{\'e}}ri, H.,
  Bertero, M.: The study of an iterative method for the reconstruction of
  images corrupted by {P}oisson and {G}aussian noise.
\newblock Inverse Problems \textbf{24}(3), 035,016 (2008).
\newblock \url{http://stacks.iop.org/0266-5611/24/i=3/a=035016}

\bibitem{TGV}
Bredies, K., Kunisch, K., Pock, T.: {Total Generalized Variation}.
\newblock SIAM Journal on Imaging Sciences \textbf{3}(3), 492--526 (2010).
\newblock \url{http://dx.doi.org/10.1137/090769521}

\bibitem{BrediesL1}
Bredies, K., Kunisch, K., Valkonen, T.: Properties of {L}$^1$-{TGV}$^{\,2}$ :
  The one-dimensional case.
\newblock Journal of Mathematical Analysis and Applications \textbf{398}(1),
  438 -- 454 (2013).
\newblock \url{http://dx.doi.org/10.1016/j.jmaa.2012.08.053}

\bibitem{journal_tvlp}
Burger, M., Papafitsoros, K., Papoutsellis, E., Sch\"onlieb, C.B.: Infimal
  convolution regularisation functionals of {BV} and $\mathrm{L}^{p}$ spaces.
  {P}art {I}: The finite $p$ case.
\newblock Journal of Mathematical Imaging and Vision \textbf{55}(3), 343--369
  (2016).
\newblock \url{http://dx.doi.org/10.1007/s10851-015-0624-6}

\bibitem{impulsegauss2008}
Cai, J., Chan, R., Nikolova, M.: Two-phase approach for deblurring images
  corrupted by impulse plus {G}aussian noise.
\newblock Inverse Problems and Imaging \textbf{2}(2), 187--204 (2008).
\newblock \url{http://dx.doi.org/10.3934/ipi.2008.2.187}

\bibitem{bilevellearning}
Calatroni, L., Chung, C., De~Los~Reyes, J.C., Sch\"onlieb, C.B., Valkonen, T.:
  Bilevel approaches for learning of variational imaging models.
\newblock In: RADON book Series on Computational and Applied Mathematics, vol.
  18. Berlin, Boston: De Gruyter (2017).
\newblock \url{https://www.degruyter.com/view/product/458544}

\bibitem{calatroni_mixed}
Calatroni, L., De~Los~Reyes, J., Sch\"onlieb, C.B.: Infimal convolution of data
  discrepancies for mixed noise removal.
\newblock SIAM Journal on Imaging Sciences \textbf{10}(3), 1196--1233 (2017).
\newblock \url{http://dx.doi.org/10.1137/16M1101684}

\bibitem{calatroni2014dynamic}
Calatroni, L., De~Los~Reyes, J.C., Sch\"onlieb, C.B.: Dynamic sampling schemes
  for optimal noise learning under multiple nonsmooth constraints.
\newblock In: C.~P\"otzsche, C.~Heuberger, B.~Kaltenbacher, F.~Rendl (eds.)
  System Modeling and Optimization, \emph{IFIP Advances in Information and
  Communication Technology}, vol. 443, pp. 85--95. Springer Berlin Heidelberg
  (2014).
\newblock \url{http://dx.doi.org/10.1007/978-3-662-45504-3\_8}

\bibitem{caselles2007discontinuity}
Caselles, V., Chambolle, A., Novaga, M.: The discontinuity set of solutions of
  the {TV} denoising problem and some extensions.
\newblock Multiscale Modeling \& Simulation \textbf{6}(3), 879--894 (2007).
\newblock \url{http://dx.doi.org/10.1137/070683003}

\bibitem{poon_TV_geometric}
Chambolle, A., Duval, V., Peyr\'e, G., Poon, C.: Geometric properties of
  solutions to the total variation denoising problem.
\newblock Inverse Problems \textbf{33}(1), 015,002 (2017).
\newblock \url{http://stacks.iop.org/0266-5611/33/i=1/a=015002}

\bibitem{chanL1}
Chan, T., Esedoglu, S.: {Aspects of Total Variation regularized ${L}^1$
  function approximation}.
\newblock SIAM Journal on Applied Mathematics pp. 1817--1837 (2005).
\newblock \url{http://dx.doi.org/10.1137/040604297}

\bibitem{Zwicknagl}
Choksi, R., Fonseca, I., Zwicknagl, B.: A few remarks on variational models for
  denoising.
\newblock Communications in Mathematical Sciences \textbf{12}(5), 843--857
  (2014).
\newblock \url{http://dx.doi.org/10.4310/CMS.2014.v12.n5.a3}

\bibitem{cristoferi}
Cristoferi, R.: Exact solutions for the denoising problem of piecewise constant
  images in dimension one.
\newblock arXiv preprint arXiv:1612.05508v2  (2017).
\newblock \url{https://arxiv.org/abs/1612.05508v2}

\bibitem{dalmasogamma}
Dal~Maso, G.: Introduction to {$\Gamma$}-convergence.
\newblock Birkh\"auser (1993)

\bibitem{Dauge1992}
Dauge, M.: Neumann and mixed problems on curvilinear polyhedra.
\newblock Integral Equations and Operator Theory \textbf{15}(2), 227--261
  (1992).
\newblock \url{https://doi.org/10.1007/BF01204238}

\bibitem{JCbook}
De~Los~Reyes, J.C.: Numerical {PDE}-constrained optimization.
\newblock SpringerBriefs in Optimization. Springer (2015).
\newblock \url{https://dx.doi.org/10.1007/978-3-319-13395-9}

\bibitem{noiselearning}
De~Los~Reyes, J.C., Sch{{\"o}}nlieb, C.B.: Image denoising: learning the noise
  model via nonsmooth {PDE}-constrained optimization.
\newblock Inverse Problems and Imaging \textbf{7}(4), 1183--1214 (2013).
\newblock \url{http://dx.doi.org/10.3934/ipi.2013.7.1183}

\bibitem{interiorpaper2015}
De~Los~Reyes, J.C., Sch{{\"o}}nlieb, C.B., Valkonen, T.: The structure of
  optimal parameters for image restoration problems.
\newblock Journal of Mathematical Analysis and Applications \textbf{434},
  464--500 (2016).
\newblock \url{https://doi.org/10.1016/j.jmaa.2015.09.023}

\bibitem{DeLosReyes2017}
De~Los~Reyes, J.C., Sch\"{o}nlieb, C.B., Valkonen, T.: Bilevel parameter
  learning for higher-order {T}otal {V}ariation regularisation models.
\newblock Journal of Mathematical Imaging and Vision \textbf{57}(1), 1--25
  (2017).
\newblock \url{https://doi.org/10.1007/s10851-016-0662-8}

\bibitem{Dong2012}
Dong, B., Ji, H., Li, J., Shen, Z., Xu, Y.: Wavelet-based blind image
  inpainting.
\newblock Applied and Computational Harmonic Analysis \textbf{32}(2), 268--279
  (2012).
\newblock \url{http://dx.doi.org/10.1016/j.acha.2011.06.001}

\bibitem{duvalL1}
Duval, V., Aujol, J., Gousseau, Y.: The {TV}-${L1}$ model: a geometric point of
  view.
\newblock SIAM Journal on Multiscale Modeling \& Simulation \textbf{8}(1),
  154--189 (2009).
\newblock \url{http://dx.doi.org/10.1137/090757083}

\bibitem{Groger1989}
Gröger, K.: A {$W$}$^{1,p}$-estimate for solutions to mixed boundary value
  problems for second order elliptic differential equations.
\newblock Mathematische Annalen \textbf{283}(4), 679--688 (1989).
\newblock \url{http://eudml.org/doc/164533}

\bibitem{stadler}
Hinterm\"uler, M., Stadler, G.: {An Infeasible Primal-Dual Algorithm for Total
  Bounded Variation--Based Inf-Convolution-Type Image Restoration}.
\newblock SIAM Journal on Scientific Computing \textbf{28}(1), 1--23 (2006).
\newblock \url{http://dx.doi.org/10.1137/040613263}

\bibitem{structuralTV}
Hinterm\"uller, M., Holler, M., Papafitsoros, K.: A function space framework
  for structural total variation regularization with applications in inverse
  problems.
\newblock Inverse Problems \textbf{34}(6), 064,002 (2018).
\newblock \url{http://stacks.iop.org/0266-5611/34/i=6/a=064002}

\bibitem{langer}
Hinterm\"uller, M., Langer, A.: Subspace correction methods for a class of
  nonsmooth and nonadditive convex variational problems with mixed
  {L}$^1$/{L}$^2$ data-fidelity in image processing.
\newblock SIAM Journal on Imaging Sciences \textbf{6}(4), 2134--2173 (2013).
\newblock \url{http://dx.doi.org/10.1137/120894130}

\bibitem{HintL1}
Hinterm\"uller, M., Monserrat Rincon-Camacho, M.: Expected absolute value
  estimators for a spatially adapted regularization parameter choice rule in
  ${L}^{1}$-{TV}-based image restoration.
\newblock Inverse Problems \textbf{26}(8), 085,005 (2010).
\newblock \url{http://stacks.iop.org/0266-5611/26/i=8/a=085005}

\bibitem{hintermuellerPartI}
Hinterm\"{u}ller, M., Rautenberg, C.N.: Optimal selection of the regularization
  function in a weighted total variation model. part {I}: Modelling and theory.
\newblock Journal of Mathematical Imaging and Vision \textbf{59}(3), 498--514
  (2017).
\newblock \url{https://doi.org/10.1007/s10851-017-0744-2}

\bibitem{hintermuellerPartII}
Hinterm\"{u}ller, M., Rautenberg, C.N., Wu, T., Langer, A.: Optimal selection
  of the regularization function in a weighted total variation model. part
  {II}: Algorithm, its analysis and numerical tests.
\newblock Journal of Mathematical Imaging and Vision \textbf{59}(3), 515--533
  (2017).
\newblock \url{https://doi.org/10.1007/s10851-017-0736-2}

\bibitem{Jalalzai2015jmiv}
Jalalzai, K.: {Some Remarks on the Staircasing Phenomenon in Total
  Variation-Based Image Denoising}.
\newblock Journal of Mathematical Imaging and Vision \textbf{54}(2), 256--268
  (2015).
\newblock \url{http://dx.doi.org/10.1007/s10851-015-0600-1}

\bibitem{poissongauss2013}
Jezierska, A., Chouzenoux, E., Pesquet, J.C., Talbot, H.: {A Convex Approach
  for Image Restoration with Exact Poisson-Gaussian Likelihood}.
\newblock SIAM Journal on Imaging Science \textbf{62}(1), 17--30 (2015).
\newblock \url{http://dx.doi.org/10.1137/15M1014395}

\bibitem{pockSVM}
Klatzer, T., Pock, T.: Continuous hyper-parameter learning for support vector
  machines.
\newblock In: 20$^{\mathrm{th}}$ Computer Vision Winter Workshop Paul Wohlhart,
  Vincent Lepetit (eds.) Seggau, Austria, February 9-11, 2015 (2015).
\newblock \url{http://dx.doi.org/10.3217/978-3-85125-388-7}

\bibitem{pockbilevel}
Kunisch, K., Pock, T.: A bilevel optimization approach for parameter learning
  in variational models.
\newblock SIAM Journal on Imaging Sciences \textbf{6}(2), 938--983 (2013).
\newblock \url{http://dx.doi.org/10.1137/120882706}

\bibitem{Langer2017}
Langer, A.: Automated parameter selection for total variation minimization in
  image restoration.
\newblock Journal of Mathematical Imaging and Vision \textbf{57}(2), 239--268
  (2017).
\newblock \url{http://dx.doi.org/10.1007/s10851-016-0676-2}

\bibitem{langerl1l2}
Langer, A.: Automated parameter selection in the {${L}^{1} \mbox{-}
  {L}^{2}$}-{TV} model for removing gaussian plus impulse noise.
\newblock Inverse Problems \textbf{33}(7), 074,002 (2017).
\newblock \url{http://stacks.iop.org/0266-5611/33/i=7/a=074002}

\bibitem{lanza2013}
Lanza, A., Morigi, S., Sgallari, F., Wen, Y.W.: Image restoration with
  {P}oisson-{G}aussian mixed noise.
\newblock Computer Methods in Biomechanics and Biomedical Engineering: Imaging
  \& Visualization \textbf{2}, 12--24 (2014).
\newblock \url{http://dx.doi.org/10.1080/21681163.2013.811039}

\bibitem{meyer2001oscillating}
Meyer, Y.: Oscillating patterns in image processing and nonlinear evolution
  equations: the fifteenth {D}ean {J}acqueline {B. L}ewis memorial lectures,
  vol.~22.
\newblock American Mathematical Society (2001)

\bibitem{nikolova04}
Nikolova, M.: A variational approach to remove outliers and impulse noise.
\newblock Journal of Mathematical Imaging and Vision \textbf{20}(1), 99--120
  (2004).
\newblock \url{http://dx.doi.org/10.1023/B:JMIV.0000011326.88682.e5}

\bibitem{histo}
Nikolova, M., Wen, Y., Chan, R.: Exact histogram specification for digital
  images using a variational approach.
\newblock Journal of Mathematical Imaging and Vision \textbf{46}(3), 309--325
  (2013).
\newblock \url{http://dx.doi.org/10.1007/s10851-012-0401-8}

\bibitem{pockBilevelNonsmooth}
Ochs, P., Ranftl, R., Brox, T., Pock, T.: Bilevel optimization with nonsmooth
  lower level problems.
\newblock In: J.F. Aujol, M.~Nikolova, N.~Papadakis (eds.) Scale Space and
  Variational Methods in Computer Vision, pp. 654--665. Springer International
  Publishing (2015).
\newblock \url{https://doi.org/10.1007/978-3-319-18461-6_52}

\bibitem{papafitsorosphd}
Papafitsoros, K.: Novel higher order regularisation methods for image
  reconstruction.
\newblock Ph.D. thesis, University of Cambridge (2014).
\newblock \url{https://www.repository.cam.ac.uk/handle/1810/246692}

\bibitem{Papafitsoros_Bredies}
Papafitsoros, K., Bredies, K.: A study of the one dimensional total generalised
  variation regularisation problem.
\newblock Inverse Problems and Imaging \textbf{9}(2), 511--550 (2015).
\newblock {\url{http://dx.doi.org/10.3934/ipi.2015.9.511}}

\bibitem{tgv_asymptotic}
Papafitsoros, K., Valkonen, T.: Asymptotic behaviour of total generalised
  variation.
\newblock In: J.F. Aujol, M.~Nikolova, N.~Papadakis (eds.) Scale Space and
  Variational Methods in Computer Vision: 5th International Conference, SSVM
  2015, Proceedings, pp. 702--714. Springer International Publishing (2015).
\newblock \url{http://dx.doi.org/10.1007/978-3-319-18461-6_56}

\bibitem{poschl2013exact}
P\"oschl, C., Scherzer, O.: Exact solutions of one-dimensional total
  generalized variation.
\newblock Communications in Mathematical Sciences \textbf{13}(1), 171--202
  (2015).
\newblock \url{http://dx.doi.org/10.4310/CMS.2015.v13.n1.a9}

\bibitem{ring2000structural}
Ring, W.: Structural properties of solutions to total variation regularization
  problems.
\newblock ESAIM: Mathematical Modelling and Numerical Analysis \textbf{34}(4),
  799--810 (2000).
\newblock \url{http://dx.doi.org/10.1051/m2an:2000104}

\bibitem{rudin1992nonlinear}
Rudin, L., Osher, S., Fatemi, E.: Nonlinear total variation based noise removal
  algorithms.
\newblock Physica D: Nonlinear Phenomena \textbf{60}(1-4), 259--268 (1992).
\newblock \url{http://dx.doi.org/10.1016/0167-2789(92)90242-F}

\bibitem{valkonen2014jump1}
Valkonen, T.: The jump set under geometric regularization. {P}art 1: Basic
  technique and first-order denoising.
\newblock SIAM Journal on Mathematical Analysis \textbf{47}(4), 2587--2629
  (2015).
\newblock \url{https://doi.org/10.1137/140976248}

\bibitem{valkonen2014jump2}
Valkonen, T.: The jump set under geometric regularisation. {P}art 2:
  Higher-order approaches.
\newblock Journal of Mathematical Analysis and Applications \textbf{453}(2),
  1044--1085 (2017).
\newblock \url{https://doi.org/10.1016/j.jmaa.2017.04.037}

\bibitem{Yan2013}
Yan, M.: Restoration of images corrupted by impulse noise and mixed gaussian
  impulse noise using blind inpainting.
\newblock SIAM Journal on Imaging Sciences \textbf{6}(3), 1227--1245 (2013).
\newblock \url{http://dx.doi.org/10.1137/12087178X}

\bibitem{pockReactionDiffusion}
Yu, W., Heber, S., Pock, T.: Learning Reaction-Diffusion Models for Image
  Inpainting, vol. 9358, pp. 356--367.
\newblock Springer International Publishing AG, Switzerland (2015).
\newblock \url{https://doi.org/10.1007/978-3-319-24947-6_29}

\bibitem{Zowe1979}
Zowe, J., Kurcyusz, S.: Regularity and stability for the mathematical
  programming problem in {B}anach spaces.
\newblock Applied Mathematics and Optimization \textbf{5}(1), 49--62 (1979).
\newblock \url{https://doi.org/10.1007/BF01442543}

\end{thebibliography}
}
\end{document}